\documentclass[11pt]{article}

\usepackage{amssymb,amsbsy,amsthm,amsmath,graphicx,mathrsfs,epsfig,graphics,times}

\input xy

\xyoption{all}

\begin{document}

\title{\textbf{Extensions of probability-preserving systems by measurably-varying homogeneous spaces and applications}}
\author{Tim Austin}
\date{}

\maketitle

\renewcommand{\thefootnote}{}

\footnote{2000 \emph{Mathematics Subject Classification}: Primary
28D15; Secondary 37A30, 37A45.}

\footnote{\emph{Key words and phrases}: Probability-preserving
system, skew product, isometric extension, ergodic decomposition,
Mackey Theory, Furstenberg-Zimmer Theory.}

\renewcommand{\thefootnote}{\arabic{footnote}}
\setcounter{footnote}{0}


\newenvironment{nmath}{\begin{center}\begin{math}}{\end{math}\end{center}}

\newtheorem{thm}{Theorem}[section]
\newtheorem{lem}[thm]{Lemma}
\newtheorem{prop}[thm]{Proposition}
\newtheorem{cor}[thm]{Corollary}
\newtheorem{conj}[thm]{Conjecture}
\newtheorem{dfn}[thm]{Definition}
\newtheorem{prob}[thm]{Problem}
\newtheorem{ques}[thm]{Question}
\theoremstyle{remark}
\newtheorem*{ex}{Example}


\renewcommand{\a}{\alpha}
\newcommand{\A}{\mathcal{A}}
\renewcommand{\b}{\beta}
\newcommand{\B}{\mathcal{B}}
\newcommand{\co}{\mathrm{co}}
\renewcommand{\d}{\mathrm{d}}
\newcommand{\e}{\mathrm{e}}
\newcommand{\eps}{\varepsilon}
\newcommand{\E}{\mathcal{E}}
\newcommand{\F}{\mathcal{F}}
\newcommand{\g}{\gamma}
\newcommand{\G}{\Gamma}
\newcommand{\hcf}{\mathrm{hcf}}
\newcommand{\frH}{\mathfrak{H}}
\newcommand{\im}{\mathrm{i}}
\renewcommand{\k}{\kappa}
\renewcommand{\l}{\lambda}
\renewcommand{\L}{\Lambda}
\newcommand{\frM}{\mathfrak{M}}
\newcommand{\Aut}{\mathrm{Aut}}
\newcommand{\Lat}{\mathrm{Lat}}
\newcommand{\Clos}{\mathrm{Clos}}
\renewcommand{\O}{\Omega}
\renewcommand{\Pr}{\mathrm{Pr}}
\newcommand{\s}{\sigma}
\renewcommand{\P}{\mathcal{P}}
\renewcommand{\S}{\Sigma}
\newcommand{\T}{\mathrm{T}}
\newcommand{\U}{\mathcal{U}}

\newcommand{\bbN}{\mathbb{N}}
\newcommand{\bbR}{\mathbb{R}}
\newcommand{\bbC}{\mathbb{C}}
\newcommand{\bbZ}{\mathbb{Z}}
\newcommand{\bbQ}{\mathbb{Q}}
\newcommand{\bbT}{\mathbb{T}}
\newcommand{\sfE}{\mathsf{E}}
\newcommand{\id}{\mathrm{id}}
\newcommand{\BSp}{\mathsf{BSp}}
\newcommand{\WBSp}{\mathsf{WBSp}}
\newcommand{\Sys}{\mathsf{Sys}}
\newcommand{\Rep}{\mathsf{Rep}}
\newcommand{\CSG}{\mathrm{CSG}}
\newcommand{\SG}{\mathrm{SG}}

\newcommand{\bfX}{\mathbf{X}}
\newcommand{\bfY}{\mathbf{Y}}
\newcommand{\bfZ}{\mathbf{Z}}
\newcommand{\bfW}{\mathbf{W}}
\newcommand{\bfV}{\mathbf{V}}

\newcommand{\bb}[1]{\mathbb{#1}}
\newcommand{\bs}[1]{\boldsymbol{#1}}
\newcommand{\fr}[1]{\mathfrak{#1}}
\renewcommand{\bf}[1]{\mathbf{#1}}
\renewcommand{\rm}[1]{\mathrm{#1}}
\renewcommand{\cal}[1]{\mathcal{#1}}
\newcommand{\uhr}{\!\!\upharpoonright}

\newcommand{\sfZ}{\mathsf{Z}}

\newcommand{\mprod}{\hbox{$\prod$}}
\newcommand{\motimes}{\hbox{$\bigotimes$}}
\newcommand{\mvee}{\hbox{$\bigvee$}}
\newcommand{\mwedge}{\hbox{$\bigwedge$}}
\renewcommand{\t}[1]{\tilde{#1}}

\newcommand{\fin}{\nolinebreak\hspace{\stretch{1}}$\lhd$}

\begin{abstract}
We study a generalized notion of a homogeneous skew-product
extension of a probability-preserving system in which the
homogeneous space fibres are allowed to vary over the ergodic
decomposition of the base. The construction of such extensions rests
on a simple notion of `direct integral' for a `measurable family' of
homogeneous spaces, which has a number of precedents in older
literature.  The main contribution of the present paper is the
systematic development of a formalism for handling such extensions,
including non-ergodic versions of the results of Mackey describing
ergodic components of such extensions~\cite{Mac66}, of the
Furstenberg-Zimmer Structure Theory~\cite{Zim76.1,Zim76.2,Fur77} and
of results of Mentzen~\cite{Men91} describing the structure of
automorphisms of such extensions when they are relatively ergodic.
We then offer applications to two structural results for actions of
several commuting transformations: firstly to describing the
possible joint distributions of three isotropy factors corresponding
to three commuting transformations; and secondly to describing the
characteristic factors for a system of double nonconventional
ergodic averages (see~\cite{Aus--nonconv} and the references listed
there).  Although both applications are modest in themselves, we
hope that they point towards a broader usefulness of this formalism
in ergodic theory.
\end{abstract}

\parskip 0pt

\tableofcontents

\parskip 7pt

\section{Introduction}

This work is concerned with probability-preserving actions $T:\G
\curvearrowright (X,\mu)$ of locally compact second countable
topological groups on standard Borel probability spaces.  We often
denote such an action by $(X,\mu,T)$ if the group is understood.

One of the more versatile constructions by which a more complicated
system may be built from comparatively simple ingredients is the
homogeneous skew-product (see, for example, Examples 2.21 in
Glasner~\cite{Gla03}). From some given $\G$-system $(Y,\nu,S)$, a
compact group $G$ and a closed subgroup $K\leq G$, and a measurable
cocycle $\rho:\G\times Y \to G$ for the action $S$, we form the
system $(Y\times G/K,\nu\otimes m_{G/K},T)$ by setting
\[T^\g(y,gK) := (S^\g y,\rho(\g,y)gK)\quad\quad\hbox{for}\ \g \in \G.\]

A well-developed theory of such systems is available in case the
base system $(Y,\nu,S)$ is ergodic (much of which can be extended to
the setting in which $\nu$ is only quasi-invariant under $S$;
see~\cite{AusLem--Rokhlincocycs}). In addition to providing a wealth
of example systems, such homogeneous skew-products over ergodic base
systems acquire a greater significance through the structure theory
developed by Zimmer in~\cite{Zim76.1,Zim76.2} and Furstenberg
in~\cite{Fur77}. This is concerned with the failure of relative weak
mixing of extensions (see, for example, Definition 9.22 in
Glasner~\cite{Gla03}). Relative weak mixing is a strengthening of
the condition of relative ergodicity which has numerous consequences
for how this extension may be joined to others.  An understanding of
these consequences and of the ways in which relative weak mixing can
fail is crucial to Furstenberg's approach to Szemer\'edi's
Theorem~(\cite{Fur77}; see also the excellent treatment
in~\cite{Fur81}). The core result of Furstenberg and Zimmer is an
inverse theorem according to which an extension of ergodic systems
fails to be relatively weakly mixing if and only if it contains a
nontrivial subextension that can be coordinatized as a homogeneous
skew-product: thus, homogeneous skew-products are identified as
precisely the `obstructions' to relative weak mixing within other
extensions.

However, in many applications in which this ergodicity of the base
system fails, this simple homogeneous skew-product construction is
not quite general enough, and the Furstenberg-Zimmer theory outlined
above is not available without modification.

In this paper we shall extend the definition of homogeneous
skew-product to a more general class of systems by the simple
artifice of allowing the fibre $G_y/K_y$, in addition to the cocycle
$\rho(y)$, to vary as a function of the base point $y \in Y$.  This
leads to a definition of an extension of $(Y,\nu,S)$ given as an
action of the group on a `direct integral' of homogeneous space
fibres over $(Y,\nu)$. It is relatively simple to settle on a
definition of `measurability' for such an assignment of fibres, and
to turn this idea into a rigorous definition.

The study of such measurably-varying groups is certainly not new. It
is already alluded to during the introduction to Section 6.1 of
Guichardet's book~\cite{Gui72} in the context of `measurable current
groups' and their representation theory, motivated in turn by
considerations from algebraic quantum field theory. A number of more
recent works have studied constructions of this nature explicitly.
For example, in~\cite{ConRau07} Conze and Raugy analyze the ergodic
decomposition of various $\s$-finite measures associated to an
extension of a non-singular ergodic base transformation by a locally
compact non-compact group, invoking for their description a
measurably-varying family of subgroups of that fibre group. However,
in their setting the need for a measurably-varying family of groups
is related to the possibly non-smooth structure of the Borel
equivalence relation of conjugacy on the collection of all such
subgroups, an issue which disappears upon restriction to the case of
compact fibre groups, and so the results that they develop are still
rather disconnected from the present paper. Perhaps closest to our
present emphasis is the emergence of measurably-varying subgroups of
a fixed compact group in the analysis of measurably-varying Mackey
groups for certain non-ergodic self-joinings of an ergodic system,
which underlies some known approaches to the study of
non-conventional ergodic averages; see, in particular, Section 3.2
in Meiri~\cite{Mei90}, the proof of Lemma 9.2 in Furstenberg and
Weiss~\cite{FurWei96} and Section 6.8 in Ziegler~\cite{Zie07}.

More generally, a need to extend known machinery for different kinds
of ergodic system to their non-ergodic relatives has been felt in
other areas: consider, for example, Downarowicz' discussion
in~\cite{Dow06} of `assignments' for his study of simplices of
invariant measures for topological systems on zero-dimensional
metric spaces, and the work of Fisher, Witte-Morris and Whyte on
cocycle superrigidity for non-ergodic systems~\cite{FisMorWhy04}.
The careful study of such non-ergodic settings also has many
parallels with the analysis of direct integrals in the
infinite-dimensional representation theory of locally compact groups
or von Neumann and C$^\ast$-algebras (nicely introduced, for
instance, in Arveson~\cite{Arv76}), although we will not explore
this connection further here.

Notwithstanding the diversity of these previous developments, within
structural ergodic theory the treatment of extensions with varying
homogeneous space fibres seems to have stayed largely immersed in
other analyses, such as those cited above.  Although it is
intuitively clear that the fundamental structural results for the
ergodic case of homogeneous skew-product extensions should admit
natural generalizations, it seems that this not yet been carried
out. In fact, after setting up the right definitions we will find
that it is largely routine to extend both the results of
Mackey~(\cite{Mac66}) on the invariant factor of the extended system
and also the Furstenberg-Zimmer Theory to this setting. We lay out
the details of this generalization in the first part of this paper.
More interesting is the extension of the results of
Mentzen~\cite{Men91} on the possible structure of an automorphism of
an isometric extension of ergodic systems: we generalize this by
presenting structure theorems for factors and automorphisms of a
relatively ergodic extension by measurably-varying compact
homogeneous spaces.

Although this generalization as much a matter of care as new ideas,
it pays off by broadening the applicability of the theory of
homogeneous skew-products to settings in which an assumption of base
ergodicity is unavailable. This arises, in particular, when
considering an action of a larger group $T:\G\curvearrowright
(X,\mu)$ restricted to some subgroup $\L \leq \G$.  Although a
routine appeal to the ergodic decomposition can often justify the
assumption that $T$ is ergodic overall, if we disintegrate further
to guarantee that the restricted action $T^{\ \uhr\L}$ is ergodic
then $T^\g$ for $\g \in\G\setminus\L$ need not preserve the
resulting disintegrands of $\mu$.  As a consequence, if we are
concerned with how the $\L$-subaction sits within the whole original
action, we may be forced to retain a system for which this subaction
of $\L$ is not ergodic.

In this paper we offer two closely-related applications meeting this
description. For both cases we specialize to $\G = \bbZ^d$. These
two applications are relatively simple, and are included largely to
illustrate the arguments made possible by the formalism described
above, but they also exemplify much more general questions on which
we suspect these methods will shed light in the future.

Given a $\bbZ^d$-system $\bfX = (X,\mu,T)$ we can consider the
$\s$-subalgebra $\S_X^{T^{\ \uhr\L}}$ of sets left invariant by the
subaction of $T$ corresponding to some subgroup $\G \leq \bbZ^d$. As
is standard in the category of standard Borel spaces, this can be
identified up to $\mu$-negligible sets with the $\s$-algebra
generated by a factor map $\zeta_0^{T^{\ \uhr\G}}:\bfX\to
\bfZ_0^{T^{\ \uhr\G}}$ to some new system on which the subaction of
$\G$ is trivial. Although individually these new systems can still
be quite complicated, a possibly more tractable task is to describe
their possible joint distributions within the original system. If
$\G,\L \leq \bbZ^d$ are two subgroups then it is easy to show that
$\zeta_0^{T^{\ \uhr\G}}$ and $\zeta_0^{T^{\ \uhr\L}}$ are relatively
independent over $\zeta_0^{T^{\ \uhr(\G + \L)}}$, but for three or
more subgroups matters become more complicated. Clearly given three
subgroups $\G_1,\G_2,\G_3 \leq \bbZ^d$ we have that $\zeta_0^{T^{\
\uhr \G_1}}$ and $\zeta_0^{T^{\ \uhr\G_2}}$ both contain
$\zeta_0^{T^{\ \uhr(\G_1 + \G_2)}}$, and similarly for other pairs,
and so a na\"\i ve candidate for a generalization of the above
result could be that the three isotropy factors $\zeta_0^{T^{\
\uhr\G_i}}$ are relatively independent over the smaller triple of
factors $\zeta_0^{T^{\ \uhr(\G_i+\G_j)}}\vee \zeta_0^{T^{\
\uhr(\G_i+\G_k)}}$ (the factor generated by $\zeta_0^{T^{\
\uhr(\G_i+\G_j)}}$ and $\zeta_0^{T^{\ \uhr(\G_i+\G_k)}}$ together)
as $(i,j,k)$ ranges over permutations of $(1,2,3)$. If we denote the
target $\bbZ^d$-system of this joint factor map by $\bfW_i$ (so this
is a joining of $\bfZ_0^{T^{\ \uhr(\G_i + \G_j)}}$ and $\bfZ_0^{T^{\
\uhr(\G_i + \G_k)}}$) and let $\a_i:\bfZ_0^{T^{\ \uhr\G_i}}\to
\bfW_i$ be the factor map defined by $\zeta_0^{T^{\
\uhr(\G_i+\G_j)}}\vee \zeta_0^{T^{\ \uhr(\G_i+\G_k)}} =
\a_i\circ\zeta_0^{T^{\ \uhr\G_i}}$, then these factors are arranged
as in the following commutative diagram:
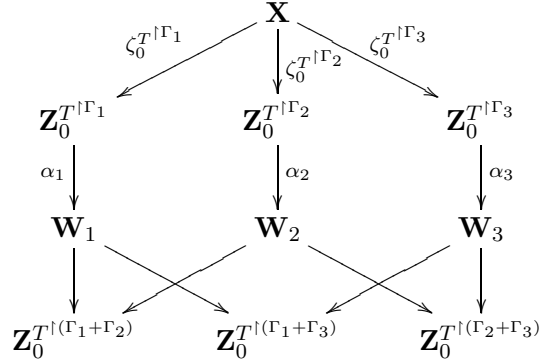
\begin{center}
$\phantom{i}$\xymatrix{
&\bfX\ar[dl]_{\zeta_0^{T^{\ \uhr\G_1}}}\ar[d]^{\zeta_0^{T^{\ \uhr\G_2}}}\ar[dr]^{\zeta_0^{T^{\ \uhr\G_3}}}\\
\bfZ_0^{T^{\ \uhr\G_1}}\ar[d]_{\a_1} & \bfZ_0^{T^{\
\uhr\G_2}}\ar[d]^{\a_2} &
\bfZ_0^{T^{\ \uhr\G_3}}\ar[d]^{\a_3}\\
\bfW_1\ar[d]\ar[dr] & \bfW_2\ar[dl]\ar[dr] & \bfW_3\ar[dl]\ar[d]\\
\bfZ_0^{T^{\ \uhr(\G_1 + \G_2)}} & \bfZ_0^{T^{\ \uhr(\G_1 + \G_3)}}
& \bfZ_0^{T^{\ \uhr(\G_2 + \G_3)}} }
\end{center}

In fact the na\"\i ve conjecture that the factors $\zeta_0^{T^{\
\uhr\G_i}}$ are relatively independent over their further factors
$\a_i$ is false, but `not by very much': we will see that it can
fail only in a very restricted way. In general, the three factors
$\zeta_0^{T^{\ \uhr\G_i}}$ are relatively independent over some
subextensions of these `natural candidate' factors $\a_i$, and these
subextensions can be coordinatized by measurable compact fibre
groups subject to certain further restrictions.

Here we will examine this when $d=3$ and $\G_i$ is the cyclic
subgroup $\bbZ\bf{e}_i$ in the direction of a basis vector
$\bf{e}_i$, but it seems clear that our methods can be extended both
to more general subgroups of Abelian groups and (probably with
considerably more work) to larger numbers of subgroups.

\begin{thm}[Joint distributions of three isotropy
factors]\label{thm:threefoldisotropy} Let $\bfX = (X,\mu,T)$ be a
$\bbZ^3$-system and write $T_i := T^{\bf{e}_i}$ for $i=1,2,3$. Let
$\bfW_i$ be the target of the joint factor map $\a_i :=
\zeta_0^{T_i,T_j}\vee\zeta_0^{T_i,T_k}$, where $\zeta_0^{T_i,T_j} :=
\zeta_0^{T^{\ \uhr(\bbZ\bf{e}_i + \bbZ\bf{e}_j)}}$, let $W_i$ be its
underlying standard Borel space and let $T_j|_{\a_i}$ be the
restriction of $T_j$ to the factor $\a_i$.

Between the single isotropy factors
$\zeta_0^{T_i}:\bfX\to\bfZ_0^{T_i}$ and the smaller factors
$\a_i:\bfX\to \bfW_i$ there are three intermediate factors
$\phi_i\circ\zeta_0^{T_i}:\bfX\to \bfV_i$, where
\[\bfX\stackrel{\zeta_0^{T_i}}{\longrightarrow}\bfZ_0^{T_i}\stackrel{\phi_i}{\longrightarrow}\bfV_i\longrightarrow\bfW_i,\]
such that
\begin{itemize}
\item the triple of factors $\zeta_0^{T_i}$ is relatively
independent over the triple $\phi_i\circ\zeta_0^{T_i}$ under $\mu$;
\item there exist compact metrizable group data $G_{i,\bullet}$ on $W_i$ invariant
under the restriction of the whole action $T$ to the factor space
$W_i$, a cocycle $\tau_{ij}:W_i \to G_{i,\bullet}$ invariant under
the restriction of $T_k$ to $W_i$ and a cocycle $\tau_{ik}:W_i \to
G_{i,\bullet}$ invariant under the restriction of $T_j$ to $W_i$
such that we can coordinatize the extension $\bfV_i\to\bfW_i$ as the
extension of $\bfW_i$ by the measurable compact fibre groups
$G_{i,\bullet}$ with the lifted actions defined by
\[T_j|_{\phi_i\circ\zeta_0^{T_i}}(w_i,g_i) = (T_j|_{\a_i}(w_i),\tau_{ij}(w_i)\cdot g_i)\]
and
\[T_k|_{\phi_i\circ\zeta_0^{T_i}}(w_i,g_i) = (T_k|_{\a_i}(w_i),g_i\cdot \tau_{ik}(w_i)).\]
\end{itemize}
\end{thm}

We will generally denote a coordinatization of the extension $\bfV_i
\to \bfW_i$ as above by the commutative diagram
\begin{center}
$\phantom{i}$\xymatrix{ \bfV_i\ar[dr]\ar@{<->}[rr]^-\cong &&
\bfW_i\ltimes (G_{i,\bullet},m_{G_{i,\bullet}},\tau_{ij},\tau_{ik}^{\rm{op}})\ar[dl]^{\rm{canonical}}\\
& \bfW_i, }
\end{center}
where we use the superscript $^{\rm{op}}$ to denote a cocycle that
acts on fibres by right-multiplication, and we have suppressed
mention of the transformation $T_i$ since by definition its
restriction to $\bfY_i$ is the identity.

Although our final conclusion here yields measurably-varying fibre
groups $G_{i,\bullet}$ that are invariant under the whole action $T$
--- and so would be constant if we had assumed that the overall action
$T$ is ergodic --- the analysis leading to this conclusion will go
via homogeneous space fibres of possibly greater variability, for
the reason described earlier that at first we will only be able to
assume that the fibres are invariant under the subaction $T_i$.

The same is true of our second application.  This is to a special
case of the problem of describing the `minimal characteristic
factors' for the nonconventional ergodic averages
\[\frac{1}{N}\sum_{n=1}^N\prod_{i=1}^df_i\circ T_i^n\]
associated to a $d$-tuple of commuting actions
$T_i:\bbZ\curvearrowright (X,\mu)$ and functions $f_1,f_2,\ldots,f_d
\in L^\infty(\mu)$.  Let us write $\bfX = (X,\mu,T)$ for the
$\bbZ^d$-system given by these one-dimensional actions in the
coordinate directions.

The question of convergence in $L^2(\mu)$ for such averages was
first settled when $d=2$ by Conze and Lesigne in~\cite{ConLes84},
and since then a number of other works have addressed other versions
or relatives of this question~\cite{Zha96,FurWei96,HosKra01},
culminating in Host and Kra's detailed analysis of the case in which
$T_i = T^i$ for some fixed $T$ in~\cite{HosKra05} (see also
Ziegler~\cite{Zie07}) and Tao's recent proof
in~\cite{Tao08(nonconv)} of convergence for arbitrary $d$. We direct
the reader to~\cite{Aus--nonconv} for a more detailed discussion of
this problem and an alternative proof of convergence.

Here we will consider in the case $d=2$ an important part of these
developments: the theory of `characteristic factor-tuples' for such
averages. In our setting, a pair of factors $\xi_i:\bfX\to \bfY_i$
is `characteristic' if
\[\frac{1}{N}\sum_{n = 1}^N(f_1\circ T_1^n)\cdot(f_2\circ T_2^n) - \frac{1}{N}\sum_{n = 1}^N\sfE_\mu(f_1\circ T_1^n\,|\,\xi_1)\cdot\sfE_\mu(f_2\circ T_2^n\,|\,\xi_2)\to 0\]
in $L^2(\mu)$ as $N\to\infty$ for any $f_1,f_2 \in L^\infty(\mu)$.
Clearly given such a pair of factors, the problem of proving
convergence reduces to the case when each $f_i$ is
$\xi_i$-measurable, and this reduction forms an important first step
in many of the known proofs of convergence. Although characteristic
factors are well-understood in some special cases, the more recent
proofs of general convergence
in~\cite{Tao08(nonconv),Tow07,Aus--nonconv} proceed by first heavily
modifying the original system and only then asking after the
characteristic factors (or their finitary analog in Tao's proof
in~\cite{Tao08(nonconv)}), and so our knowledge of the
characteristic factors of the original system remains incomplete
except in some special cases~\cite{Zha96,HosKra05,FraKra05,Zie07}.
More is known in the case $d=2$ from the work of Conze and
Lesigne~\cite{ConLes84}, and in addition the following very precise
description of the characteristic factors when $d=2$ has achieved
folkloric currency since that work appeared. However, a complete
proof seems to be surprisingly subtle, and we shall give such a
proof as our second application of our non-ergodic machinery for
extensions by homogeneous spaces.

\begin{thm}[Characteristic factors for double nonconventional
averages]\label{thm:charfactors} Given a $\bbZ^2$-system $\bfX =
(X,\mu,T_1,T_2)$, let $\bfW_i$ be the target system of the factor
$\a_i := \zeta_0^{T_i}\vee\zeta_0^{T_1T_2^{-1}}$ with underlying
standard Borel space $W_i$. Then $\bfX$ admits a characteristic pair
of factors $\xi_i:\bfX\to\bfY_i$ with underlying standard Borel
spaces $Y_i$ that extend the factors $\bfX\to\bfW_i$ and can be
described as follows: there are
\begin{itemize}
\item a $T$-invariant measurable family of compact fibre groups
$G_\bullet$,
\item a $T_1$-ergodic cocycle $\s:\bfZ_0^{T_1T_2^{-1}}
\to G_\bullet$ that is ergodic for the restricted action of $T_1$,
\item and a pair of cocycles $\tau_i:\bfZ_0^{T_{3-i}}\to G_\bullet$ ergodic for
the restricted action of $T_i$
\end{itemize}
such that we can coordinatize these probability spaces as
\begin{center}
$\phantom{i}$\xymatrix{
(Y_1,(\xi_1)_\#\mu)\ar[dr]\ar@{<->}[rr]^-\cong &&
(W_1,(\a_1)_\#\mu)\ltimes (G_\bullet,m_{G_\bullet})\ar[dl]^{\rm{canonical\ map}}\\
& (W_1,(\a_1)_\#\mu) }
\end{center}
so that the restricted actions are given by
\begin{itemize}
\item[] restriction of $T_1$: $(w_1,g)\mapsto \big(T_1|_{\a_1}(w_1),\s(\zeta_0^{T_1T_2^{-1}}(w_1))\cdot
g\big)$,
\item[] restriction of $T_2$: $(w,g)\mapsto \big(T_2|_{\a_1}(w_1),\s(\zeta_0^{T_1T_2^{-1}}(w_1))\cdot
g\cdot \tau_2(\zeta_0^{T_1}(w_1))\big)$,
\end{itemize}
and similarly
\begin{center}
$\phantom{i}$\xymatrix{
(Y_2,(\xi_2)_\#\mu)\ar[dr]\ar@{<->}[rr]^-\cong &&
(W_2,(\a_2)_\#\mu)\ltimes (G_\bullet,m_{G_\bullet})\ar[dl]^{\rm{canonical\ map}}\\
& (W_2,(\a_2)_\#\mu) }
\end{center}
with
\begin{itemize}
\item[] restriction of $T_1$: $(w_2,g)\mapsto \big(T_1|_{\a_2}(w_2),\s(\zeta_0^{T_1T_2^{-1}}(w_2))\cdot
g\cdot \tau_1(\zeta_0^{T_2}(w_2))\big)$,
\item[] restriction of $T_2$: $(w_2,g)\mapsto \big(T_2|_{\a_2}(w_2),\s(\zeta_0^{T_1T_2^{-1}}(w_2))\cdot g\big)$.
\end{itemize}
\end{thm}

We suspect that our methods should extend to offer at least some
description of characteristic factor-tuples for larger numbers of
commuting transformations, although we also suspect that it will
become rapidly more complicated.

In summary, the body of this paper is organized as follows.

In Section~\ref{sec:background} we recall some definitions and
standard results from group theory, measure theory and ergodic
theory that we will need later in the paper, and in doing so set up
some convenient notation.

Section~\ref{sec:dirint} introduces our definitions of measurable
families of homogeneous space data and their direct integrals.

In Section~\ref{sec:Mackey} we cover quite briskly the main
definitions and results of the non-ergodic Mackey Theory, and then
in Section~\ref{sec:FZ} we treat similarly the non-ergodic version
of the Furstenberg-Zimmer inverse theory.

In Section~\ref{sec:autos} we pursue a slightly less standard
consequence of the Mackey Theory, using it first to describe the
possible factors and groups of automorphisms of an extension by
homogeneous space data, and then translating this into conditions on
an automorphism of a base system that it be liftable to an
automorphism of an extension. This generalizes the classical work of
Mentzen~\cite{Men91} in the case of ergodic systems, and will be
important for the applications of the theory that follow.

In Section~\ref{sec:app} we present our two applications, to the
joint three-fold distributions of isotropy factors and to double
characteristic factors.

Finally, in Section~\ref{sec:furtherques} we discuss some further
possible applications of this machinery.

\textbf{Acknowledgements}\quad My thanks go to Vitaly Bergelson,
John Griesmer, Bernard Host, Keith Kearnes, Bryna Kra, Alexander
Leibman, Terence Tao and Tamar Ziegler for several helpful
discussions and communications, and to the Mathematical Sciences
Research Institute (Berkeley) 2009 program on Ergodic Theory and
Additive Combinatorics, during which large parts of this work were
completed.

\section{Background and notation}\label{sec:background}

\subsection{Measurable functions and probability kernels}

We will work exclusively in the category of standard Borel
probability spaces $(X,\S_X,\mu)$, and so will often suppress
mention of their $\s$-algebras.

Any Borel map $\phi:X\to Y$ specifies a $\s$-subalgebra of $\S_X$ in
the form of $\phi^{-1}(\S_Y)$. Two such maps $\phi:X\to Y$ and
$\psi:X\to Z$ are \textbf{equivalent} if these $\s$-subalgebras of
$\S_X$ that they generate are equal up to $\mu$-negligible sets, in
which case we shall write $\phi \simeq \psi$; this clearly defines
an equivalence relation among Borel maps with domain $X$. As it
standard, in the category of standard Borel spaces equivalence
classes of such Borel maps are in bijective correspondence with
equivalence classes of $\s$-subalgebras under the relation of
equality modulo the $\s$-ideal of $\mu$-negligible sets. A treatment
of these classical issues may be found, for example, in Chapter 2 of
Glasner~\cite{Gla03}.

A measure-preserving Borel map $\pi:(X,\mu) \to (Y,\nu)$
\textbf{contains} another such map $\psi:(X,\mu)\to (Z,\theta)$ if
$\pi^{-1}(\S_Y) \supseteq \psi^{-1}(\S_Z)$ up to $\mu$-negligible
sets.  In this case we shall write $\pi \succsim \psi$ or $\psi
\precsim \pi$, and sometimes that $\psi$ is \textbf{$\mu$-virtually
a function of $\phi$} or that it is \textbf{$\mu$-virtually
$\phi^{-1}(\S_Y)$-measurable}. It is again a classical fact that in
the category of standard Borel spaces this notion of containment is
equivalent to the existence of a \textbf{factorizing} Borel map
$\phi:(Y,\nu) \to (Z,\theta)$ with $\psi = \phi\circ\pi$
$\mu$-almost everywhere, and that a measurable analog of the
Schroeder-Bernstein Theorem holds: $\pi \simeq \psi$ if and only if
in each direction such a $\phi$ may be chosen that is invertible
away from some negligible subsets of the domain and target. It is
clear that (up to set-theoretic niceties) this defines a partial
order on the class of $\simeq$-equivalence classes of Borel maps out
of the given space $(X,\mu)$.

Measure-respecting Borel maps from one probability space to another
comprise the simplest class of morphisms between such spaces, but in
this paper we shall sometimes find ourselves handling also a weaker
class of morphisms. Suppose that $Y$ and $X$ are standard Borel
spaces. Then by a \textbf{probability kernel from $Y$ to $X$} we
understand a function $P:Y\times\S_X \to [0,1]$ such that
\begin{itemize}
\item the map $y\mapsto P(y,A)$ is $\S_Y$-measurable for every $A \in
\S_X$;
\item the map $A \mapsto P(y,A)$ is a probability measure on $\S_X$
for every $y \in Y$.
\end{itemize}
Intuitively, such a kernel amounts to a `randomized map' from $Y$ to
$X$: rather than specify a unique image in $X$ for each point $y \in
Y$, it specifies only a probability distribution $P(y,\,\cdot\,)$
from which a point of $X$ could be chosen.  The first of the above
conditions is then the natural sense in which this assignment of a
probability distribution is measurable in $y$; indeed, a popular
alternative definition of probability kernel is as a measurable
function from $Y$ to the set $\Pr\,X$ of Borel probability measures
on $X$.  In ergodic theory this notion lies behind that of a
`quasifactor' (which assumes also a certain equivariance of this
map): see, for example, Chapter 8 of Glasner~\cite{Gla03}, where
this alternative convention and notation are used. We will write
$P:Y \stackrel{\rm{p}}{\to} X$ when $P$ is a probability kernel from
$Y$ to $X$.

Given a kernel $P:Y \stackrel{\rm{p}}{\to} X$ and a probability
measure $\nu$ on $Y$, we define the measure $P_\# \nu$ on $X$ by
\[P_\# \nu(A) := \int_YP(y,A)\,\nu(\d y);\]
this measure on $X$ can be interpreted as the law of a member of $X$
selected randomly by first selecting a member of $Y$ with law $\nu$
and then selecting a member of $X$ with law $P(y,\,\cdot\,)$. By
analogy with the case of a function between measurable spaces, we
will refer to this as the \textbf{pushforward} of $\nu$ by $P$. This
extends standard deterministic notation: given a measurable function
$\phi:Y \to X$, we may associate to it the deterministic probability
kernel given by $P(y,\,\cdot\,) = \delta_{\phi(y)}$ (the point mass
at the image of $y$ under $\phi$), and now $P_\#\nu$ is the usual
push-forward measure $\phi_\#\nu$.

Certain special probability kernels naturally serve as adjoints to
factor maps, in the sense of the following theorem.

\begin{thm}
Suppose that $Y$ and $X$ are standard Borel spaces, that $\mu$ is a
probability measure on $X$ and that $\phi:X \to Y$ is a measurable
factor map. Then, denoting the push-forward $\phi_\# \mu$ by $\nu$,
there is a $\nu$-almost surely unique probability kernel $P:Y
\stackrel{\rm{p}}{\to} X$ such that $\mu = P_\#\nu$ and which
\textbf{represents the conditional expectation with respect to
$\phi$}: for any $f \in L^1(\mu)$, the function
\[x_1 \mapsto \int_Xf(x)\,P(\phi(x_1),\d x)\]
is a version of the $\mu$-conditional expectation of $f$ with
respect to $\phi^{-1}(\S_Y)$.

We also write that this $P$ \textbf{represents the disintegration of
$\mu$ over $\phi$}.  A general probability kernel
$P:Y\stackrel{\rm{p}}{\to} X$ represents the disintegration over
$\phi$ of some measure that pushes forward onto $\nu$ if and only if
$\int_AP(x,\,\cdot\,)\,\nu(\d y)$ and $\int_BP(y,\,\cdot\,)\,\nu(\d
y)$ are mutually singular whenever $A\cap B = \emptyset$.
\end{thm}

\textbf{Proof}\hspace{5pt} See Theorem 6.3 in
Kallenberg~\cite{Kal02}. \qed

\subsection{Systems, subactions and factors}

In this paper we shall spend a great deal of time passing up and
down from systems to extensions or factors.  Moreover, sometimes one
system will appear as a factor of a `larger' system in \emph{several
different ways} (most obviously, when we work with a system that
appears under each coordinate projection from some self-joining).
For this reason the notational abuse of referring to one system as a
factor of another but leaving the relevant factor map to the
understanding of the reader, although popular and useful in modern
ergodic theory, seems dangerous here, and we shall carefully avoid
it.  In its place we substitute the alternative abuse, slightly
safer in our circumstances, of often referring only to the factor
maps we use, and leaving either their domain or target systems to
the reader's understanding.  Let us first set up some notation to
support this practice.

If $\G$ is a locally compact second countable topological group, by
a \textbf{$\G$-system} (or, if $\G$ is clear, just a
\textbf{system}) we understand a jointly measurable
probability-preserving action $T:\G\curvearrowright (X,\mu)$ on a
standard Borel probability space.  We will often alternatively
denote this space and action by $(X,\mu,T)$, or by a corresponding
single boldface letter such as $\bfX$. If $\L\leq \G$ we denote by
$T^{\ \uhr\L}:\L\curvearrowright (X,\mu)$ the action defined by
$(T^{\ \uhr\L})^\g := T^\g$ for $\g \in \L$, and refer to this as a
\textbf{subaction}, and if $\bfX = (X,\mu,T)$ is a $\G$-system we
write similarly $\bfX^{\ \uhr\L}$ for the system $(X,\mu,T^{\
\uhr\L})$ and refer to it as a \textbf{subaction system}.

A \textbf{factor} from one system $(X,\mu,T)$ to another $(Y,\nu,S)$
is a Borel map $\pi:X\to Y$ with $\pi_\#\mu = \nu$ and $\pi\circ T =
S\circ \pi$. Given such a factor, we sometimes write $T|_{\pi}$ to
denote the action $S$ with which $T$ is intertwined by $\pi$.

In this paper, given a globally invariant $\s$-subaglebra in $\bfX$,
a choice of factor $\pi:\bfX\to \bfY$ generating that
$\s$-subalgebra will sometimes be referred to as a
\textbf{coordinatization} of the $\s$-subalgebra. Importantly for
us, some choices of a coordinatizing factor $\pi$ may reveal some
underlying structure more clearly than others, and so we will
sometimes need to pass between coordinatizing factors. Given one
coordinatization $\pi:\bfX \to \bfY$ and an isomorphism $\psi: \bfY
\to \bfX$, we shall sometimes refer to the composition
$\psi\circ\pi$ as a \textbf{recoordinatization} of $\pi$. We will
also extend this terminology to that of \textbf{coordinatizations}
and \textbf{recoordinatizations} of families of factors of a system
in the obvious way in terms of the appropriate commutative diagram
of isomorphisms.

Given a $\G$-system $\bfX = (X,\mu,T)$, the $\s$-algebra $\S_X^T$ of
sets $A\in\S_X$ for which $\mu(A\triangle T^\g(A))=0$ for all $\g
\in \G$ is $T$-invariant, so defines a factor of $\bfX$. More
generally, if $\G$ is Abelian and $\L \leq\G$ then we can consider
the $\s$-algebra $\S_X^{T^{\ \uhr\L}}$ generated by all $T^{\
\uhr\L}$-invariant sets: we refer to this as the
\textbf{$\L$-isotropy factor} and write $\bfZ_0^{T^{\ \uhr\L}}$ for
some new system that we adopt as the target for a factor map
$\zeta_0^{T^{\ \uhr\L}}$ that generates $\S_X^{T^{\ \uhr\L}}$, and
$Z_0^{T^{\ \uhr\L}}$ for the standard Borel space underlying
$\bfZ_0^{T^{\ \uhr\L}}$.  Note that in this case the Abelianness
condition (or, more generally, the condition that $\L\unlhd \G$) is
needed for this to be a globally $T$-invariant factor. If $T_1$ and
$T_2$ are two commuting actions of the same Abelian group $\G$ on
$(X,\mu)$ then we can define a third action $T_1T_2^{-1}$ by setting
$(T_1T_2^{-1})^\g := T_1^\g T_2^{\g^{-1}}$, and in this case we may
write $\zeta_0^{T_1 = T_2}:\bfX\to\bfZ_0^{T_1 = T_2}$ in place of
$\zeta_0^{T_1T_2^{-1}}:\bfX\to \bfZ_0^{T_1T_2^{-1}}$. If $S\subseteq
\G$ and $\L$ is the group generated by $S$, we will sometimes write
$\bfZ_0^{T^{\ \uhr S}}$ in place of $\bfZ_0^{T^{\ \uhr\L}}$, and
similarly.

An important construction of new systems from old is that of
\textbf{relatively independent products}.  If $\bfY = (Y,\nu,S)$ is
some fixed system and $\pi_i:\bfX_i = (X_i,\mu_i,T_i)\to \bfY$ is an
extension of it for $i=1,2,\ldots,k$ then we define the relatively
independent product of the systems $\bfX_i$ over their factor maps
$\pi_i$ to be the system
\[\prod_{\{\pi_1 = \pi_2 = \ldots = \pi_k\}}\bfX_i = \Big(\prod_{\{\pi_1 = \pi_2 = \ldots = \pi_k\}}X_i,\bigotimes_{\{\pi_1 = \pi_2 = \ldots = \pi_k\}}\mu_i,T_1\times T_2\times\cdots\times T_k\Big)\]
where
\begin{multline*}
\prod_{\{\pi_1 = \pi_2 = \ldots = \pi_k\}}X_i :=
\{(x_1,x_2,\ldots,x_k)\in X_1\times X_2\times\cdots\times X_k:\\
\pi_1(x_1) = \pi_2(x_2) = \ldots = \pi_k(x_k)\},
\end{multline*}
\[\bigotimes_{\{\pi_1 = \pi_2 = \ldots = \pi_k\}}\mu_i = \int_Y \bigotimes_{i=1}^kP_i(y,\,\cdot\,)\,\nu(\d y)\]
and $P_i:Y\stackrel{\rm{p}}{\to} X_i$ is a probability kernel
representing the disintegration of $\mu_i$ over $\pi_i$. In case
$k=2$ we will write this instead as $\bfX_1\times_{\{\pi_1=
\pi_2\}}\bfX_2$, and in addition if $\bfX_1 = \bfX_2 = \bfX$ and
$\pi_1 = \pi_2 = \pi$ then we will abbreviate this further to
$\bfX\times_\pi\bfX$, and similarly for the individual spaces and
measures.

\subsection{Measurable selectors}

At several points in this paper we need to appeal to some basic
results on the existence of measurable selectors, often as a means
of making rigorous a selection of representatives of one or another
kind of data above the ergodic components of a non-ergodic system.

\begin{thm}\label{thm:meas-select}
Suppose that $(X,\S_X)$ and $(Y,\S_Y)$ are standard Borel spaces,
that $A \subseteq X$ is Borel and that $\pi:X\to Y$ is a Borel
surjection. Then the image $\pi(A)$ lies in the
$\nu^\rm{c}$-completion of $\S_Y$ for every Borel probability
measure $\nu$ on $(Y,\S_Y)$ with completion $\nu^\rm{c}$, and for
any such $\nu$ there is a map $f:B\to A$ with domain $B \in \S_Y$
such that $B \subseteq \pi(A)$, $\nu^\rm{c}(\pi(A)\setminus B) = 0$
and $\pi\circ f = \id_B$. \qed
\end{thm}

\textbf{Proof}\quad See, for example, 423O and its consequence
424X(h) in Fremlin~\cite{FreVol4}. \qed

\begin{dfn}[Measurable selectors]
We refer to a map $f$ as given by the above theorem as a
\textbf{measurable selector} for the set $A$.
\end{dfn}

\textbf{Remark}\quad We should stress that this is only one of
several versions of the `measurable selector theorem', due variously
to von Neumann, Jankow, Lusin and others. Note in particular that in
some other versions a map $f$ is sought that select points of $A$
for \emph{strictly} all points of $\pi(A)$. In the above generality
we cannot guarantee that a strictly-everywhere selector $f$ is
Borel, but only that it is Souslin-analytic and hence universally
measurable (of course, from this the above version follows at once).
On the other hand, if the map $\pi|_A$ is countable-to-one, then a
version of the result due to Lusin does guarantee a
strictly-everywhere Borel selector $f$. This version has already
played a significant r\^ole in our corner of ergodic theory in the
manipulation of the Conze-Lesigne equations (see, for
example,~\cite{ConLes84,FurWei96,BerTaoZie09}), and so we should be
careful to distinguish it from the above.  A thorough account of all
these different results and their proofs can be found in Sections
423, 424 and 433 of Fremlin~\cite{FreVol4}. \fin

In the right circumstances it is possible to strengthen
Theorem~\ref{thm:meas-select} to obtain a Borel selector that is
invariant under a group of transformations, by making use of a
coordinatization of the invariant factor.

\begin{prop}\label{prop:invar-meas-select}
Suppose that $(X,\S_X)$ and $(Y,\S_Y)$ are standard Borel spaces,
$A\subseteq X$ is Borel and $\pi:X\to Y$ is a surjective Borel map,
and in addition that $T:\G\curvearrowright (X,\S_X)$ is a jointly
measurable action of a locally compact second countable group such
that $\pi$ is a factor map, so $\pi\circ T^\g = S^\g\circ\pi$ for
some jointly measurable action $S:\G\curvearrowright (Y,\S_Y)$, and
that $A$ is $T$-invariant. Then for any $S$-invariant probability
measure $\nu$ on $(Y,\S_Y)$ with completion $\nu^\rm{c}$ there are
an $S$-invariant set $B\in \S_Y$ such that $B\subseteq \pi(A)$ and
$\nu^\rm{c}(\pi(A)\setminus B) = 0$ and an $S$-invariant map $f:B\to
A$ such that $\pi\circ f = \id_B$.
\end{prop}

\textbf{Proof}\quad Let $f_0:B_0 \to A$ be an ordinary measurable
selector as given by Theorem~\ref{thm:meas-select}, and let $\nu$ be
any $S$-invariant probability measure on $(Y,\S_Y)$. This $B_0$ must
be $\nu$-almost $S$-invariant, simply because $\pi(A)$ is
$T$-invariant and $\nu^{\rm{c}}(B\triangle \pi(A)) = 0$. Using local
compactness and second countability, let $(F_i)_{i\geq 1}$ be a
countable compact cover of $\G$, and also let $m_\G$ be a
left-invariant Haar measure on $\G$. From the joint measurability of
$T$ it follows that the set
\begin{multline*}
B := \big\{y\in Y:\ m_\G\{\g \in \G:\ S^\g(y) \in Y\setminus B_0\}
= 0\big\}\\
= \bigcap_{i\geq 1}\big\{y\in Y:\ m_\G\{\g \in F_i:\ S^\g(y) \in
B_0\} = m_\G(F_i)\big\}
\end{multline*}
is Borel, $T$-invariant and satisfies $\nu(B_0\triangle B) = 0$.

We now let $\zeta:(Y,\S_Y,\nu)\to (Z,\S_Z,\theta)$ be any
coordinatization of the invariant factor $\S_Y^T$; it is easy to see
that this may be chosen so that there exists some $C \in \S_Z$ such
that $B =\zeta^{-1}(C)$.  We can now use $B$ and $\zeta$ to `tidy
up' our original selector $f_0$.  Indeed, by the $S$-invariance of
$\zeta$ and the fact that for every $y \in B$ we have $S^\g(y)\in
B_0$ for some (indeed, almost all) $\g \in \G$, we must have
$\zeta(B_0) \supseteq C$.  Therefore by applying the ordinary
Measurable Selector Theorem a second time we can find a Borel subset
$D \in \S_Z$ with $D \subseteq C$ and $\theta(C\setminus D) = 0$ and
a Borel section $\eta:D\to B_0$ such that $\zeta\circ\eta = \id_D$;
and so now replacing $B_0$ with $B$ and the map $f_0$ with
$f:y\mapsto f_0(\eta(\zeta(y)))$ completes the proof. \qed

\begin{dfn}[Invariant measurable selectors]
We refer to a map $f$ as given by the above proposition as a
\textbf{$T$-invariant measurable selector} for the set $A$.
\end{dfn}

\subsection{Background from group theory}\label{subs:group}

We collect here some standard group theoretic definitions and
results for future reference.

\begin{dfn}[Core]
If $G$ is a group and $H \leq G$ we denote by $\rm{Core}_G(H)$ the
\textbf{core of $H$ in $G$}: the largest subgroup of $H$ that is
normal in $G$.  It is clear that this exists and equals
$\bigcap_{g\in G}g^{-1}Hg$.  If $G$ is compact and $H$ is closed
then so is $\rm{Core}_G(H)$.

If $\rm{Core}_G(H) = \{1_G\}$ we shall write that $H$ is
\textbf{core-free} in $G$.
\end{dfn}

\begin{dfn}[Full one-dimensional projections; slices]
Given two groups $G_1$ and $G_2$ and a subgroup $M \leq G_1\times
G_2$, and writing $\pi_i:G_1\times G_2\to G_i$ for the two
coordinate projections, we say that $M$ has \textbf{full
one-dimensional projections} if $\pi_i(M) = G_i$ for $i=1,2$.

We refer to the subgroups
\[H_1 := \pi_1(M\cap(G_1\times\{1_{G_2}\}))\]
and
\[H_2 := \pi_2(M\cap(\{1_{G_1}\}\times G_2))\]
as the \textbf{first} and \textbf{second slices} of $M$
respectively.
\end{dfn}

It is a classical observation of Goursat (see, for example, Section
1.6 of Schmidt~\cite{Sch94}) that $M$ has full one-dimensional
projections and trivial first and second slices if and only if it is
the graph of an isomorphism $\Phi:G_2\to G_2$.  If the slices are
non-trivial, we do at least have the following.

\begin{lem}\label{lem:groupcorrespondence}
If $M \leq G_1\times G_2$ has full one-dimensional projections then
its slices satisfy $H_i \unlhd G_i$ for $i=1,2$.
\end{lem}

\textbf{Proof}\quad By symmetry it suffices to treat the case $i=1$.
Let $r_1 \in G_1$.  Since $\pi_1(M) = G_1$ we can find $r_2 \in G_2$
such that $(r_1,r_2)\in M$.  It is now easy to check that
\begin{eqnarray*}r_1H_1 &=& \{g\in G_1:\ (r_1^{-1}g,e) \in M\}\\
&=& \{g\in G_1:\ (r_1,r_2)(r_1^{-1}g,e) \in M\}\\
&=&\{g\in G_1:\ (g,r_2) \in M\}\\
&=& \{g\in G_1:\ (gr_1^{-1},e)(r_1,r_2) \in M\}\\
&=& \{g\in G_1:\ (gr_1^{-1},e) \in M\} = H_1r_1.
\end{eqnarray*}
Since $r_1$ was arbitrary, $H_1$ is normal, as required. \qed

Given a compact group $G$ or one of its homogeneous spaces $G/H$ we
shall always consider it endowed with its usual Borel structure and
Haar probability measure, which we shall denote by $m_G$ or
$m_{G/H}$.

If $U$ is a compact metrizable group then we write $\Clos\,U$ for
its collection of closed subsets endowed with the Vietoris topology
(see, for example, 2.7.20, 3.12.27 and 4.5.23 of
Engelking~\cite{Eng89}; as is standard, this is also compact and
metrizable) and the associated standard Borel structure and $\Lat\,U
\subseteq \Clos\,U$ for the further Veitoris-closed subfamily of
closed sub\emph{groups} with its induced standard Borel structure.
In this setting of subgroups of compact metrizable groups, the
Vietoris topology is easily seen to coincide with the Fell topology
and the Chabauty topology, both of which also commonly appear in the
study of lattices of closed subgroups; see Subsection 2.1 of Conze
and Raugy~\cite{ConRau07} and the references given there. This
topology and Borel structure can be understood in terms of Haar
measures in the following standard way.

\begin{lem}
The Vietoris topology and measurable structure on $\Lat\,U$ coincide
with the pullbacks of the vague topology and measurable structure
under the Haar-measure map $H \mapsto \mu_H$. \qed
\end{lem}

\section{Direct integrals of homogeneous-space
data}\label{sec:dirint}

In this section we give the rigorous definition of a `direct
integral' of measurably-varying homogeneous spaces and of the lifted
transformation acting on it, and establish some of their elementary
properties. We build such an extension $X$ as a union of different
fibres $G_y/K_y$ above each $y \in Y$, the fibre actually depending
only on $\zeta_0^S(y) \in Z_0^S$, and we extend $S$ to an action $T$
on $X$ using a cocycle constrained to lie at (almost) every point in
the relevant fibre. We enforce a suitable measurable structure by
drawing $G_y$ and $K_y$ from among the compact subgroups of some
fixed `repository' group, subject to the condition $K_y \leq G_y$,
measurably for the Vietoris measurable structure on such subgroups.

\begin{dfn}[Measurable homogeneous space data]
Let $Y$ be a standard Borel space and $U$ a fixed compact metrizable
group. By \textbf{measurable compact group data on $Y$ with fibre
repository $U$} we understand a map $Y\to \Lat\,U:y\mapsto G_y$ that
is measurable for the Vietoris Borel structure on $\Lat\,U$.  We
shall usually denote such a map by $G_\bullet$, and will often omit
explicit mention of the fibre repository $U$. More generally, by
\textbf{measurable compact homogeneous space data on $Y$ with fibre
repository $U$} we understand a pair $(G_\bullet,K_\bullet)$ of
measurable compact group data with repository $U$ such that $K_y
\leq G_y$ for every $y$.  We shall usually denote this pair instead
by $G_\bullet/K_\bullet$, and think of it as a measurable assignment
of the compact homogeneous space $G_y/K_y$ to each point $y \in Y$.
\end{dfn}

\begin{dfn}[Direct integral of measurable homogeneous space data]
Given a standard Borel probability space $(Y,\nu)$ and measurable
compact homogeneous space data as above, we shall define their
\textbf{direct integral} to be the subset
\[\{(y,gK_y):\ y\in Y,\,g \in G_y\} \subseteq Y \times \rm{Clos}\,U,\]
which we denote by $Y\ltimes G_\bullet/K_\bullet$.  This is easily
verified to be standard Borel for the relevant product measurable
structure, and we will always assume it to be endowed with the
restriction of that measurable structure.

On this space we define the \textbf{direct integral} measure
$\nu\ltimes m_{G_\bullet/K_\bullet}$ by
\[\nu\ltimes m_{G_\bullet/K_\bullet}(A) := \int_Y\delta_y\otimes m_{G_y/K_y}\big(A\cap(\{y\}\times G_y/K_y)\big)\,\nu(\d
y).\]

Given another measurable assignment of subgroup data $H_\bullet \leq
G_\bullet$, we define analogously the direct integrals $(Y\ltimes
H_\bullet\backslash G_\bullet,\nu\ltimes m_{H_\bullet\backslash
G_\bullet})$ of the spaces of right-cosets and $\big(Y\ltimes
(H_\bullet\backslash G_\bullet/K_\bullet),\nu\ltimes
m_{H_\bullet\backslash G_\bullet/K_\bullet}\big)$ of the spaces of
double cosets.
\end{dfn}

\textbf{Remark}\quad We will rarely remark again on the assumption
that the fibre repository $U$ be metrizable, but this will always be
implicit. This ensures that the above construction keeps us within
the category of standard Borel spaces (and it will be a natural
consequence of the non-ergodic Furstenberg-Zimmer theory applied to
such spaces), and will occasionally be important for proofs (such as
in Lemma~\ref{lem:enable-quotienting} below). One could attempt to
construct an extended theory that allows non-metrizable fibre groups
and works instead in the larger category of perfect measure spaces
(see 342K of Fremlin~\cite{FreVol3}), but we will not do so here.
\fin

\begin{dfn}[Cocycle-sections]
Suppose that $(Y,\nu)$, $U$ and $G_\bullet/K_\bullet$ are as above,
that $\G$ is a locally compact second countable group and that
$S:\G\curvearrowright (Y,\nu)$, and suppose further that the group
data $y\mapsto G_y$ and $y\mapsto K_y$ are $S$-invariant. Then a
\textbf{cocycle-section} of $G_\bullet$ over $S$ is a measurable
cocycle $\rho:\G\times Y\to U$ over $S$ such that $\rho(\g,y) \in
G_y$ for every $\g \in \G$ and $y \in Y$. We shall denote such a
cocycle-section by $\rho:\G\times Y \to G_\bullet$.
\end{dfn}

\textbf{Remark}\quad Note that in the setting of a general locally
compact second countable group $\G$, the definition that $\rho$ be a
cocycle over $S$ demands only that $\rho(\g_1\g_2,y) =
\rho(\g_1,S^{\g_2}y)\cdot\rho(\g_2,y)$ for $\nu$-almost every $y \in
Y$ for strictly every $\g_1$ and $\g_2$ (see, for instance, Section
4.2 of Zimmer~\cite{Zim84}), where the negligible set of `bad' $y$
is allowed to vary with $(\g_1,\g_2)$; and that by convention two
cocycles are equivalent if they agree $\nu$-almost surely for
strictly every $\g$. In view of this, we lose no generality in
asking that $\rho(\g,y) \in G_y$ for strictly every $y$ and $\g$,
rather than for almost every $y$ for strictly every $\g$, since in
the latter case we may simply adjust $\rho$ to equal $1$ on the
Borel set where is falls outside the specified repository, and this
changes each $\rho(\g,\,\cdot\,)$ on only a $\nu$-negligible set for
strictly every $\g$. \fin

Finally, we can define our class of extensions.

\begin{dfn}[Extensions by measurable homogeneous space data]
Suppose that $\bfY = (Y,\nu,S)$ and $G_\bullet/K_\bullet$ are as
above, that the group data $y\mapsto G_y$ and $y\mapsto K_y$ are
$S$-invariant and that $\rho:\G\times Y\to G_\bullet$ is a
cocycle-section over $S$. Then the \textbf{extension of $\bfY$ by
the data $(G_\bullet/K_\bullet,\rho)$} is the action $T$ of $\G$ on
$(Y\ltimes G_\bullet/K_\bullet,\nu\ltimes m_{G_\bullet/K_\bullet})$
given by
\[T^\g(y,gK_y) := (S^\g y,\rho(\g,y)gK_y);\]
it is routine to verify that this is measurable and
measure-preserving.

We will often denote this extended system by $\bfY\ltimes
(G_\bullet/K_\bullet,m_{G_\bullet/K_\bullet},\rho)$. It clearly
admits $\bfY$ as a factor simply by projecting out the fibre
coordinate; we will refer to this as the \textbf{canonical} factor
map. The data $G_\bullet/K_\bullet$ and cocycle-section $\rho$ are
together \textbf{relatively ergodic} if the extension $\bfY\ltimes
(G_\bullet/K_\bullet,m_{G_\bullet/K_\bullet},\rho) \to \bfY$ through
the canonical map is relatively ergodic.
\end{dfn}

\textbf{Remarks}\quad\textbf{1.}\quad In light of the Peter-Weyl
Theorem (treated in most standard texts on compact group
representations, such as in Section III.3 of Br\"ocker and tom
Dieck~\cite{Brotom85}) all compact metrizable groups can be realized
isomorphically, albeit highly non-uniquely, as closed subgroups of a
suitably large direct product of unitary groups, say $U :=
\prod_{n\geq 1}\rm{U}(n)^{\bbN}$.  This suggests that such a direct
product should suffice as a compact repository for all purposes, and
indeed this can be proved with just a little work; however, this
result seems to contribute little to the theory, and so we will not
present it here.  Note, however, that it is also precisely such
direct products of unitary groups that will emerge naturally as
repositories in the non-ergodic Furstenberg-Zimmer inverse theory of
Section~\ref{sec:FZ} below.

\quad\textbf{2.}\quad In view of the condition that $G_\bullet$ and
$K_\bullet$ are $S$-invariant, given a coordinatization
$\zeta_0^S:Y\to Z_0^S$ of the $S$-isotropy factor we could
alternatively work with compact measurable group data defined
initially as functions on the space $Z_0^S$ and then lifted through
$\zeta_0^S$. We will occasionally use this alternative description
when it is notationally convenient. \fin

The following related definition will also occasionally be useful.

\begin{dfn}[Opposite extensions by measurable group
data]\label{dfn:opp} If $\bfY$, $G_\bullet$ and $\rho$ are as above,
then they also define an extended action $T$ on $(Y\ltimes
G_\bullet,\nu\ltimes m_{G_\bullet})$ by
\[T^\g(y,g) := (S^\g y,g\rho(\g,y)^{-1}):\]
this is the \textbf{opposite extension of $\bfY$ by the data
$(G_\bullet,\rho)$}, and we will denoted this $T$ by $S\ltimes
\rho^{\rm{op}}$.
\end{dfn}

\textbf{Remark}\quad In fact we always have
\begin{center}
$\phantom{i}$\xymatrix{\bfY\ltimes
(G_\bullet,\rho)\ar[dr]_{\rm{canonical}}\ar@{<->}[rr]^-\cong &
&\bfY\ltimes (G_\bullet,\rho^{\rm{op}})\ar[dl]^{\rm{canonical}}\\ &
\bfY }
\end{center}
through the fibrewise isomorphism $(y,g)\mapsto (y,g^{-1})$.  The
use of opposite extensions will matter to us in situations where we
have two different actions on the extended space, one by a cocycle
and one by an opposite cocycle. \fin

Before leaving this section, it is worth noting one way in which
some redundancy in the above definition can be removed.

\begin{dfn}
Homogeneous space date $G_\bullet/K_\bullet$ over $(Y,\nu)$ is
\textbf{core-free} if $K_y$ is core-free in $G_y$ almost everywhere.
\end{dfn}

\begin{lem}\label{lem:enable-quotienting}
Suppose $\bfY := (Y,\nu,S)$ is a $\G$-system, that
$G_\bullet/K_\bullet$ are measurable $S$-invariant homogeneous space
data on $Y$ with repository $U$ and that $\rho:\G\times Y\to
G_\bullet$ is a cocycle-section over $S$.  If in addition the group
$K_y$ is normal in $G_y$ for $\nu$-almost every $y$, then there are
a fibre repository $U'$, measurable $S$-invariant group data $G'$ on
$Y$ and a measurable $S$-invariant family of isomorphisms
$\Psi_y:G_y/K_y\to G'_y$ such that the map $(y,gK_y)\mapsto
(y,\Psi_y(gK_y))$ defines an isomorphism of extensions
\begin{center}
$\phantom{i}$\xymatrix{ \bfY\ltimes
(G_\bullet/K_\bullet,m_{G_\bullet/K_\bullet},\rho)\ar[dr]_{\rm{canonical}}\ar@{<->}[rr]^-\cong
& & \bfY\ltimes
(G'_\bullet,m_{G'_\bullet},\rho')\ar[dl]^{\rm{canonical}}\\
& \bfY }
\end{center}
with $(\g,y)\mapsto\rho'(\g,y):= \Psi_y(\rho(\g,y)):\G\times Y\to
G'$.
\end{lem}

\textbf{Proof}\quad This rests on the construction of the new fibre
repository for the quotient groups $G_\bullet/K_\bullet$. For $y \in
Y$ let $\frH_y \leq L^2(m_U)$ be the separable Hilbert subspace of
square-integrable functions on $U$ invariant under left-rotation by
$K_y$.  This is an $S$-invariant measurable family of separable
Hilbert spaces (in the sense familiar from the analysis of group
representations and von Neumann algebras; see, for instance,
Mackey~\cite{Mac76}), and so we can partition $Y$ into $S$-invariant
measurable subsets $A_1$, $A_2$, \ldots, $A_\infty$ and for each
$n\in\bbN\cup\{\infty\}$ select an $S$-invariant measurable family
of isomorphisms $\Phi_y:\frH_y \to \frH_n'$ for $y \in A_n$, where
$\frH_n'$ is some fixed $n$-dimensional reference complex Hilbert
space and $\frH_1' \leq \frH_2'\leq \ldots \leq \frH_\infty'$.

Now let $\pi_y:G_y \curvearrowright \frH_\infty$ for $y \in A_n$ be
the representation that results from first restricting the
left-regular representation of $G_y$ on $L^2(m_U)$ to $\frH_y$
(which is possible when $K_\bullet \unlhd G_\bullet$, hence almost
everywhere), then composing with $\Phi_y$ to obtain a representation
on $\frH_n'$ and finally extending this to act on $\frH_\infty'$ by
acting trivially on $\frH_\infty' \ominus \frH_n'$.

This defines an $S$-invariant measurable family of representations
$\pi_y$ of $G_y$ for $y \in Y$ outside some $\nu$-conegligible
subset, and such that $K_y = \ker \pi_y$ almost surely. Next, it is
easy to see that the decomposition of $\pi_y$ into
finite-dimensional representations given by the Peter-Weyl Theorem
is measurable in $y$ (for example, since they may recovered as the
spectral projections of each of a countable dense subfamily of all
the measurable selections over $y$ of $\pi_y$-invariant compact
operators on $\frH_\infty$). Hence this decomposition gives a
measurable family of continuous homomorphic embeddings $G_\bullet
\longrightarrow \prod_{n\geq 1}\rm{U}(n)^\bbN$ with kernels
$K_\bullet$, and so letting $G'_\bullet$ be the image group data of
these embeddings they define a measurable family of isomorphisms
$\Psi_\bullet$ such that defining $\rho'$ as above and applying
$\Psi_\bullet$ fibrewise on $Y\ltimes G_\bullet/K_\bullet$ gives the
desired isomorphism of extensions. \qed

\begin{cor}
If $\bfX = \bfY\ltimes
(G_\bullet/K_\bullet,m_{G_\bullet/K_\bullet},\rho)$ is an extension
by homogeneous space data, then it is isomorphic (as an extension of
$\bfY$ through the canonical map) to an extension by core-free
homogeneous space data.
\end{cor}

\textbf{Proof}\quad Let $U$ be the repository and let $L_y :=
\bigcap_{g\in G_y}g^{-1}K_y g$ be the pointwise core of $K_y$ in
$G_y$.  First observe that for any $u \in U$ the set
\[\Big\{(G,K):\ u \in \bigcap_{g\in G}g^{-1}Kg\Big\} = \{(G,K):\ K \ni gug^{-1}\ \forall g\in G\}\]
it open in $(\Lat\,U)^2$, since if $(K,G)$ does not lie in this set
then there are a closed set $V_1 \subseteq U$ with nonempty interior
and an open set $V_2 \subseteq U$ such that $K\cap V_1 = \emptyset$,
$G\cap V_2 \neq \emptyset$ and $V_2uV_2^{-1} \subseteq V_1$.  It
follows that the map $(\Lat\,U)^2\to \Lat\,U:(G,K)\mapsto
\bigcap_{g\in G}g^{-1}Kg$ is measurable, and hence that $L_\bullet$
is measurable group data.

Now by the preceding lemma we can select a measurable family of
embeddings of the groups $G_y/L_y$ into a suitably-modified
repository to obtain an isomorphism of systems
\begin{center}
$\phantom{i}$\xymatrix{ \bfY\ltimes
(G_\bullet/L_\bullet,m_{G_\bullet/L_\bullet},\rho)\ar[dr]_{\rm{canonical}}\ar@{<->}[rr]^-\cong
& & \bfY\ltimes
(G'_\bullet,m_{G'_\bullet},\rho')\ar[dl]^{\rm{canonical}}\\
& \bfY }
\end{center}
corresponding to a continuous group isomorphism $G_y/L_y\to G'_y$ at
almost every $y$.  Under these isomorphisms the subgroups $K_y \leq
G_y$ correspond measurably to some $K'_y \leq G'_y$ so that
$(G'_y,K'_y)\cong (G_y/L_y,K_yL_y/L_y)$, and so observing from its
definition that $L_yK_y/L_y$ is always core-free in $G_y/L_y$, this
completes the proof. \qed

\section{Mackey Theory in the non-ergodic setting}\label{sec:Mackey}

We will now move on to a more detailed analysis of extensions by
homogeneous space data, and more specifically of their invariant
factors and relatively ergodic measures.  Many of the ideas that
follow are nearly direct translates to our setting of those of
Mackey in the case of an ergodic base system, and we will follow
quite closely their treatment in Section 3.5 of
Glasner~\cite{Gla03}.

In fact, more is true: earlier work on multiple recurrence and
nonconventional ergodic averages has already encountered the
possibility of a measurably-varying Mackey group within an extension
of a non-ergodic base system by a fixed overall group. This
technicality arises in the work of Meiri~\cite{Mei90} on correlation
sequences arising from probability-preserving systems, of
Furstenberg and Weiss~\cite{FurWei96} on certain polynomial
nonconventional ergodic averages and more recently in Ziegler's
approach in~\cite{Zie07} to convergence of linear nonconventional
averages for powers of a single transformation.  For example, during
the analysis in~\cite{FurWei96} a homogeneous skew-product extension
of ergodic systems $(X,\mu,T) = (Y,\nu,S)\ltimes (G,m_G,\rho)$ is
three-fold joined to itself, to give a measure on $X^3$ that is
invariant for a transformation of the form $T^r\times T^s\times T^t$
but which is not ergodic for that transformation. This system is now
coordinatized as an extension of an action on $Y^3$ by the group and
cocycle $(G^3,(\rho^{(r)},\rho^{(s)},\rho^{(t)}))$, but since the
base is no longer ergodic the description of the ergodic components
of the overall system requires the possibility that the Mackey group
can vary among the closed subgroups of $G^3$ (a possibility that is
then discounted by an argument showing that they are all actually
conjugate, and so may in fact be taken to be constant; we shall see
a similar trick in Subsection~\ref{subs:charfactors} below).

The only extra subtlety for which we must allow here is that the
overall group $G_\bullet$ now also varies measurably.  This will
require us to work rather harder in setting up the proof, although
the overall idea is very similar to those mentioned above.  For this
reason, although we have included complete proofs here, we refer the
reader to these other sources, and also Section 3.5 of
Glasner~\cite{Gla03}, for relevant background.

\subsection{Ergodic decompositions and Mackey group data}

The Mackey Theory describes the invariant factor of an extension
$\bfX = \bfY\ltimes (G_\bullet,m_{G_\bullet},\rho)$ in terms of the
invariant factor of $\bfY$ and the data
$(G_\bullet,m_{G_\bullet},\rho)$ of the extension.  Here it will
prove convenient to treat $G_\bullet$ as varying over the factor
space $Z_0^S$ of ergodic components, lifted to $Y$ for the purpose
of defining the extended system.

\begin{thm}[Mackey Theorem in the non-ergodic
case]\label{thm:nonergMackey} Suppose that $\bfX$ is the group-data
extension $\bfY\ltimes (G_\bullet,m_{G_\bullet},\rho)$ and that
$\zeta_0^S:\bfY \to \bfZ_0^S$ is a coordinatization of the base
isotropy factor, and let $\pi:X = Y\ltimes G_\bullet\to Y$ and
$\theta:Z_0^S\ltimes G_\bullet \to Z_0^S$ be the canonical factor
maps. Then there are measurable subgroup data $K_\bullet \leq
G_\bullet$ on $Z_0^S$, a $T$-invariant map
\[\phi:(X,\mu) \to
(Z_0^S,{\zeta_0^S}_\#\nu)\ltimes (K_\bullet\backslash
G_\bullet,m_{K_\bullet\backslash G_\bullet})\] and a section $b:Y
\to G_\bullet$ such that
\begin{enumerate}
\item[(1)] $\phi$ coordinatizes $\bfZ_0^T$ and the following
diagram commutes:
\begin{center}
$\phantom{i}$\xymatrix{ Y\ltimes
G_\bullet\ar[r]^-\phi\ar[d]_{\rm{canonical}} & Z_0^S\ltimes
K_\bullet\backslash
G_\bullet\ar[d]^{\rm{canonical}}\\
Y\ar[r]_{\zeta_0^S} & Z_0^S }
\end{center}
\item[(2)] $\phi(y,g) = (\zeta_0^S(y),K_{\zeta_0^S(y)}b(y)g)$ for
$\mu$-almost every $(y,g)$; \item[(3)] the cocycle-section
$(\g,y)\mapsto b(S^\g y)\rho(\g,y)b(y)^{-1}$ takes a value in
$K_{\zeta_0^S(y)}$ for $\nu$-almost every $y$ for every $\g$;
\item[(4)](\textbf{Conjugate minimality}) if $K'_\bullet \leq G_\bullet$ is another measurable
assignment of compact subgroup data on $Z_0^S$ and $b':Y\to
G_\bullet$ another section such that the cocycle-section
$(\g,y)\mapsto b'(S^\g y)\rho(\g,y)b'(y)^{-1}$ takes a value in
$K'_{\zeta_0^S(y)}$ for $\nu$-almost every $y$ for every $\g$, then
there is a section $c:Z_0^S\to G_\bullet$ such that
\[c(s) \cdot K'_s \cdot c(s)^{-1} \geq K_s\] for
$(\zeta_0^S)_\#\nu$-almost every $s$;
\item[(5)] if $P:Z_0^S \stackrel{\rm{p}}{\longrightarrow} Y$ is a version of
the disintegration of $\nu$ over $\zeta_0^S$, then the probability
kernel $(s,K_sg') \stackrel{\rm{p}}{\mapsto} P(s,\,\cdot\,)\ltimes
m_{b(\bullet)^{-1}K_sg'}$ is a version of the disintegration of
$\mu$ over $\phi$.
\end{enumerate}
\end{thm}

\textbf{Remark}\quad Clearly with hindsight we can take the property
(2) above as defining $\phi$; the point, however, is that we will
obtain $K_\bullet$ and $\phi$ first and then show that $\phi$ takes
this form for some $b$. \fin

The proof of this theorem will require some initial constructions
and an enabling lemma concerning the measurable selection of generic
points.

First let $\zeta_0^T:\bfX \to \bfZ_0^T$ be some coordinatization of
the isotropy factor of the large system with the property that
$\zeta_0^S\circ\pi$ factorizes through the natural factor
$\xi:\bfZ_0^T\to \bfZ_0^S$ in the sense of the commutative diagram
\begin{center}
$\phantom{i}$\xymatrix{ \bfX\ar[r]^{\zeta_0^T}\ar[d]_\pi &
\bfZ_0^T\ar[d]^{\xi}\\ \bfY\ar[r]_{\zeta_0^S} & \bfZ_0^S}
\end{center}
(it is a standard fact that this is possible; see, for example,
Section 2.2 of Glasner~\cite{Gla03}).

We need a formal way to work with the action by right-multiplication
of the fibre $G_y$ on itself (the difficulty being, of course, that
global constructions based on this pointwise-varying action need to
kept measurable).  To this end we define the map
\[\tau_0:X \times_{\{\zeta_0^S\circ \pi = \theta\}}(Z_0^S\ltimes
G_\bullet) \to X:((y,g'),(s,g))\mapsto (y,g'g);\] intuitively,
$\tau_0(x,(s,g))$ gives the image of $x$ under right-multiplication
by $g$, which is an element of the group $G_s$ over $s =
\zeta_0^S(\pi(x))$ that acts on the fibre above $\pi(x)$. This map
$\tau_0$ can be well-defined only for those tuples with $s =
\zeta_0^S(\pi(x))$, hence the need for the relative self-product in
the specification of its domain; however, with this restriction in
place it is easily seen to be measurable.

We can now work with the domain of $\tau_0$ as a system in its own
right under the action $(T\times \id_{Z_0^S\ltimes G_\bullet})$. It
is clear that the measure $\mu\otimes_{\{\zeta_0^S\circ \pi =
\theta\}}({\zeta_0^S}_\#\nu\ltimes m_{G_\bullet})$ is $(T\times
\id_{Z_0^S\ltimes G_\bullet})$-invariant and that $\tau_0$ itself is
a factor map from the resulting system onto $(X,\mu,T)$.

From $\tau_0$ we now define the composition
\[\tau:X\times_{\{\zeta_0^S\circ \pi=\theta\}} (Z_0^S\ltimes
G_\bullet) \stackrel{\tau_0}{\longrightarrow}
X\stackrel{\zeta_0^T}{\longrightarrow} Z_0^T.\]  Heuristically this
assigns to the pair $(x,(s,g))$ the ergodic component of $(X,\mu,T)$
that contains the image of $x$ under the right-multiplication by $g$
acting on its fibre.  Of course, the whole point is that different
points within a single fibre of $\pi$ will generally lie in
different ergodic components, and this map $\tau$ reports on this
dependence.

Informally, the proof of Theorem~\ref{thm:nonergMackey} now proceeds
by selecting a representative $p(s) \in X$ above each $s \in Z_0^S$
and then defining $K_s$ to be the subgroup of those $g \in G_s$ such
that the $T$-ergodic component of $p(s)$ does not change upon
right-multiplication by $g$ inside the $\pi$-fibre of $s$: that is,
such that $\tau(p(s),(s,g)) = \zeta_0^T(p(s))$.  We need $p$ to
select points that are sufficiently `generic' in the fibres above
$s$, in the sense made precise by the following lemma.

\begin{lem}
In the setting of Theorem~\ref{thm:nonergMackey}, we can find a
Borel measurable section $p:Z_0^S \to X$ of the factor map
$\zeta_0^S\circ \pi$ such that the probability kernel
\[P(s,A) := m_{G_s}\{g\in G_s:\ \tau(p(s),(s,g)) \in A\},\quad A \in \S_{Z_0^T},\]
is a version of the disintegration of ${\zeta_0^T}_\#\mu$ over
$\xi:Z_0^T\to Z_0^S$.
\end{lem}

\textbf{Proof}\quad We first define a probability kernel $P':X
\stackrel{\rm{p}}{\to} Z_0^T$ by
\[P'((y,g'),A) := m_{G_{\zeta_0^S(y)}}\{g\in G_{\zeta_0^S(y)}:\ \tau((y,g'),(\zeta_0^S(y),g))\in A\}.\]
Intuitively this takes a point $x = (y,g')$, replaces it with an
average over the fibre $\pi^{-1}\{y\}$, regarded as a copy of the
fibre group $G_{\zeta^S_0(y)}$, and then pushes the resulting
probability measure down to $Z_0^T$.

It is easy to see that this satisfies the measurability conditions
of a probability kernel, since we have
\[P'((y,g'),A) = P''((y,g'),(\zeta_0^T)^{-1}(A))\]
where $P'':X\stackrel{\rm{p}}{\to} X$ is the manifestly measurable
probability kernel
\[P''((y,g'),\,\cdot\,) := (\delta_y\otimes m_{G_{\zeta_0^S(y)}})(\,\cdot\,).\]

In addition the above definition implies that
\[P'((y,g'),\,\cdot\,) = P'((y,g'g),\,\cdot\,) = P'((S^\g
y,\rho(\g,y)g'),\,\cdot\,)\] for any $g \in G_{\zeta_0^S(y)}$ and
$\g \in \G$, firstly in view of the averaging over
$G_{\zeta_0^S(y)}$ and secondly because we take the image under the
invariant function $\zeta_0^T$.  Therefore the function $x \mapsto
P'(x,\,\cdot\,)$, regarded as a Borel map from $X$ to the space of
Borel probability measures on $Z_0^T$ with its usual Borel
structure, both factorizes through $\pi$ and is then $S$-invariant.
Therefore there is some $P:Z_0^S\stackrel{\rm{p}}{\to}Z_0^T$ such
that $P' = P\circ\zeta_0^S\circ\pi$, $\mu$-almost surely.  This $P$
must be a version of the disintegration of ${\zeta_0^T}_\#\mu$ over
$\xi$.

Setting
\[A_s := \big\{x \in X:\ \zeta_0^S(\pi(x)) = s\ \hbox{and}\ P'(x,\,\cdot\,) = P(s,\,\cdot\,)\big\},\]
this is now the set of `generic points' above $s$ from which we need
to select $p(s)$. From the above relation between $P'$ and $P$ we
deduce that that $P(s,A_s) = 1$ for ${\zeta_0^S}_\#\nu$-almost every
$s$. In addition, since $A_s$ is the section above $s$ of the Borel
set
\[A := \big\{(s,x) \in Z_0^S\times X:\
\zeta_0^S(\pi(x)) = s\ \hbox{and}\ P'(x,\,\cdot\,) =
P(s,\,\cdot\,)\big\},\] by the Measurable Selector
Theorem~\ref{thm:meas-select} we can choose a measurable selector
$p$ of $\zeta_0^S\circ\pi$ such that $(s,p(s)) \in A$ for
${\zeta_0^S}_\#\nu$-almost every $s$.  This selector now has the
properties claimed. \qed

\textbf{Proof of
Theorem~\ref{thm:nonergMackey}}\quad\textbf{(1)}\quad Given the
measurable selector of the above lemma, define \[K_s := \{g \in
G_s:\ \tau(p(s),(s,g)) = \zeta_0^T(p(s))\}.\] It is clear from the
definition of $\tau$ that $\tau(p(s),(s,g)) = \zeta_0^T(p(s))$ if
and only if $\tau(p(s),(s,gg')) = \tau(p(s),(s,g'))$ for every $g'
\in G_s$, and hence $K_s$ is a closed subgroup of $G_s$ for almost
every $s$.

Also, we have
\[A := \{(s,g) \in Z_0^S\times G_\bullet:\ g \in K_s\} = \{(s,g) \in Z_0^S\times G_\bullet:\ \tau(p(s),(s,g)) = \zeta_0^T(p(s))\},\]
and so this is a Borel subset of $Z_0^S\times G_\bullet$. Letting
$U$ be the repository, it follows that for any open $V_1$, $V_2$,
\ldots, $V_k \subseteq U$ and closed $W \subseteq U$ we have
\begin{multline*}
\{s\in Z_0^S:\ K_s\cap V_i \neq \emptyset\ \forall i \leq k\ \&\ K_s
\cap W = \emptyset\}\\ = \theta\Big(A\cap\bigcap_{i\leq
k}(Z_0^S\times V_i)\Big)\Big\backslash \theta(A\cap (Z_0^S\times
W)),
\end{multline*}
(recalling that $\theta$ is the projection onto the first
coordinate) and so from the Measurable Selector
Theorem~\ref{thm:meas-select} this is a universally measurable set.
Allowing $V_1$, $V_2$, \ldots, $V_k$ and $U\setminus W$ to run over
all finite strings drawn from some countable collection of open
subsets of $U$ that generates the whole topology, we deduce that the
map $s \mapsto K_s$ is universally measurable, and so after
modifying it on a $\nu$-negligible set if necessary we may assume it
is Borel.

Finally, it also follows from the fact that $\tau(p(s),(s,gg')) =
\tau(p(s),(s,g'))$ for every $g' \in G_s$ and $g \in K_s$ that the
map $(s,g')\mapsto\tau(p(s),(s,g'))$ virtually factorizes through
the canonical factor
\[Z_0^S\ltimes G_\bullet \to Z_0^S\ltimes K_\bullet\backslash G_\bullet\]
to leave a map
\[\a:Z_0^S\ltimes K_\bullet\backslash G_\bullet \to Z_0^T,\]
and that this is injective away from some negligible set, since if
$\a(s_1,K_{s_1}g_1) = \a(s_2,K_{s_2}g_2)$ then $s_1 =
\xi(\a(s_1,g_1)) = \xi(\a(s_2,g_2)) = s_2$ and now
\[\tau(p(s_1),(s_1,g_1)) =\tau(p(s_1),(s_1,(g_1g_2^{-1})g_2)) =
\tau(p(s_1),(s_1,g_2)),\] which implies that $g_1g_2^{-1} \in
K_{s_1}$ provided $s_1$ did not lie in the negligible subset of
$Z_0^S$ on which we modified $K_\bullet$ above.  By another appeal
to Theorem~\ref{thm:meas-select} the map $\a$ has a Borel virtual
inverse, say
\[\b:Z_0^T\to Z_0^S\ltimes K_\bullet\backslash
G_\bullet.\] Now $\phi:= \b\circ\zeta_0^T$ is a coordinatization of
$Z_0^T$ for which the desired diagram is commutative.

\quad\textbf{(2)}\quad From the construction of $\b$ and another
measurable selection there is a measurable map $u:Y \to
K_\bullet\backslash G_\bullet$ such that
\[\phi(y,g) =
\phi\big(\tau_0\big((y,1_{G_{\zeta_0^S(y)}}),(\zeta_0^S(y),g)\big)\big)=(\zeta_0^S(y),u(y)g);\]
composing $u$ with a measurable selector $Y\ltimes
K_\bullet\backslash G_\bullet \to Y\ltimes G_\bullet$ gives a
measurable map $b:Y \to G_\bullet$ such that $\phi(y,g) =
\big(\zeta_0^S(y),K_{\zeta_0^S(y)}b(y)g\big)$ almost everywhere.

\quad\textbf{(3)}\quad Since $\phi$ is $T$-invariant we have
\begin{eqnarray*}
&&\big(\zeta_0^S(y),K_{\zeta_0^S(y)}b(S^\g
y)\rho(\g,y)1_{G_{\zeta_0^S(y)}}\big)\\ &&= \big(\zeta_0^S(S^\g
y),K_{\zeta_0^S(S^\g y)}b(S^\g y)\rho(\g,y)1_{G_{\zeta_0^S(S^\g y)}}\big)\\
&&=
\phi(S^\g y,\rho(\g,y)1_{G_{\zeta_0^S(y)}})\\
&&= \phi(T(y,g)) = \phi(y,g) =
\big(\zeta_0^S(y),K_{\zeta_0^S(y)}b(y)1_{G_{\zeta_0^S(y)}}\big)
\end{eqnarray*}
$\nu$-almost surely for any $\g \in \G$, from which the required
cohomology condition follows at once.

\quad\textbf{(4)}\quad Suppose that $K'_\bullet$ and $b'$ have the
asserted properties, and let $U$ be the overall repository for our
compact group data. Then the map
\[Y\ltimes G_\bullet\to \Clos\,U:(y,g)\mapsto K'_{\zeta_0^S(y)} b'(y)g\]
is a measurable map into a standard Borel space, and it is
$\mu$-almost surely $T$-invariant from the property that $b'(S^\g
y)\rho(\g,y)b'(y)$ lies almost surely in $K'_{\zeta_0^S(y)}$.
Therefore, because $\phi$ coordinatizes the $T$-invariant factor, we
know that there is a Borel map $f:Z_0^S\ltimes (K_\bullet\backslash
G_\bullet) \to \rm{Clos}\,U$ such that
\[K'_{\zeta_0^S(y)} b'(y)g = f\big(\zeta_0^S(y),K_{\zeta_0^S(y)}b(y)g\big)\]
for $\mu$-almost every $(y,g)$.  In particular, it follows that for
$\nu$-almost every $y$, it is the case that for
$m_{G_{\zeta_0^S(y)}}$-almost every $g \in G_{\zeta_0^S(y)}$ and
$m_{b(y)^{-1} K_{\zeta_0^S(y)}b(y)}$-almost every $h \in b(y)^{-1}
K_{\zeta_0^S(y)}b(y)$ we have
\[f\big(\zeta_0^S(y),K_{\zeta_0^S(y)}b(y)g\big) = f\big(\zeta_0^S(y),K_{\zeta_0^S(y)}b(y)hg\big),\]
and hence
\[K'_{\zeta_0^S(y)} b'(y)g = K'_{\zeta_0^S(y)} b'(y)hg\quad\Rightarrow\quad K'_{\zeta_0^S(y)} b'(y) = K'_{\zeta_0^S(y)} b'(y)h.\]
Since $(g,h)$ were chosen arbitrarily from a Haar-conegligible
subset of $G_{\zeta_0^S(y)}\times b(y)^{-1} K_{\zeta_0^S(y)}b(y)$,
it follows that for $\mu$-almost every $y$ we have
\[b(y)b'(y)^{-1}\cdot K'_{\zeta_0^S(y)}\cdot (b(y)b'(y)^{-1})^{-1} \geq K_{\zeta_0^S(y)}.\]
This tells us that the set
\[\{(s,c)\in Z_0^S\ltimes G_\bullet:\ c\cdot K'_s\cdot c^{-1} \geq K_s\}\]
(which is Borel for Borel versions of the measurable assignments
$K_\bullet$ and $K'_\bullet$) has nonempty fibre above
$(\zeta_0^S)_\#\nu$-almost every $s \in Z_0^S$, and so letting $c$
be a measurable selector for this set completes the argument.

\quad\textbf{(5)}\quad Finally, observe that for almost every
$(s,K_sg')$ that parameterizes a $T$-ergodic component we have
$\phi(y,g) = (s,K_sg')$ if and only if $\zeta_0^S(y) = s$ and
$K_sb(y)g = K_s g'$, hence if and only if $g \in b(y)^{-1}K_s g'$.
From this the last conclusion follows at once. \qed

\begin{dfn}[Mackey data]
We refer to the measurable group data $K_s$ given by the above
theorem as \textbf{Mackey group data} of $\rho$ over $(Y,\nu,S)$,
and to the section $b$ as a \textbf{Mackey section} (note that in
general the Mackey group data is not unique, but is so up to
$S$-invariant conjugacy, by part (4) of the theorem).
\end{dfn}

From the above result for extensions by group data we can easily
generalize to extensions by homogeneous space data.

\begin{cor}\label{cor:homo-nonergMackey}
Suppose that $(X,\mu,T) = (Y,\nu,S)\ltimes
(G_\bullet/H_\bullet,m_{G_\bullet/H_\bullet},\rho)$, $\zeta_0^S:Y\to
Z_0^S$ a coordinatization of the base isotropy factor and $P:Z_0^S
\stackrel{\rm{p}}{\longrightarrow} Y$ a version of the
disintegration of $\nu$ over $\zeta_0^S$. Then there are subgroup
data $K_\bullet \leq G_\bullet$ and a cocycle-section $b:Y\to
G_\bullet$ such that the factor map
\[\phi:X \to Z_0^S\ltimes (K_\bullet\backslash G_\bullet/H_\bullet): (y,gH_{\zeta_0^S(y)})\mapsto (\zeta_0^S(y),K_{\zeta_0^S(y)}b(y)gH_{\zeta_0^S(y)})\]
is a coordinatization of the isotropy factor $\zeta_0^T:X\to Z_0^T$,
and the probability kernel
\[(s,K_sg'H_s)\stackrel{\rm{p}}{\mapsto} P(s,\,\cdot\,)\ltimes m_{b(\bullet)^{-1}K_sg'H_s/H_s}\]
is a version of the disintegration of $\mu$ over $\phi$, where for
any subset $S \subseteq G_s$ we write $S/H_s := \{gH_s:\ g\in S\}$.
\end{cor}

\textbf{Proof}\quad Let $(X',\mu',T') := (Y,\nu,S)\ltimes
(G_\bullet,\rho)$ and $\pi:(X',\mu',T') \to (X,\mu,T)$ the covering
factor map, and now let $K_\bullet$, $b$ and $\phi'$ be given by
applying Theorem~\ref{thm:nonergMackey} to the canonical factor map
$\pi':(X',\mu',T')\to (Y,\nu,S)$. Let $\phi$ be the map given by the
above formula.

Since $\phi$ is manifestly $T$-invariant we need only show that it
recoordinatizes the whole of $Z_0^T$, and not a properly smaller
factor. If $f \in L^\infty(\mu)$ is $T$-invariant then $f\circ\pi$
is $T'$-invariant, and so factorizes through the map $\phi'$ given
by Theorem~\ref{thm:nonergMackey}. It follows that $f$ virtually
virtually factorizes through $\phi$. \qed

\subsection{More general lifted measures on homogeneous space extensions}

Theorem~\ref{thm:nonergMackey} describes the components of the
ergodic decomposition of $\mu\ltimes m_{G_\bullet}$ under $S\ltimes
\rho$, but in fact the same ideas can be used to describe \emph{all}
relatively $(S\ltimes \rho)$-ergodic lifts of $\nu$. This stronger
result, and its corollary for extensions by homogeneous space data,
will be important for both the inverse theory to be developed in the
next section and the study of automorphisms of extensions in the
section after that.

We first set up some simple enabling results concerning the
collection of lifts of a given probability measure on an extension
by homogeneous space data to measures on a covering extension by
group data.

\begin{dfn}[Lift topology]
Suppose that $S:\G \curvearrowright (Y,\nu)$, $G_\bullet/H_\bullet$
is $S$-invariant measurable homogeneous space data with repository
$U$, $\rho:\G\times Y\to G_\bullet$ is a cocycle-section and $\mu$
is an $(S\ltimes \rho)$-invariant probability measure on $Y\ltimes
G_\bullet/H_\bullet$ that lifts $\mu$, and let $Q$ be the set of all
further lifts of $\mu$ to $(S\ltimes \rho)$-invariant probability
measures on $Y\ltimes G_\bullet$.  We define the \textbf{lift
topology} on $Q$ as the weakest topology with respect to which the
evaluation functionals
\[\mu' \mapsto \int_{Y\ltimes G_\bullet}1_A(y)\cdot f(g)\,\mu'(\d(y,g))\]
are continuous for all $A \in \S_Y$ and continuous functions
$f:U\to\bbR$.
\end{dfn}

The following is now routine.

\begin{lem}\label{lem:lifts}
Under the lift topology $Q$ is a nonempty compact convex set.
\end{lem}

\textbf{Proof}\quad It is easy to witness one member of $Q$ (and so
see that it is nonempty):
\[\mu' := \int_{Y\ltimes G_\bullet/H_\bullet} m_{gH_y}\,\mu(\d(y,gH_y)).\]
Convexity is obvious, so we need only verify compactness. However,
$Q$ is a closed subset of the larger convex set $Q_0$ containing all
lifts of $\mu$ to $Y\times U$, and this set is easily seen to be a
closed subset of the unit ball of $L^\infty(\mu;\cal{M}(U))$ (where
$\cal{M}(U)$ is the Banach space of signed measures on $U$) in the
weak$^\ast$ topology.  Now the Banach-Alaoglu Theorem tells us that
this larger set $Q_0$ is compact and the proof is complete. \qed

\textbf{Remark}\quad Of course, it is easy to construct examples of
general extensions of Borel actions $(X,T) \to (Y,S)$ for which a
given invariant probability measure on $Y$ has no invariant
extension to $X$ (indeed, with $Y = \{\ast\}$ any Borel action on
$X$ with no invariant probability measure will do).  The
nonemptyness assertion of the above lemma, though simple, is very
much a consequence of the isometric structure of the extensions in
question. \fin

\textbf{Remark}\quad An alternative route to topologizing $Q$ (to be
found, for example, in Furstenberg's paper~\cite{Fur77} and
Glasner's book~\cite{Gla03}) is to choose a coordinatization of
$(Y,\nu,S)$ as a homeomorphic action on a compact space with an
invariant Borel probability measure, and then simply introduce the
usual vague topology on our convex set of lifts.  It is not hard to
see that the resulting topology is the same; we have chosen the
present approach only because it seems more intrinsic. \fin

We can now approach the main results of this subsection.

\begin{prop}\label{prop:relergmeasures}
Suppose that $S:\G\curvearrowright (Y,\nu)$, that $G_\bullet$ are
$Z_0^S$-measurable group data and $\rho:\G\times Y\to G_\bullet$ is
a cocycle-section over $S$ and that $X$ is the space $Y\ltimes
G_\bullet$ but equipped with some unknown $(S\ltimes
\rho)$-invariant and relatively ergodic lift $\mu$ of $\nu$. Then
there are subgroup data $K_\bullet \leq G_\bullet$ and a section
$b:Y\to G_\bullet$ such that $\mu = \nu\ltimes
m_{b(\bullet)^{-1}K_\bullet}$.
\end{prop}

\textbf{Remark}\quad Once again, the case with $S$ ergodic is
classical: it can be found as Theorem 3.26 in Glasner~\cite{Gla03}.
\fin

\textbf{Proof}\quad Let $T:= S\ltimes \rho$ and consider again the
extension of the Borel system $T:\G\curvearrowright X$ given by
\[\tau_0: X \times_{\{\zeta_0^S\circ\pi=\theta\}} (Z_0^S\ltimes G_\bullet)\to X.\]
In addition, let $P:Z_0^S \stackrel{\rm{p}}{\longrightarrow} X$ be a
version of the disintegration of the unknown lift $\mu$ over
$\zeta_0^S\circ\pi$. Now let
\[K'_s := \{g \in G_s:\ \tau_0(P(s,\,\cdot\,),(s,g))=P(s,\,\cdot\,)\};\]
this is a closed subgroup of $G_s$ that is universally measurable in
$s$ by just the same argument as in part (1) of the proof of
Theorem~\ref{thm:nonergMackey}, and is such that
\[\tau_0(P(s,\,\cdot\,),(s,gg')) = \tau_0(P(s,\,\cdot\,),(s,g'))\]
whenever $g' \in G_s$ and $g \in K'_s$.  Adjusting $K'_\bullet$ on a
negligible subset of $Y$ so that it is Borel, we still obtain that
the composed kernel $(s,g) \stackrel{\rm{p}}{\mapsto}
\tau_0(P(s,\,\cdot\,),(s,g))$ virtually factorizes through the
canonical factor $Z_0^S\ltimes G_\bullet \to Z_0^S\ltimes
K'_\bullet\backslash G_\bullet$ to a kernel $P':Z_0^S\ltimes
K'_\bullet\backslash G_\bullet \stackrel{\rm{p}}{\to} X$ such that
$P'((s,K'_s),\,\cdot\,) = P(s,\,\cdot\,)$.

Now define another kernel $P'':Z_0^S \stackrel{\rm{p}}{\to} X$ by
\[P''(s,A) := \int_{K_s\backslash G_s} P'((s,K'_sg),A)\,\d(K'_s g),\quad A\in\S_X.\]
This is an `averaged out' version of $P'$. It is also clearly
$S$-invariant, satisfies $\pi_\# P''_\#{\zeta_0^S}_\#\nu = \nu$, and
now also satisfies
\[\tau_0(P''(s,\,\cdot\,),(s,g)) = P''(s,\,\cdot\,)\quad\quad\hbox{for all}\ g \in G_s;\]
hence $P''_\#{\zeta_0^S}_\#\nu$ must simply be equal to $\nu\ltimes
m_{G_\bullet}$.

Now, the measure $P'((s,K'_sg),\,\cdot\,)$ is $T$-ergodic for almost
every $(s,g)$ (since it is a fibrewise right-translate of
$P(s,\,\cdot\,)$ by some fixed element of $G_s$), and so it follows
that the integral
\[\nu\ltimes m_{G_\bullet} = P''_\#{\zeta_0^S}_\#\nu = \int_{Z_0^S}\int_{K'_s\backslash G_s} P'((s,K'_sg),\,\cdot\,)\,\d(K'_s g)\,\nu(\d s)\]
is a version of the ergodic decomposition of $\nu\ltimes
m_{G_\bullet}$.  By Theorem~\ref{thm:nonergMackey} $K'_s$ must be a
version of the Mackey group data for $\rho$ over $(Y,\nu,S)$, and
hence by part (3) of that theorem it follows that there is a section
$b':Y\to G_\bullet$ such that $\rho'(\g,y):= b'(S^\g y)\cdot
\rho(\g,y)\cdot b'(y)^{-1}$ takes values in $K'_{\zeta_0^S(y)}$
almost surely.

Now we finish the proof simply by applying the fibrewise
recoordinatizing isomorphism $\psi:(y,g)\mapsto (y,b'(y)g)$ from our
original system $(X,\mu,T)$ to the system
$(X,\psi_\#\mu,S\ltimes\rho')$. This map $\psi$ must carry the
$(S\ltimes \rho)$-ergodic decomposition of $\nu\ltimes
m_{G_\bullet}$ to its $(S\ltimes \rho')$-ergodic decomposition, and
hence each of the measures $P'((s,K'_sg),\,\cdot\,)$ to ergodic
measures supported on the disjoint sets $\{s\}\times K'_s g'$. Since
$\rho'$ almost surely takes values in $K_s$ and these components
must integrate up to $\nu\ltimes m_{G_\bullet}$, it follows that
$\psi_\#P_\#{\zeta_0^S}_\#\nu = \nu\ltimes m_{K'_\bullet
g(\bullet)}$ for some $g:Z_0^S \to G_\bullet$, and now applying
$\psi^{-1}$ to this equation and replacing $K_\bullet :=
g(\bullet)^{-1}K'_\bullet g(\bullet)$ and $b := b'\cdot g$ gives the
result. \qed

Once again, the result for extensions by group data implies a
version for extensions by homogeneous space data through lifting to
a covering group extension and then descending again, just as for
Corollary~\ref{cor:homo-nonergMackey}.  The proof is essentially the
same, and so we omit it here.

\begin{cor}
Suppose that $S:\G\curvearrowright (Y,\nu)$, $H_\bullet \leq
G_\bullet$ are $Z_0^S$-measurable group data and $\rho:\G\times Y\to
G_\bullet$ is a cocycle-section over $S$ and $X$ is the space
$Y\ltimes G_\bullet/H_\bullet$ but equipped with some unknown
$(S\ltimes \rho)$-invariant and relatively ergodic lift $\mu$ of
$\nu$. Then there are subgroup data $K_\bullet \leq G_\bullet$ on
$Z_0^S$ and a section $b:Y\to G_\bullet$ such that $\mu = \nu\ltimes
m_{b(\bullet)^{-1}K_\bullet H_\bullet/H_\bullet}$. \qed
\end{cor}

Arguing exactly as in the classical case of an ergodic base system
by replacing some given group data $G_\bullet$ with the Mackey group
data $K_\bullet$ and recoordinatizing (see Corollary 3.27 in
Glasner~\cite{Gla03}), we obtain the following corollary.

\begin{cor}\label{cor:homo-strong-nonergMackey}
Given a $\G$-system $\bfY = (Y,\nu,S)$, measurable $S$-invariant
homogeneous space data $G_\bullet/K_\bullet$ over $Y$ and a
cocycle-section $\rho:\G\times Y\to G_\bullet$, and defining $X:=
Y\ltimes G_\bullet/K_\bullet$ and $T:= S\ltimes\rho$, any $(S\ltimes
\rho)$-relatively ergodic lift $\mu$ of $\nu$ admits a
recoordinatization
\begin{center}
$\phantom{i}$\xymatrix{
(X,\mu,T)\ar[dr]_{\rm{canonical}}\ar@{<->}[rr]^-\cong
& & \bfY\ltimes(G'_\bullet/H'_\bullet,m_{G'_\bullet/H'_\bullet},\rho')\ar[dl]^{\rm{canonical}}\\
& \bfY,}
\end{center}
so that the implicit covering group extension $\bfY\ltimes
(G'_\bullet,m_{G'_\bullet},\rho')\to\bfY$ is also relatively
ergodic. \qed
\end{cor}

During the development of the inverse theory of the next section we
will use the preceding results in conjunction with the following
elementary lemma.

\begin{lem}\label{lem:joinstillisometric}
Suppose that $\pi:(X,\mu,T)\to(Y,\nu,S)$ is a relatively ergodic
extension of a not-necessarily ergodic system, and that
\[(X,\mu,T)\stackrel{\pi_{(n)}}{\longrightarrow}(Z_{(n)},\mu_{(n)},T_{(n)})\stackrel{\xi_{(n)}}{\longrightarrow}
(Y,\nu,S)\] for $n=1,2,\ldots$ is a sequence of intermediate
extensions that are all coordinatizable as extensions by
homogeneous-space data. Then the resulting joint extension
\[\big(Z_{(1)}\times Z_{(2)} \times \cdots,(\pi_{(1)}\vee\pi_{(2)}\vee\cdots)_\#\mu,T|_{\pi_{(1)}\vee\pi_{(2)}\vee\cdots}\big) \to
(Y,\nu,S)\] is also coordinatizable as an extension by
homogeneous-space data.
\end{lem}

\textbf{Remark}\quad This is a straightforward extension of Lemma
8.4 in Furstenberg~\cite{Fur77}. \fin

\textbf{Proof}\quad We know that for each $n \geq 1$ there are some
$Z_0^S$-measurable homogeneous-space data
$G_{n,\bullet}/H_{n,\bullet}$, a cocycle-section $\rho_n:\G\times
Y\to G_{n,\bullet}$ and a coordinatizing isomorphism
\[\a_n:(Z_{(n)},\mu_{(n)},T_{(n)}) \stackrel{\cong}{\longrightarrow} (Y,\nu,S)\ltimes (G_{n,\bullet}/H_{n,\bullet},m_{G_{n,\bullet}/H_{n,\bullet}}\rho_n).\]
Each $\a_n\circ\pi_{(n)}$ is a factor map of $X$ that gives a
recoordinatization of the factor associated to $\pi_{(n)}$ and takes
the form $(\pi,\theta_n)$ for a suitable map $\theta_n:X\to
G_{n,\bullet}/H_{n,\bullet}$. Now we simply set $G_\bullet :=
\prod_{n\geq 1}G_{n,\bullet}$, $H_\bullet := \prod_{n\geq
1}H_{n,\bullet}$, $\rho := (\rho_n)_{n\geq 1}$ and $\theta :=
(\theta_n)_{n\geq 1}$; it is clear that the map $(\pi,\theta)$ now
gives the desired recoordinatization of $(Z_{(1)}\times Z_{(2)}
\times
\cdots,(\pi_{(1)}\vee\pi_{(2)}\vee\cdots)_\#\mu,T|_{\pi_{(1)}\vee\pi_{(2)}\vee\cdots})$
as an extension of $Y$ by homogeneous-space data. \qed

\section{Relative weak non-mixing and isometric extensions of non-ergodic
systems}\label{sec:FZ}

In this section we shall recount the main results of our non-ergodic
version of the Furstenberg-Zimmer inverse theory. Although it seems
that these non-ergodic analogs do not formally follow from their
ergodic predecessors, their proofs largely follow the original
arguments of Furstenberg and Zimmer, with a few judicious
invocations of measurable selectors along the way.  For this reason
our presentation here, as in the preceding section, will be quite
terse. The original papers~\cite{Zim76.1,Zim76.2} and~\cite{Fur77}
remain clear and thorough references for the classical results in
the presence of ergodicity, and we direct the reader to these for
many of the original ideas.

The theory developed by Furstenberg and Zimmer considers an ergodic
extension of an ergodic system $\bfY = (Y,\nu,S)$. Given such an
extension $\pi:\bfX\to\bfY$, this theory gives an account of the
possible failure of ergodicity of the relatively independent
self-joining $\bfX\times_\pi\bfX$ (that is, of the `relative weak
mixing' of $\bfX$ over $\pi$): it turns out that this occurs if and
only if the extension contains a nontrivial subextension that can be
coordinatized as a homogeneous skew-product.  It is this result that
we shall presently extend by dropping the assumption that $S$ be
ergodic, and by working instead with extension by (possibly
variable) homogeneous space data.  Note, however, that we will
continue to assume \emph{relative} ergodicity of $\bfX\to\bfY$: if
this fails then the arguments that follow derail quite quickly, and
the best account of the structure of the extension that can be given
in this case seems to result from simply considering the relatively
invariant subextension first, and then working with the remaining
(necessarily relatively ergodic) extension over that.

In fact, here as in the ergodic setting just a little extra work
will show that once the failure of relative weak mixing is
understood in this way, the same structures account for the
non-ergodicity of other relatively independent joinings: given two
extensions $\pi_i:\bfX_i = (X_i,\mu_i,T_i) \to \bfY_i$, $i=1,2$, a
joining $\nu$ of $\bfY_1$ and $\bfY_2$ and a lift $\mu$ of $\nu$ to
a joining of $\bfX_1$ and $\bfX_2$ under which the copies of these
two factors are relatively independent over the copies of $\bfY_1$
and $\bfY_2$, then this larger joining $\mu$ can fail to be
$(T_1\times T_2)$-ergodic only if each of the extensions
$\bfX_i\to\bfY_i$ contains a nontrivial subextension that is
coordinatizable as a homogeneous skew-product, and the homogeneous
spaces and cocycles of these skew-product are suitably related to
each other.  The result above simply corresponds to the case $\bfX_1
= \bfX_2$, $\bfY_1 = \bfY_2$ and $\l$ the diagonal self-joining of
the smaller system.

\subsection{Generalized eigenfunctions and finite-rank modules}

Key to the reduction from the failure of relative weak mixing to
nontrivial extensions by homogeneous space data are the notions of
finite rank modules and isometric extensions. These definitions are
taken almost unchanged from the papers of Furstenberg~\cite{Fur77}
and Zimmer~\cite{Zim76.2}.

\begin{dfn}[Modules over factors and their rank]
If $\pi:(X,\mu)\to (Y,\nu)$ is an extension of standard Borel
probability spaces and $\frM$ is a closed subspace of $L^2(\mu)$,
then we shall refer to $\frM$ as a \textbf{$\pi$-module} if
$(h\circ\pi)\cdot f \in \frM$ whenever $f \in \frM$ and $h \in
L^\infty(\nu)$.

If $\int^\oplus_Y\frH_y\,\nu(\d y)$ is the direct integral
decomposition of $L^2(\mu)$ over $\pi$, then subordinate to this we
may form the direct integral decomposition $\int^\oplus_Y
\frM_y\,\nu(\d y)$ of $\frM$ over $\pi$; and now we shall write that
$\frM$ has \textbf{rank $r$ over $\pi$} if $r \in
\{1,2,\ldots,\infty\}$ is minimal such that $\dim \frM_y \leq r$ for
$\nu$-almost every $y$. If $r < \infty$ we shall write that $\frM$
has \textbf{finite $\pi$-rank}.
\end{dfn}

\begin{dfn}[Isometric extension]\label{dfn:isometric}
A system extension $\pi:\bfX\to\bfY$ is \textbf{isometric} if
$L^2(\mu)$ is generated as $\pi$-module by its finite rank
$T$-invariant $\pi$-submodules.
\end{dfn}

The following first step towards representing finite rank modules
will be crucial.

\begin{lem}[Orthonormal basis for a module]\label{lem:ortho}
Suppose that $\frM$ is a rank-$r$ $\pi$-module for some $r <
\infty$. Then there is a tuple $\phi_1$, $\phi_2$, \ldots, $\phi_r$
of functions in $\frM$, none of them vanishing everywhere, and a
measurable function $R:Y \to \{1,2,\dots,r\}$ such that
$\sfE_\mu(\phi_i\cdot\phi_j\,|\,\pi) = \delta_{i,j}1_{\{i,j \leq
R\}}$ and \[(L^\infty(\nu)\circ\pi)\cdot \phi_1 + \cdots +
(L^\infty(\nu)\circ\pi)\cdot \phi_r\] is $L^2$-dense in $\frM$.
\end{lem}

\textbf{Proof}\quad Let $R(y) := \dim \frM_y$ and decompose $Y$ as
$\bigcup_{s \leq r}\{y:\ R(y) = s\}$.  It is easy to check that $R$
must be measurable, and the existence of suitable $\phi_i$ now
follows just as in the classical case by considering the cells
$\{y:\ R(y) = s\}$ separately: see Lemma 9.4 of
Glasner~\cite{Gla03}. It is clear that none of the $\phi_i$ can
vanish everywhere, else the module $\frM$ would actually have rank
at most $r-1$, contradicting our assumptions. \qed

\begin{dfn}
We refer to the function $R$ above as the \textbf{local rank
function} of the module $\frM$. \fin
\end{dfn}

The following lemma is trivial, but it will prove convenient to be
able to call on it explicitly.

\begin{lem}
If $(X,\mu)\stackrel{\pi}{\to}(Z,\theta)\stackrel{\pi'}{\to}(Y,\nu)$
is a tower of probability-preserving maps and $\frM\leq L^2(\mu)$ is
a $(\pi'\circ\pi)$-module of rank $r$, then
$\overline{L^\infty(\theta)\cdot \frM}$ is a $\pi$-module of rank at
most $r$. \qed
\end{lem}

\subsection{The non-ergodic Furstenberg-Zimmer inverse
theorem}\label{subs:FurZim}

We now present our non-ergodic extension of the Furstenberg-Zimmer
inverse theorem. This follows naturally from two separate
propositions.

\begin{prop}[From relative non-weak-mixing to finite rank modules]\label{prop:non-w.m.-to-mods} If
$\pi:\bfX\to \bfY$ is relatively ergodic but not relatively weakly
mixing, then $L^2(\mu)$ contains a nontrivial $T$-invariant
finite-rank $\pi$-submodule.
\end{prop}

\textbf{Proof}\quad This first proposition is proved just as in the
ergodic case, so we shall only sketch its proof, referring the
reader to Chapter 9 in Glasner~\cite{Gla03} or Section 7 of
Furstenberg~\cite{Fur77} for a more careful treatment. Form the
relatively independent self-product $\bfX\times_\pi\bfX$ with its
natural coordinate projections $\pi_1,\pi_2:\bfX\times_\pi\bfX \to
\bfX$. Note that by construction these maps quotient to give a
single factor copy of $\bfY$ through $\pi$, and so up to
$(\mu\otimes_\pi\mu)$-almost-everywhere equality of functions we
have $L^2(\nu)\circ\pi\circ\pi_1 = L^2(\nu)\circ\pi\circ\pi_2$. The
key idea is to choose a $(T\times T)$-invariant function $H$ that
lies in $L^2(\mu\otimes_\pi\mu)\backslash
(L^2(\nu)\circ\pi\circ\pi_1)$, and then define from it a bounded
operator $A$ on $L^2(\mu)$ by
\[A\psi := \sfE_{\mu\otimes_\pi\mu}\big((\psi\circ\pi_2)\cdot H\,\big|\,\pi_1\big).\]
It is easy to see that this cannot be the identity if $H \not\in
L^2(\nu)\circ\pi\circ\pi_1$. Replacing $H(x,x')$ by either $H(x,x')
+ \overline{H(x',x)}$ or $H(x,x') - \rm{i}\overline{H(x',x)}$ if
necessary, we may also arrange that $A$ be self-adjoint.  In
defining this $A$ we have produced a `relative Hilbert-Schmidt
operator', acting as a Hilbert-Schmidt operator separately on each
fibre $\frH_y$ of the Hilbert space direct integral decomposition
$L^2(\mu) = \int_Y^\oplus \frH_y\,\nu(\d y)$. It now admits a
spectral decomposition into finite-dimensional eigenspaces in each
of these fibres, relative over the factor $\pi$. In principle the
list of corresponding eigenvalues can vary with $y \in Y$, but it is
standard (see Section 9.3 of Glasner~\cite{Gla03}) that they do so
measurably and are $S$-invariant. Now we can simply make a
measurable selection of one of these non-zero eigenvalues over each
$\zeta_0^S(y) \in Z_0^S$ and associate to that eigenvalue its
corresponding finite-dimensional eigenspace $\frM_y$, to produce a
nontrivial $\pi$-submodule $\frM$ of $L^2(\mu)$ such that
$\dim\frM_y$ is $S$-invariant and almost surely finite, and now
truncating this by retaining the nontrivial space $\frM_y$ only on
the $S$-invariant set $\{y:\ \dim\frM_y < M\}$ for some sufficiently
large $M$ gives a true finite-rank module. Its $T$-invariance
follows immediately from that of $H$ and hence of $A$, completing
the proof. \qed

Our machinery of direct integrals of homogeneous spaces becomes
necessary for the second stage of the argument.

\begin{prop}[Coordinatization of isometric
extensions]\label{prop:isom-coord} If $\pi:\bfX\to \bfY$ is a
relatively ergodic extension for which $L^2(\mu)$ is generated by
its $T$-invariant finite-rank $\pi$-submodules, then there are
$Z_0^S$-measurable homogeneous space data $G_\bullet/K_\bullet$ with
fibre repository $U:= \prod_{n\geq 1}\rm{U}(n)^\bbN$ and a
cocycle-section $\rho:\G\times Y\to G_\bullet$ such that
\begin{center}
$\phantom{i}$\xymatrix{ \bfX\ar[dr]_\pi\ar@{<->}[rr]^-\cong & &
\bfY\ltimes (G_\bullet/K_\bullet,m_{G_\bullet/K_\bullet},\rho)\ar[dl]^{\rm{canonical}}\\
& \bfY. }
\end{center}
\end{prop}

\textbf{Proof}\quad Let $\frM^{(n)}$ for $n\geq 1$ be a sequence of
finite-rank $\pi$-modules generating $L^2(\mu)$; since $L^2(\mu)$ is
separable we need only countably many.

Suppose $\frM^{(n)}$ has rank $r$ and let $\phi_1$, $\phi_2$,
\ldots, $\phi_r$ be an orthonormal basis for it as guaranteed by
Lemma~\ref{lem:ortho}.  It is easy to check that $\phi_j\circ T^\g$,
$j\leq r$, also form an orthonormal basis for $\frM^{(n)}$ for each
$\g \in \G$, and so we can find a measurable cocycle $\Phi:\G\times
Y \to \rm{U}(R(\bullet))$ with values in the finite-dimensional
unitary group such that
\[\big(\phi_1|_{\pi^{-1}(S^\g y)}\circ T^\g,\ldots,\phi_{R(y)}|_{\pi^{-1}(S^\g y)}\circ T^\g\big) = \Phi(\g,y)\big(\phi_1|_{\pi^{-1}(y)},\ldots,\phi_{R(y)}|_{\pi^{-1}(y)}\big)\]
in $\frH_y\times\cdots\times \frH_y$ for $\nu$-almost every $y$
(indeed, this equation serves as the definition of $\Phi$ and
witnesses its measurability as a function of $(\g,y)$). From this it
follows that $\sum_{j \leq r}|\phi_j(x)|^2$ is $T$-invariant, and by
relative ergodicity can therefore be factorized through $\pi$. It
must therefore equal $R(\pi(x))$ almost everywhere by the relative
orthonormality of the $\phi_j$.

Letting
\[\phi:x \mapsto \frac{1}{R(\pi(x))}\big(\phi_1(x),\phi_2(x),\ldots,\phi_{R(y)}(x)\big),\]
it follows that $x\mapsto (\pi(x),\phi(x))$ is a map $X \to Y\ltimes
\rm{S}^{2R(\bullet)-1}$ (where $\rm{S}^{2r-1}$ denotes the unit
sphere in $\bbC^r$) which intertwines $T$ with $S\ltimes \Phi$.
Letting $\xi:Y\to \rm{S}^{2R(\bullet)-1}$ be a measurable selection
(in this case we could take it to be constant on the level-sets of
$R$) and then noting that $\rm{U}(R(\bullet))\curvearrowright
\rm{S}^{2R(\bullet)-1} \cong
\rm{U}(R(\bullet))/\rm{Stab}_{\rm{U}(R(\bullet))}(\xi(\bullet))$, it
follows that $(\pi,\phi)$ coordinatizes a subextension of $\bfX \to
\bfY$ as the homogeneous-space data extension $S\ltimes
\Phi\curvearrowright Y\ltimes
\rm{U}(R(\bullet))/\rm{Stab}_{\rm{U}(R(\bullet))}(\xi(\bullet))$
carrying some invariant measure.  Finally, by
Corollary~\ref{cor:homo-strong-nonergMackey} we know that this can
be adjusted to a genuine extension by homogeneous-space data
carrying the associated direct integral measure.

Writing this subextension as $\bfX \stackrel{\pi^{(n)}}{\to}
\bfZ^{(n)}\stackrel{\pi^{(n)\prime}}{\to} \bfY$, it is clear that
$L^2(\pi^{(n)}_\#\mu)\circ\pi^{(n)} \geq \frM^{(n)}$.  Hence the
target system of the factor $\pi^{(1)}\vee \pi^{(2)}\vee\ldots$
contains every $\frM^{(n)}$, and so must be equivalent to the whole
system $\bfX$. Finally, Lemma~\ref{lem:joinstillisometric} assures
us that this can still be coordinatized as an extension by
homogeneous space data, completing the proof. \qed

We should also check the converse of the preceding proposition in
the non-ergodic setting.

\begin{lem}\label{lem:FZ-direct}
A relatively ergodic extension by compact homogeneous space data
$\bfX := \bfY\ltimes
(G_\bullet/K_\bullet,m_{G_\bullet/K_\bullet},\rho)$ with canonical
factor map $\pi:\bfX\to \bfY$ is generated by its finite-rank
$\pi$-submodules.
\end{lem}

\textbf{Proof}\quad Let $U$ be the compact fibre repository for the
data $G_\bullet$, and let $\frH_y := L^2(m_{G_y/K_y})$ for $y \in
Y$; it is easy to check that this defines a measurable family of
Hilbert spaces.  Let $L:G_y \curvearrowright \frH_y$ be the left
regular unitary action.  Now for any continuous function $\psi$ on
$U$ the associated operators
\[A_y := \int_{G_y} \psi(g)L_g\,m_{G_\bullet}(\d g)\]
clearly form a measurable family in $y$ and are each a compact
operator on $\frH_y$.  Moreover, their eigenspaces are
$L_{G_y}$-invariant, and so different measurable selections of these
eigenspaces now combine to form finite-rank $\pi$-submodules of
$L^2(\mu) \cong \int_Y^\oplus \frH_y\,\nu(\d y)$ that are
$T$-invariant. Finally choosing a single sequence of continuous
mollifiers $(\psi_n)_{n\geq 1}$ on $U$ (that is, of continuous
functions on $U$ such that $m_U \llcorner \psi_n \to \delta_{1_U}$
in the vague topology as $n\to\infty$), all possible measurable
selections of eigenspaces of their associated compact operators
together generate the whole of $L^2(\mu)$, completing the proof.
\qed

Given an arbitrary relatively ergodic extension $\pi:\bfX\to \bfY$
that is not relatively weakly mixing,
Proposition~\ref{prop:non-w.m.-to-mods} guarantees that its
subextension generated by all finite-rank $\pi$-submodules is
nontrivial; and now applying Proposition~\ref{prop:isom-coord} to
this subextension immediately gives our full version of the inverse
theorem.  Combined with the `direct' result of
Lemma~\ref{lem:FZ-direct} this gives the following.

\begin{thm}[Furstenberg-Zimmer Theorem in the non-ergodic setting]\label{thm:FurZim} Suppose that the extension $\pi:\bfX\to\bfY$
is relatively ergodic.  Then there is a unique maximal subextension
$\bfX\to \bfZ\to \bfY$ that can be coordinatized as an extension of
$\bfY$ by homogeneous-space data, it equals the subextension
generated by all finite-rank $T$-invariant $\pi$-modules, and $\pi$
fails to be relatively weakly mixing for $T$ if and only if this
subextension contains $\pi$ strictly. \qed
\end{thm}

\begin{dfn}[Maximal isometric subextension]
The subextension $\bfZ\to \bfY$ given by the preceding theorem is
the \textbf{maximal isometric subextension} of $\pi:\bfX\to\bfY$.
\end{dfn}

\textbf{Remark}\quad Let us digress to locate the need for relative
ergodicity in the above arguments. This occurred during the proof of
Proposition~\ref{prop:isom-coord} when we argued that our
orthonormal basis $\phi_1$, $\phi_2$, \ldots, $\phi_r$ must have
$\sum_{j\leq r}|\phi_j(x)|^2$ measurable with respect to $\pi$ in
view of its $T$-invariance, and from this that we could synthesize
from this basis a map $X\to \rm{S}^{2R(\bullet)-1}$ which together
with $\pi$ would lead to an explicit coordinatization by homogeneous
space data.

It is clear that given an extension that is not relatively ergodic,
we can still derive some structural consequences from the failure of
relative weak mixing as follows.  A simple check shows that the
argument that converts a nontrivial invariant function on
$\bfX\times_\pi \bfX$ to a finite-rank module over $\pi$ does not
require relative weak mixing, and so we can still sensibly define
the subextension $\bfX \stackrel{\xi}{\to} \bfZ\stackrel{\a}{\to}
\bfY$ generated by all finite-rank $\pi$-modules.  Next, any
function on $X$ that is actually $T$-invariant clearly defines a
rank-$1$ such module, and so $\xi$ certainly contains the factor
$\zeta_0^T\vee \pi$.  Now, in addition, any finite-rank module
$\frM$ over the factor $\pi$ gives a finite-rank module
$\overline{L^\infty((\zeta_0^T\vee \pi)_\#\mu)\cdot \frM}$ over
$\zeta_0^T\vee \pi$.  In light of this and the inverse theory for
the relatively ergodic case, $\xi:\bfX\to\bfZ$ must actually be
contained in the maximal isometric subextension of the joint factor
map $\zeta_0^T\vee \pi$.  Letting $\bfW$ be the target system of
$\zeta_0^T\vee \pi$, the maximal isometric subextension of the
relatively ergodic extension $\zeta_0^T\vee \pi$ now \emph{is}
coordinatizable as an extension by compact homogeneous space data,
and from the results of the next section it will follow that the
resulting subextension $\bfZ\to \bfW$ is also coordinatizable by
compact homogeneous space data. It may, however, be properly
contained in the maximal isometric subextension of $\zeta_0^T\vee
\pi$, since there may be finite-rank modules over $\zeta_0^T\vee
\pi$ that cannot be obtained as above from finite-rank modules over
$\pi$.

We suspect that more can be said in general about which
subextensions of $\zeta_0^T\vee\pi$ can be obtained from finite-rank
modules over $\pi$, but we will not explore this matter further
here. Note that in several previous works, such as Furstenberg and
Katznelson's proof of multidimensional multiple recurrence
in~\cite{FurKat78} and their later applications of similar ideas to
prove other results in density Ramsey theory in~\cite{FurKat85}
and~\cite{FurKat91}, this very concrete analysis in terms of
extensions by homogeneous space data is avoided altogether.  In its
place is used a much softer property of extension called `relative
compactness', which is also a consequence of its being generated by
finite-rank submodules and turns out to be enough to enable a proof
of the relevant multiple recurrence results without the more precise
information offered by a coordinatization. We will not explore
relations with this idea further here, but refer the reader to
Furstenberg's book~\cite{Fur81} for a treatment of such arguments.
\fin

We can now extend the following definition from the ergodic-base
case.

\begin{dfn}[Distal extensions and distal towers]
Given a system extension $\pi:\bfX\to\bfY$ we note that
$\zeta_0^T\vee \pi$ coordinatizes the maximal factor of $\bfX$ that
is relatively invariant over $\pi$, and now write
$\zeta_{1/\pi}^T:\bfX\to \bfZ_1^T(\bfX/\pi)$ to denote any choice of
factor map that coordinatizes the maximal isometric subextension of
$\zeta_0^T\vee\pi$. In general we define recursively an increasing
transfinite sequence of factors $\zeta_{\eta/\pi}^T:\bfX\to
\bfZ_{\eta/\pi}^T(\bfX/\pi)$ indexed by all ordinals $\eta$ by
letting $\zeta_{\eta+1/\pi}^T:\bfX\to \bfZ_{\eta+1/\pi}^T(\bfX/\pi)$
denote any choice of factor map that coordinatizes the maximal
isometric subextension of $\zeta_{\eta/\pi}^T$ for each $\eta$, and
letting $\zeta_{\eta/\pi}^T := \bigvee_{\k <
\eta}\zeta_{\eta/\pi}^T$ when $\eta$ is a limit ordinal. Note that
for any fixed system $\bfX$ there must be some ordinal $\leq
\omega_1$ at which this tower stabilizes. If $\bfY$ is a trivial
system $(\{\ast\},\delta_\ast,\id_{\{\ast\}})$ we simplify this
notation to $\zeta_{\eta}^T:\bfX\to\bfZ_{\eta}^T$.

The extension $\pi:\bfX\to \bfY$ is \textbf{distal} if $\bigvee_\eta
\zeta_{\eta/\pi}^T \simeq \id_X$.  We refer to
$\bfZ_\eta^T(\bfX/\pi)\to\bfY$ as the \textbf{maximal $\eta$-step
distal subextension of $\pi$}, and to the totally ordered collection
of all of these subextensions as the \textbf{distal tower} of $\pi$.
\end{dfn}

As suggested at the beginning of this section, our inverse theorem
can quite easily be extended to account for relative non-ergodicity
of more general relatively independent self-joinings. The standard
proof of this in the case of ergodic base (Theorem 9.21 in
Glasner~\cite{Gla03}), which reduces this situation to that of the
self-joining treated above, does not rely on the ergodicity of the
base system and so carries over essentially unchanged.  We only
state the result here.

\begin{thm}\label{thm:rel-ind-joinings}
Suppose that $\pi_i:\bfX_i\to\bfY_i$ are relatively ergodic
extensions for $i=1,2,\ldots,n$ and that $\nu$ is a joining of
$\bfY_1$, $\bfY_2$, \ldots, $\bfY_n$ forming the system $\bfY =
(Y,\nu,S) := (Y_1\times Y_2\times\cdots\times Y_n,\nu,S_1\times
S_2\times \cdots\times S_n)$.  Suppose further that $\bfX =
(X,\mu,T)$ is similarly a joining of $\bfX_1$, $\bfX_2$, \ldots,
$\bfX_n$ that extends $\nu$ through the coordinatewise factor map
$\pi:\bfX\to\bfY$ assembled from the $\pi_i$, and such that under
$\mu$ the coordinate projections $\a_i:\bfX\to\bfX_i$ are relatively
independent over the tuple of further factors $\pi_i\circ\a_i$. Then
the intermediate factor map
\[\zeta_{1/\pi_1}^{T_1}\vee\zeta_{1/\pi_2}^{T_2}\vee\ldots\vee\zeta_{1/\pi_n}^{T_n}:\bfX\to\bfZ\]
whose target $\bfZ$ is a joining of the systems
$\bfZ_1^{T_i}(\bfX_i/\pi_i)$ is equivalent to
$\zeta_{1/\pi}^T:\bfX\to\bfZ_1^T(\bfX/\pi)$.  In particular, it
contains the relatively invariant extension $\bfZ_0^T\vee \bfY \to
\bfY$, which may be nontrivial. \qed
\end{thm}

We will call on this version of the inverse theorem when we come to
our applications in Section~\ref{sec:app}.

\section{Factors and automorphisms of isometric extensions}\label{sec:autos}

In this subsection we examine the possible forms of factors and
automorphisms of extensions by homogeneous space data, using as our
main tool the non-ergodic Mackey Theory of Section~\ref{sec:Mackey}.

\subsection{Some more notation}

We first need to set up some additional notation that will help us
to describe algebraic transformations between the fibres $G_y/K_y$
of our extensions.

Our first important convention is that given a Polish group $U$ and
two compact subgroups $G,G' \leq U$, we identify a continuous
isomorphism $\Phi:G\stackrel{\cong}{\longrightarrow} G'$ with its
graph $\{(g,\Phi(g)):\ g\in G\}\leq U\times U$.  In view of the
results of Subsection~\ref{subs:group}, this sets up a bijective
correspondence between continuous isomorphisms and compact subgroups
$M \in \Lat\,(U\times U)$ with the property that $M$ has first and
second projections equal to $G$ and $G'$ respectively and $M\cap
(U\times \{1_U\}) = M\cap (\{1_U\}\times U) = \{1_{U\times U}\}$. We
write $\rm{Isom}(G,G')$ for this collection of isomorphisms, and
interpret composition $\Phi\circ \Phi'$ for $\Phi \in
\rm{Isom}(G,G')$ and $\Phi'\in\rm{Isom}(G',G'')$ in the obvious way.
Note that after fixing a complete separable metric on $U$, the
resulting strong topology on $\rm{Isom}(G,G')$ is in general
strictly stronger than the Vietoris topology inherited from
$\Lat(U\times U)$, but the resulting Borel structures are the same.
Exactly similarly we can also interpret the collection
$\rm{Hom}(G,G')$ of homomorphisms $G\to G'$, $\rm{Aut}\,G =
\rm{Isom}(G,G)$ and $\rm{End}\,G = \rm{Hom}(G,G)$ as collections of
subgroups of $U\times U$.

This identification made, the necessary definition follows very
naturally.

\begin{dfn}[Homomorphisms sections]\label{dfn:isom-sections}
Suppose that $(Y,\nu)$ is a standard Borel probability space and
that $G_{i,\bullet}$, $i=1,2$, are two different measurable families
of compact group data over $Y$ with compact metrizable repositories
$U$ and $V$ respectively.  Then a \textbf{homomorphism section}
associated to this data is a map $y\mapsto \Phi_y:Y\to \Lat(U\times
V)$ that is Borel for the Vietoris measurable structure on
$\Lat(U\times V)$ and is such that $\Phi_y \in
\rm{Hom}(G_{1,y},G_{2,y})$ for $\nu$-almost every $y\in Y$. We will
sometimes denote this situation by $\Phi:Y\to
\rm{Hom}(G_{1,\bullet},G_{2,\bullet})$.

We extend this definition in the obvious way to epimorphism
sections, isomorphism sections and automorphism sections.
\end{dfn}

The benefit of formulating our notion of `isomorphism section' as
above will become clear shortly, when we use the Mackey theory in a
self-joining of a given extension to produce a measurable family of
compact subgroups of $U\times U$ that we can then immediately
re-interpret as an isomorphism cocycle in this sense.  In doing this
we will also benefit from having the following notation.

\begin{dfn}[Extensions by epimorphism sections]
If $\b:(Y_1,\nu_1) \to (Y_2,\nu_2)$ is a probability-preserving map,
$G_{i,\bullet}/H_{i,\bullet}$ is measurable compact homogeneous
space data on $Y_i$ for $i=1,2$ and $\Phi:Y\to
\rm{Hom}(G_{1,y},G_{2,\beta(y)})$ is a measurable epimorphism
section such that $\Phi_y(H_{1,y}) \subseteq H_{2,\b(y)}$ almost
surely then we write $\b\ltimes \Phi_\bullet$ for the
probability-preserving map \[Y_1\ltimes
G_{1,\bullet}/H_{1,\bullet}\to Y_2\ltimes
G_{2,\bullet}/H_{2,\bullet}: (y,gH_{1,y})\mapsto
(\b(y),\Phi_y(g)H_{2,\b(y)}).\]

More generally, if $\rho:Y\to G_{2,\b(\bullet)}$ is a section then
we will sometimes write $L_{\rho(y)}$ (resp. $R_{\rho(y)}$) for the
translation of $G_{2,\b(\bullet)}$ by left-rotation by $\rho(y)$
(resp. right-rotation by $\rho(y)$), and define
\[\b\ltimes
(L_{\rho(\bullet)}\circ\Phi_\bullet)|^{H_{1,\bullet}}_{H_{2,\b(\bullet)}}:(y,gH_{1,\bullet})\mapsto
(\b(y),\rho(y)\cdot\Phi_y(g)H_{2,\b(y)}).\] In case
$H_{i,\bullet}\equiv \{1_{G_{i,\bullet}}\}$ we simplify this to
\[\b\ltimes
(L_{\rho(\bullet)}\circ\Phi_\bullet):(y,g)\mapsto
(\b(y),\rho(y)\cdot\Phi_y(g)).\]
\end{dfn}

Note that if $(Y_1,\nu_1) = (Y_2,\nu_2)$ then $\b\ltimes
L_{\rho(\bullet)} = \b \ltimes \rho$ and $\b\ltimes
R_{\rho(\bullet)} = \b \ltimes \rho^{\rm{op}}$ in the original
notation of Section~\ref{sec:dirint}.  We introduce the above class
of fibrewise transformations, together with their new notation, to
help us describe more explicitly certain maps and joinings between
isometric extensions, in the sense of the following definition, and
to help differentiate them from those isometric extensions
themselves, among which these new maps serve as morphisms.

\begin{dfn}[Fibrewise automorphism and affine recoordinatizations]
Suppose that $\bfY = (Y,\nu,S)$ is a $\G$-system,
$G_\bullet/K_\bullet$ are $S$-invariant compact homogeneous space
data over $Y$ and $\s:\G\times Y\to G_\bullet$ is an ergodic
cocycle-section for $S$, and that
$\Phi_\bullet:Y\to\rm{Aut}(G_\bullet)$ is an $S$-invariant
automorphism section and $\rho:Y\to G_\bullet$ a section. Then the
map $R:= \id_Y\ltimes
(L_{\rho(\bullet)}\circ\Phi_\bullet)|^{K_\bullet}_{\Phi_\bullet(K_\bullet)}$
defines a recoordinatization
\begin{center}
$\phantom{i}$\xymatrix{ \bfY\ltimes
(G_\bullet/K_\bullet,m_{G_\bullet/K_\bullet},\s)\ar[dr]_{\rm{canonical}}\ar@{<->}[rr]^R
& & \bfY\ltimes
(G_\bullet/K'_\bullet,m_{G_\bullet/K'_\bullet},\s')\ar[dl]^{\rm{canonical}}\\
&\bfY }
\end{center}
with $K'_\bullet := \Phi_\bullet(K_\bullet)$ and
\[\s'(\g,y) := \rho(S^\g y)\cdot \Phi_y(\s(\g,y))\cdot \rho(y)^{-1}\]
We refer to such a recoordinatization as a \textbf{fibrewise affine
recoordinatization}. If $\rho \equiv 1_{G_\bullet}$ then it is a
\textbf{fibrewise automorphism recoordinatization}.
\end{dfn}

\textbf{Remarks}\quad\textbf{1.}\quad The condition that
$\Phi_\bullet$ be $S$-invariant is needed in order that $R$
intertwine $S\ltimes\s$ with another cocycle extension; without this
condition the new transformation will in general still involve also
a nontrivial automorphism-valued cocycle.

\quad\textbf{2.}\quad If $G_\bullet$ and $\Phi_\bullet$ are as above
then we have $\b\ltimes (L_{\rho(\bullet)}\circ\Phi_\bullet) =
\b\ltimes (\rm{Co}_{\rho(\bullet)}\circ\Phi_\bullet\circ
R_{\rho(\bullet)})$, where $\rm{Co}_{\rho(y)} \in \rm{Aut}(G_y)$ is
the inner automorphism of conjugation by $\rho(y)$. It follows that
we can always choose between expressing a fibrewise affine
recoordinatization by using fibrewise left-rotations or fibrewise
right-rotations (or even some mixture of the two!). In general we
will prefer to write fibrewise affine recoordinatizations with
rotation part acting on the left, since this is the form in which
the $S$-invariance of the automorphism part $\Phi_\bullet$ is
directly visible, but occasionally this ability to conjugate between
the two forms will give some useful flexibility in how we write a
fibrewise affine recoordinatization. \fin

Note that transformations of the form $\b\ltimes
(L_{\rho(\bullet)}\circ\Phi_\bullet)\curvearrowright U\times V$ for
compact Abelian groups $U$ and $V$, $\beta$ a rotation of $U$,
$\rho:U\to V$ and $\Phi:U\to \rm{Aut}\,V$ are already objects of
study in ergodic theory: examples arise from explicitly
coordinatizing a flow on a two-step solvmanifold (see, for example,
Auslander, Green and Hahn~\cite{AusGreHah63} or Starkov~\cite{Sta00}
for background), with rather special conditions on the cocycles
$\rho$ and $\Phi$ that result (although working explicitly with such
coordinatizations would probably be a cumbersome way to handle
solvflows).  More generally, it is possible that an enlarged theory
of extensions by compact homogeneous space data could be constructed
to treat actions lifted by affine-valued cocycle comprising both
fibrewise rotations and automorphisms. Certainly, $\bbZ^d$-actions
by affine maps on compact groups have recently begun to receive
greater attention from ergodic theorists (see Schmidt~\cite{Sch95},
in particular), and it would be interesting to know how far the
results that are now known for such systems could be `relativized'.
In this paper, however, fibrewise automorphisms and affine maps will
play a strictly auxiliary r\^ole.

\subsection{The structure theorems}

The main results of this section amount to structure theorems for
factors and automorphisms of isometric extensions, generalizing the
classical result (and also the proof) of Mentzen~\cite{Men91} in the
ergodic case (see also earlier work of Newton~\cite{New79}, and
compare with the classical argument of Veech in~\cite{Vee82} for his
condition for an extension to be coordinatizable as a compact group
skew-product).

\begin{thm}[Relative Factor Structure Theorem]\label{thm:RFST}
Suppose that $\bfY_i = (Y_i,\nu_i,S_i)$ for $i=1,2$ are$\G$-systems
and $\b:\bfY_1\to\bfY_2$ is a factor map, that
$G_{i,\bullet}/H_{i,\bullet}$ are $S_i$-invariant core-free
homogeneous space data on $Y_i$ and that $\s_i:\G\times Y_i\to
G_{i,\bullet}$ are ergodic cocycle-sections for the action $S_i$,
and let $\bfX_i = (X_i,\mu_i,T_i) := \bfY_i\ltimes
(G_{i,\bullet}/H_{i,\bullet},m_{G_{i,\bullet}/H_{i,\bullet}},\s_i)$
with canonical factor $\xi_i:\bfX_i\to \bfY_i$. Suppose further that
$\b$ admits extension to a factor map $\a:\bfX_1 \to \bfX_2$:
\begin{center}
$\phantom{i}$\xymatrix{ \bfX_1\ar[r]^\a\ar[d]_{\xi_1} &
\bfX_2\ar[d]^{\xi_2}\\ \bfY_1\ar[r]_\b & \bfY_2 }
\end{center}
Then there are an $S_1$-invariant measurable family of epimorphisms
$\Phi_\bullet:G_{1,\bullet} \to G_{2,\b(\bullet)}$ such that
$\Phi_\bullet(H_{1,\bullet}) \subseteq H_{2,\b(\bullet)}$ almost
surely and a section $\rho:Y_1\to G_{2,\b(\bullet)}$ such that $\a =
\b\ltimes (L_{\rho(\bullet)}\circ
\Phi_\bullet)|^{H_{1,\bullet}}_{H_{2,\bullet}}$, and then
\[\s_2(\g,\b(y)) = \rho(S_1^\g y)\cdot \Phi_y(\s_1(\g,y))\cdot \rho(y)^{-1}\]
for $\nu_1$-almost all $y$ for all $\g \in \G$.
\end{thm}

\textbf{Remark}\quad When $S_1$ is ergodic this tells us that any
extension of $\b$ to a factor map of $\bfX_1$ must take the form of
fibrewise application of a \emph{fixed} group automorphism and then
left-multiplication by a cocycle: this is the result of
Mentzen~\cite{Men91}. \fin

\textbf{Proof}\quad We will deduce this by considering the joining
of $\bfX_1$ and $\bfX_2$ defined by the graph of $\a$, proceeding in
two steps.

\quad\textbf{Step 1}\quad Suppose first that $H_{i,\bullet} \equiv
\{1_{G_{i,\bullet}}\}$ for $i=1,2$. Setting $\l :=
(\rm{id}_{X_1},\a)_\#\mu_1$, this is a $T_1\times T_2$-invariant
probability measure on
\[X_1\times_{\{\b\circ\xi_1 = \xi_2\}} X_2 \cong Y_1\ltimes
(G_{1,\bullet}\times G_{2,\b(\bullet)}),\] which is an extension of
$\bfY_1$ via the natural factor map. Let us denote by $\pi_i$,
$i=1,2$ the two coordinate projections $X_1\times_{\{\b\circ\xi_1 =
\xi_2\}} X_2 \to X_i$, and, slightly abusively, also the coordinate
projections $G_{1,\bullet}\times G_{2,\b(\bullet)} \to
G_{1,\bullet},G_{2,\b(\bullet)}$. Since this is a graph joining the
first coordinate projection almost surely determines the second, and
so the joined system that results is actually isomorphic to
$\bfX_1$; it follows that $\l$ is $(T_1\times T_2)$-relatively
ergodic over $\bfY_1$.

Applying Proposition~\ref{prop:relergmeasures} we obtain some
$S_1$-invariant Mackey group data $M_\bullet \leq
G_{1,\bullet}\times G_{2,\b(\bullet)}$ on $Y$ and a section
$b:Y_1\to G_{1,\bullet}\times G_{2,\b(\bullet)}$ such that $\l =
\nu\ltimes m_{b(\bullet)^{-1}M_\bullet}$. The measure $\l$ must
project onto $\mu_i$ under $\pi_i$, and so
\[{\pi_1}_\#(\nu\ltimes m_{b(\bullet)^{-1}M_\bullet}) =
\nu\ltimes m_{\pi_1(b(\bullet))^{-1}\pi_1(M_\bullet)} = \nu\ltimes
m_{G_{1,\bullet}}\] and hence $\pi_1(M_\bullet) = G_{1,\bullet}$
almost surely, and similarly $\pi_2(M_\bullet) = G_{2,\b(\bullet)}$
almost surely.

On the other hand, if $L_{1,\bullet}$ and $L_{2,\bullet}$ are the
first and second slices of $M_\bullet$, then we see that the
coordinate factors $\pi_1$ and $\pi_2$ are actually relatively
independent over the further canonical factors
\[(X_1,\mu_1,T_1)\to(Y_1,\nu_1,S_1)\ltimes
(G_{1,\bullet}/L_{1,\bullet},m_{G_{1,\bullet}/L_{1,\bullet}},\s_1)\]
and
\begin{multline*}
(Y_1,\nu_1,S_1)\otimes_{\{\b = \xi_2\}}(X_2,\mu_2,T_2)\\
\to(Y_1,\nu_1,S_1)\ltimes
(G_{2,\b(\bullet)}/L_{2,\bullet},m_{G_{2,\b(\bullet)}/L_{2,\bullet}},\s_2\circ\b).
\end{multline*}
Of course, under a graphical joining of a factor map such as $\l$
the second coordinate projection is almost surely determined by the
first, so we must have $L_{2,\bullet} \equiv
\{1_{G_{2,\b(\bullet)}}\}$ almost surely. Combined with the property
of having full projections, this shows that $M_\bullet$ is almost
surely the graph of a measurably-varying epimorphism
$\Phi_\bullet:G_{1,\bullet}\to G_{2,\b(\bullet)}$.

Also, since $M_\bullet$ has full one-dimensional projections we can
multiply $b$ by some $M_\bullet$-valued cocycle if necessary to
assume that $b(\bullet) = (1_{G_{1,\bullet}},\rho(\bullet))$. We can
now simply read off that when $\a(y,g) = (\b(y),g')$, then almost
surely
\[(g,g') = (m,\rho(y)\Phi_y(m))\]
for some $m\in G_{1,y}$. Hence we must have $m = g$, and so $g' =
\rho(y)\Phi_y(g)$: that is, $\a = \b\ltimes
(L_{\rho(\bullet)}\circ\Phi_\bullet)$, as required. Finally, with
this expression in hand it is immediate to check that the
commutative diagram relating $\a$, $\b$, $S$ and $T$ is equivalent
to the requirement that
\[\rho(S^\g(y))\cdot\Phi_y(\s_1(\g,y)) =
\s_2(\g,\b(y))\cdot\rho(y)\] for $\nu_1$-almost every $y$ for all
$\g \in \G$, which re-arranges into the equation stated.

\quad\textbf{Step 2}\quad Now consider the case of general core-free
$H_{i,\bullet}$ and let $\l$ be the graphical self-joining given
previously, so as in Step 1 the system $(X_1\times_{\{\b\circ\xi_1
=\xi_2\}}X_2,\l,T_1\times T_2)$ has a natural factor isomorphic to
$\bfY_1$ and is isomorphic to $\bfX_1$ through the first coordinate
projection, which also virtually determines the second coordinate
projection.

Now in addition let $\t{\bfX_i}\to \bfX_i$ be the implied covering
group-data extensions, and let $\t{\l}$ be any relatively ergodic
lift of $\l$ to an invariant joining of $\t{\bfX_1}$ and
$\t{\bfX_2}$. Arguing as in Step 1 now gives Mackey group data
$M_\bullet \leq G_{1,\bullet}\times G_{2,\b(\bullet)}$ and a section
$b:Y_1\to G_{1,\bullet}\times G_{2,\b(\bullet)}$ such that
\[\t{\l} =
\nu_1\ltimes m_{b(\bullet)^{-1}M_\bullet (H_{1,\bullet}\times
H_{2,\b(\bullet)})}.\] Since the cocycles $\s_i$ are ergodic, it
follows as before that $M_\bullet$ has full one-dimensional
projections. The condition that the first coordinate almost surely
determine the second under $\l$ becomes more subtle. Firstly, it
requires that the second slice $L_{2,\bullet}$ of $M_\bullet$
satisfy $L_{2,\bullet}\cdot H_{2,\b(\bullet)} = H_{2,\b(\bullet)}$;
but on the other hand Lemma~\ref{lem:groupcorrespondence} tells us
that $L_{2,\bullet} \unlhd G_{2,\b(\bullet)}$, and by the core-free
assumption $H_{2,\b(\bullet)}$ does not contain any nontrivial
normal subgroup, so in fact we must still have $L_{2,\bullet} \equiv
\{1_{G_{2,\b(\bullet)}}\}$.  It follows that $M_\bullet$ still
defines the graph of an epimorphism $\Phi_\bullet:G_{1,\bullet}\to
G_{2,\b(\bullet)}$. Secondly, this same condition on the
determination of the second coordinate requires that $\Phi_\bullet$
have a well-defined quotient between the spaces
$G_{1,\bullet}/H_{i,\bullet}$ and
$G_{2,\b(\bullet)}/H_{2,\b(\bullet)}$, and hence that
$\Phi_\bullet(H_{1,\bullet}) \subseteq H_{2,\b(\bullet)}$ almost
surely.

This leads as before to the expression $\b\ltimes
(L_{\rho(\bullet)}\circ\Phi_\bullet)|^{H_{1,\bullet}}_{H_{2,\b(\bullet)}}$
for $\a$, and the cocycle equation also follows as before from the
condition that $\t{\l}$ is $\t{T}_1\times \t{T}_2$-invariant
(equivalent to the intertwining property in the previous case). \qed

Specializing the above now gives a structure theorem for groups of
automorphisms of an extension.

\begin{thm}[Relative Automorphism Structure Theorem]\label{thm:RAST}
Suppose that $\bfY = (Y,\nu,S)$ is a $\G$-system, that
$G_\bullet/H_\bullet$ are $S$-invariant core-free homogeneous space
data on $Y$ and that $\s:\G\times Y\to G_\bullet$ is an ergodic
cocycle-section for the action $S$, and let $\bfX = (X,\mu,T) :=
\bfY\ltimes (G_\bullet/H_\bullet,m_{G_\bullet/H_\bullet},\s)$.
Suppose further that $\L$ is a discrete group and
$R:\L\curvearrowright (X,\mu)$ is another action that commutes with
$T$ and respects the canonical factor $\pi:\bfX\to \bfY$ (so it
defines an action of $\L$ by automorphisms of the extension $\pi$).
Then for each $h\in\L$ there are an $S$-invariant measurable family
of isomorphisms $\Phi_{h,\bullet}:G_\bullet \to
G_{R|^h_\pi(\bullet)}$ such that $\Phi_{h,\bullet}(H_\bullet) =
H_{{R|^h_{\pi}(\bullet)}}$ almost surely and a section $\rho_h:Y\to
G_{R|^h_{\pi}(\bullet)}$ such that \[R^h = R|^h_\pi\ltimes
(L_{\rho_h(\bullet)}\circ
\Phi_{h,\bullet})|^{H_\bullet}_{H_{R|^h_{\pi}(\bullet)}}\] for each
$h\in \L$, and then
\begin{itemize}
\item we have
\[\s(\g,R|_\pi^h(y)) = \rho_h(S^\g y)\cdot \Phi_{h,y}(\s(\g,y))\cdot \rho_h(y)^{-1}\]
for $\nu$-almost all $y$ for all $\g \in \G$ and $h \in \L$, and
\item we have
\[\Phi_{h_1h_2,y} = \Phi_{h_1,R|_{\pi}^{h_2}(y)}\circ \Phi_{h_2,y}\]
and
\[\rho_{h_1h_2}(y) = \rho_{h_1}(R|_\pi^{h_2}(y))\cdot \Phi_{h_1,R|_\pi^{h_2}(y)}(\rho_{h_2}(y))\]
for $\nu$-almost all $y$ for all $h_1,h_2 \in \L$.
\end{itemize}
\end{thm}

\textbf{Proof}\quad Consider first $R^h$ for some fixed $h\in\L$.
Treating $R^h$ as a factor map of $\bfX$ and applying
Theorem~\ref{thm:RFST} gives immediately the representation of $R^h$
as \[R|^h_\pi\ltimes (L_{\rho_h(\bullet)}\circ
\Phi_{h,\bullet})|^{H_\bullet}_{H_{R|^h_\pi(\bullet)}}\] for some
section $\rho_h$ and measurable family of continuous epimorphisms
$\Phi_{h,\bullet}$, and now the condition that $R^h$ actually be
equivalent to $\id_X$ (that is, it is not a proper factor map) gives
that $\Phi_{h,\bullet}$ is an isomorphism and
$\Phi_\bullet(H_\bullet) = H_{R|_\pi(\bullet)}$ almost surely.

This establishes the existence of $\rho_h$ and $\Phi_{h,\bullet}$,
and also the first of the two additional conclusions above.  To
deduce the second we need only compare the resulting
coordinatizations of each side of the equation $R^{h_1h_2} =
R^{h_1}\circ R^{h_2}$ defining $R$ as a $\L$-action: substituting
from above this becomes
\begin{multline*}
R|_\pi^{h_1h_2}\ltimes
(L_{\rho_{h_1h_2}(\bullet)}\circ\Phi_{h_1h_2,\bullet})|^{H_\bullet}_{H_{R|_\pi^{h_1h_2}(\bullet)}}\\
= (R|_\pi^{h_1}\circ
R|_\pi^{h_2})\ltimes(L_{\rho_{h_1}(R|_\pi^{h_2}(\bullet))}\circ\Phi_{h_1,R|_\pi^{h_2}(\bullet)}\circ
L_{\rho_{h_2}(\bullet)}\circ
\Phi_{h_2,\bullet})|^{H_\bullet}_{H_{R|_\pi^{h_1h_2}(\bullet)}},
\end{multline*}
and so we must have
\begin{eqnarray*}
&&\rho_{h_1h_2}(y)\cdot \Phi_{h_1h_2,y}(g)\cdot
H_{R|_\pi^{h_1h_2}(y)}\\
&&=\rho_{h_1h_2}(y)\cdot \Phi_{h_1h_2,y}(gH_y)\\
&&= \rho_{h_1}(R|_\pi^{h_2}(y))\cdot
\Phi_{h_1,R|_\pi^{h_2}(y)}(\rho_{h_2}(y)\cdot \Phi_{h_2,y}(gH_y))\\
&&= \rho_{h_1}(R|_\pi^{h_2}(y))\cdot
\Phi_{h_1,R|_\pi^{h_2}(y)}(\rho_{h_2}(y))\cdot
\Phi_{h_1,R|_\pi^{h_2}(y)}(\Phi_{h_2,y}(g))\cdot
H_{R|_\pi^{h_1h_2}(y)}
\end{eqnarray*}
for all $g\in G_y$ for $\nu$-almost every $y\in Y$. Since
$H_{R|_\pi^{h_1h_2}(y)}$ is core-free in $G_y$ almost surely, the
validity of this equation for all $g \in G_y$ implies the two parts
of the second additional conclusion above, completing the proof.
\qed

\textbf{Remark}\quad It should be possible to enhance the above
theorem further by allowing an arbitrary locally compact second
countable group $\L$ and imposing suitable continuity assumptions on
the assignments $h\mapsto \rho_h$ of measurable sections and
$h\mapsto \Phi_{h,\bullet}$ of measurable families of isomorphisms.
The additional arguments required seem to more fiddly than
enlightening, however, and so we leave the details to the interested
reader. \fin

Although Theorem~\ref{thm:RAST} shows that the $S$-ergodic fibre
systems above the points $\zeta_0^S(y) \in Z_0^S$ and
$\zeta_0^S(R|^h_\pi y)\in Z_0^S$ are isomorphic for all $h\in \L$
for almost every $y \in Y$, it need not follow that these fibre
systems are almost all isomorphic to a single model system. The
following simple example has long been a part of ergodic-theoretic
folklore.

\textbf{Example}\quad Let $(Y,\nu) := (\bbT^2,m_{\bbT^2})$, and form
the direct integral space $X:= Y\ltimes \bbT^2$ with constant fibre
$\bbT^2$ (so this is really just the direct product $Y\times
\bbT^2$) and measure $\mu:= m_{\bbT^2}\otimes m_{\bbT^2}$. Define
$T:\bbZ\curvearrowright X$ by $T(y,z) := (y,y+z)$, so overall
$(X,\mu)$ is the direct integral of the individual Kronecker systems
$(\bbT^2,m_{\bbT^2},R_y)$, writing $R_y$ for the rotation by
$y\in\bbT^2$.

In addition, suppose that $S\curvearrowright\bbT^2$ is any ergodic
toral automorphism.  Then it is easy to check that $S\times S$
commutes with $T$; in particular, it carries fibres of the obvious
factor map $X \to Y$ onto fibres, and so acts as an automorphism of
the fibre system $(\bbT^2,m_{\bbT^2},R_y)$ onto
$(\bbT^2,m_{\bbT^2},R_{Sy})$. However, the fibre systems
$(\bbT^2,m_{\bbT^2},R_y)$ are not almost all isomorphic for
different $y$: the map \[y \mapsto
(\bbT^2,m_{\bbT^2},R_y)/\sim_{\rm{Isomorphism}}\] is not almost
surely constant, even though it is invariant under the ergodic
transformation $S$.  This is possible because the isomorphism
equivalence relation on the space of all Kronecker systems (suitably
interpreted as pairs comprising a monothetic compact metrizable
subgroup of a suitable fixed repository and a distinguished element
for the rotation) is non-smooth. \fin

\begin{ques}
Can an example be found for which the group fibres $G_y$ themselves
are not almost all continuously isomorphic above each
$R|_\pi$-ergodic component of $\nu$?
\end{ques}

This may relate to the work of Conze and Raugy~\cite{ConRau07} on
the behaviour of measurable families of (not-necessarily compact)
groups related by measurable cocycles, but we have not been able to
answer the above as a direct corollary of their work.

The following question may also be related to the above:

\begin{ques}
Can an example be found in which for \emph{no} coordinatization of
the extension is it possible that each $\Phi_y$ can be extended from
$G_y$ to an automorphism of the whole repository group $U$? \fin
\end{ques}

The following corollaries concerning the extendability of
automorphisms will also prove useful later, and may be of some
independent interest.

\begin{cor}[Condition for lifting an automorphism to a group-data extension]
An action $R$ of $\L$ by automorphisms of $\bfY$ can be lifted to a
$\L$-action by automorphisms of an ergodic group-data extension
$\bfY\ltimes (G_\bullet,m_{G_\bullet},\s)$ if and only if for every
$h\in\L$ the cocycle $\G\times Y\to G_\bullet\times
G_{R(\bullet)}:(\g,y)\mapsto (\s(\g,y),\s(\g,R^hy))$ has relativized
Mackey group data over $\bfY$ that is the graph of an isomorphism
almost everywhere, and in this case any such extended action is of
the form $h\mapsto R^h\ltimes
(L_{\rho_h(\bullet)}\circ\Phi_{h,\bullet})$ for some families of
sections $\rho_h:Y\to G_\bullet$ and $S$-invariant cocycles
$\Phi_{h,\bullet}:Y\to\rm{Isom}(G_\bullet,G_{R^h(\bullet)})$ and
$\Phi_h$ is unique up to composition with an arbitrary $S$-invariant
inner automorphism cocycle. \qed
\end{cor}

\begin{cor}[Automorphisms can always be lifted to core-free ergodic covering group extensions]\label{cor:homogextiso}
Suppose that $\bfY$ is a $\G$-system, $G_\bullet/H_\bullet$ are
$S$-invariant core-free homogeneous space data and $\s:\G\times Y\to
G_\bullet$ is an ergodic cocycle-section. Set $\bfX := \bfY\ltimes
(G_\bullet/H_\bullet,m_{G_\bullet/H_\bullet},\s)$ and $\tilde{\bfX}
:= \bfX\ltimes (G_\bullet,m_{G_\bullet},\s)$. Then any action of a
discrete group by automorphisms of the canonical extension $\bfX \to
\bfY$ lifts to an action by automorphisms of the tower $\tilde{\bfX}
\to \bfX \to \bfY$.
\end{cor}

\textbf{Proof}\quad Let $\pi:\bfX\to\bfY$ and
$\t{\pi}:\tilde{\bfX}\to\bfY$ be the canonical factor maps and
suppose that $R$ is an automorphism of the extension
$\pi:\bfX\to\bfY$. Theorem~\ref{thm:RAST} allows us to write $R$
explicitly as $R|_\pi\ltimes
(L_{\rho(\bullet)}\circ\Phi_\bullet)|^{H_\bullet}_{H_{R|_\pi(\bullet)}}$
for some $\rho:Y\to G_{R|_\pi(\bullet)}$ and $S$-invariant
$\Phi:Y\to \rm{Isom}(G_\bullet,G_{R|_\pi(\bullet)})$ such that
$\Phi_y(H_y) = H_{R|_\pi(y)}$ almost surely; and, having done this,
we have that $\rm{graph}(\Phi_\bullet)$ is the Mackey group data of
the cocycle-section $(\g,y)\mapsto (\s(\g,y),\s(\g,Ry))$, so that
\[\s(\g,R(y)) = b(S^\g y)\cdot\Phi_y(\s(\g,y))\cdot b(y)^{-1}\]
almost surely. This equation immediately tells us that we can lift
$R$ to the transformation $\tilde{R} := R|_\pi\ltimes
(L_{\rho(\bullet)}\circ\Phi_\bullet)$ on $\tilde{X}$, and that this
still commutes with $\tilde{T} = S\ltimes\s$. Given a whole
$\L$-action of automorphisms $R^h$, applying this argument to each
$h\in\L$ individually and considering the consistency equations
promised by Theorem~\ref{thm:RAST} shows that the lifted maps still
define a $\L$-action, and hence completes the proof. \qed

\section{Applications}\label{sec:app}

In this section we offer two closely-related applications of the
theory developed above.

We first study the possible joint distribution of the isotropy
factors $\zeta_0^{T_i}$ corresponding to three commuting
transformations $T_1$, $T_2$ and $T_3$.  This will require some
quite careful analysis in terms of Mackey group data, cocycles, and
representations given by the Relative Automorphism Structure
Theorem.  We will then show that this analysis can also be brought
to bear on a detailed description of characteristic factors of the
double nonconventional ergodic averages associated to a pair of
commuting transformations (see, for example,~\cite{Aus--nonconv} and
the references listed there).

Throughout this section we specialize to the setting of $\G :=
\bbZ^d$, and will write $\bf{e}_1$, $\bf{e}_2$, \ldots, $\bf{e}_d$
for its standard basis.

\subsection{Application to joint distributions of isotropy
factors}\label{subs:isotropy}

For a generic $\bbZ^d$-action on a fixed atomless $(X,\mu)$ the
isotropy factors $\zeta_0^{T^{\ \uhr\L}}:\bfX\to \bfZ_0^{T^{\
\uhr\L}}$ corresponding to subgroups $\L \leq \bbZ^d$ are all
trivial (indeed, it is a classical result that a generic such action
is totally weakly mixing). However, if they are not all trivial then
they generate a sublattice of the lattice of all factors of
$(X,\mu,T)$ that can exhibit some quite rich structure.

Letting $T_i := T^{\bf{e}_i}$, we will here consider only the
further sublattice generated by the isotropy factors $\zeta_0^{T_I}
:= \zeta_0^{T_{i_1},T_{i_2},\ldots, T_{i_r}}$ corresponding to the
possible choices of subset $I := \{i_1, i_2,\ldots,i_r\} \subseteq
[d]$, where $[d] := \{1,2,\ldots,d\}$.

Clearly in general the action of each $T_j$ for $j \in [d]\setminus
I$ on the sets of $\S_X^{T_I}$ can still be quite arbitrary, and so
we cannot hope to say anything about the structure of each isotropy
factor as a system in its own right. Instead we will focus on their
\emph{joint} distribution within the original system.

\textbf{Example}\quad Let $(X,\mu,T_1,T_2)$ be the $\bbZ^2$-system
$(\bbT^2,\rm{Haar},R_{(\a,0)},R_{(0,\a)})$, where $R_q$ denotes the
rotation of the compact Abelian group $\bbT^2$ by an element $q \in
\bbT^2$ and we choose $\a \in \bbT$ irrational.  In this case we
have natural coordinatizations
\[\zeta_0^{T_i}:X \to \bbT:(t_1,t_2)\to t_{3-i},\]
and similarly, since $T_1T_2 = R_{(\a,\a)}$,
\[\zeta_0^{T_1T_2}:X\to \bbT:(t_1,t_2)\to t_1 - t_2.\]
It follows that in this example any two of $\zeta_0^{T_1}$,
$\zeta_0^{T_2}$ and $\zeta_0^{T_1T_2}$ are independent, but also
that any two of them generate the whole system (and so overall
independence fails). \fin

In this section we will employ the general machinery of non-ergodic
isometric extensions and the non-ergodic Furstenberg-Zimmer and
Mackey theories to describe this joint distribution in the case
$d=3$.  It will turn out that these factors are always relatively
independent outside certain special `obstruction' factors, which are
in turn only a little more general than the above example.

\begin{thm}\label{thm:triple-isotropy-1}
Suppose that $T_i:\bbZ\curvearrowright (X,\mu)$, $i=1,2,3$, are
three commuting actions. Then
\begin{enumerate}
\item[(1)] The triple of factors $\zeta_0^{T_1,T_2}$, $\zeta_0^{T_1,T_3}$, $\zeta_0^{T_2,T_3}$ is relatively independent
over $\zeta_0^T$;
\item[(2)] The triple of factors $\zeta_0^{T_1}$, $\zeta_0^{T_2}$, $\zeta_0^{T_3}$ is
relatively independent the further triple of factors
\[\zeta_0^{T_1}\wedge (\zeta_0^{T_2}\vee \zeta_0^{T_3}),\quad \zeta_0^{T_2}\wedge (\zeta_0^{T_3}\vee \zeta_0^{T_1}),\quad \zeta_0^{T_3}\wedge (\zeta_0^{T_1}\vee \zeta_0^{T_2}).\]
\end{enumerate}
\end{thm}

\begin{thm}\label{thm:triple-isotropy-2}
We have
\[\zeta_0^{T_3}\wedge (\zeta_0^{T_1}\vee \zeta_0^{T_2}) \succsim \zeta_0^{T_1,T_3}\vee\zeta_0^{T_2,T_3},\]
and the extension of systems
\[(\zeta_0^{T_1,T_3}\vee\zeta_0^{T_2,T_3})\big|_{\zeta_0^{T_3}\wedge (\zeta_0^{T_1}\vee \zeta_0^{T_2})}:\big(\zeta_0^{T_3}\wedge (\zeta_0^{T_1}\vee \zeta_0^{T_2})\big)(\bfX) \to \zeta_0^{T_1,T_3}\vee\zeta_0^{T_2,T_3}(\bfX)\]
can be coordinatized as the group extension
\begin{center}
$\phantom{i}$\xymatrix{
(\zeta_0^{T_1,T_3}\vee\zeta_0^{T_2,T_3})(\bfX)\ltimes(G_{3,\bullet},m_{G_{3,\bullet}},(\tau_{3,1}\circ\zeta_0^{T_2,T_3}),(\tau_{3,2}\circ\zeta_0^{T_1,T_3})^{\rm{op}},1)\ar[d]^{\rm{canonical}}\\
\zeta_0^{T_1,T_3}\vee\zeta_0^{T_2,T_3}(\bfX) }
\end{center}
for some $T|_{\zeta_0^{T_1,T_3}\vee\zeta_0^{T_2,T_3}}$-invariant
compact group data $G_{3,\bullet}$ and cocycle-sections
$\tau_{3,1}:Z_0^{T_2,T_3} \to G_{3,\bullet}$ and
$\tau_{3,2}:Z_0^{T_1,T_3}\to G_{3,\bullet}$ and similarly for the
extension
$(\zeta_0^{T_i,T_j}\vee\zeta_0^{T_i,T_k})\big|_{\zeta_0^{T_i}\wedge
(\zeta_0^{T_j}\vee \zeta_0^{T_k})}$ for any other permutation
$i,j,k$ of the indices $1,2,3$ (in general with different group data
$G_{i,\bullet}$ and cocycle-sections $\tau_{i,j}$).
\end{thm}

These two results together amount to
Theorem~\ref{thm:threefoldisotropy} of the Introduction.
Heuristically, they assert that the joint distribution of the
factor-triple $(\zeta_0^{T_i})_{i=1}^3$ fails to be relatively
independent over the natural candidate factor-triple
$(\zeta_0^{T_i,T_j}\vee \zeta_0^{T_i,T_k})_{i=1}^3$ only up to
single isometric extensions, and give fairly explicit
coordinatizations of those extensions.

It seems likely that these results could be extended with only
routine modifications to treat a triple of commuting actions
$T_i:\G_i\curvearrowright (X,\mu)$ of other locally compact second
countable groups $\G_i$ (the key feature being that the different
actions commute). We have restricted to a triple of $\bbZ$-actions
for notational simplicity.

On the other hand, although we naturally expect
Theorem~\ref{thm:threefoldisotropy} to be a special case of a result
for larger numbers of commuting transformations (or actions), the
analysis of the corresponding isotropy factors based on the
Furstenberg-Zimmer and Mackey theories becomes quickly much more
complicated, and we shall not pursue this generalization any further
at present.

Observe also that while Theorem~\ref{thm:triple-isotropy-2}
describes the structure of each system $\big(\zeta_0^{T_i}\wedge
(\zeta_0^{T_j}\vee \zeta_0^{T_k})\big)(\bfX)$ as an extension of
$\zeta_0^{T_i,T_j}\vee\zeta_0^{T_i,T_k}(\bfX)$ (which, by
Theorem~\ref{thm:triple-isotropy-1}, is itself just a relatively
independent joining of $\bfZ_0^{T_i,T_j}$ and $\bfZ_0^{T_i,T_k}$
over $\bfZ_0^T$), it does not describe the joint distribution of the
factor maps $\zeta_0^{T_i}\wedge (\zeta_0^{T_j}\vee \zeta_0^{T_k})$.
This would require a further analysis, using the relative
independence of the isotropy factors $\bfZ_0^{T_i,T_j}$ over
$\bfZ_0^T$ to understand first the joint distribution of the
$\zeta_0^{T_i,T_j}\vee\zeta_0^{T_i,T_k}(\bfX)$ and then working
upwards, and would proceed using very similar ideas to those below
but with relatively smaller returns; we omit the details.

Our basic approach rests on an appeal to the Furstenberg-Zimmer
inverse theory to reduce the problem to the study of certain
isometric extensions, followed by a detailed analysis of the
possible structure of an associated Mackey group to obtain finer
information about these extensions.

This strategy is already well-established in the literature from
studies of other questions working under more restrictive ergodicity
assumptions. Indeed, Furstenberg's original paper~\cite{Fur77}
developing an ergodic-theoretic approach to Szemer\'edi's Theorem,
for which much of the abovementioned machinery was originally
developed, uses similar ideas to analyze the structure of a certain
self-joining of a given ergodic $\bbZ$-system on route to the proof
of multiple recurrence.  That paper has since lead to a considerably
more detailed study of the `nonconventional ergodic averages' that
appear in this connexion, which we will revisit in the next
section~\cite{ConLes84,ConLes88.1,ConLes88.2,Zha96,FurWei96,HosKra05,Zie07}.

In addition, Rudolph has given in~\cite{Rud93} an analysis of a
different question rather more closely related to the study of
isotropy factors: he obtains a description of the possible
eigenfunctions of the product system $S\times T\curvearrowright
(Y\times X,\nu\otimes \mu)$ built from ergodic transformations
$S\curvearrowright (Y,\nu)$ and $T\curvearrowright (X,\mu)$,
effectively using for the earlier stages of his work a special case
of the analysis to be given below applied to this latter product
system.

We depart from these previous works in our use of the non-ergodic
versions of the basic machinery.  In this more general setting we
will find that the resulting structures are considerably more
complex, even though the description of the extensions ultimately
obtained in Theorem~\ref{thm:triple-isotropy-2} involves only data
$G_{i,\bullet}$ that are invariant for the whole action $T$. In
particular, many of these older works have ultimately reduced their
subjects to the study of factors that lie in a very special class of
systems, the `pronilsystems' (see, in particular, Rudolph's
work~\cite{Rud93} and the papers of Host and Kra~\cite{HosKra05} and
Ziegler~\cite{Zie07} on nonconventional averages).  Already in the
cases considered in this section we find that we must look beyond
that class.

\subsection{Reduction to compositions of isotropy factors}

We first prove the (rather simpler)
Theorem~\ref{thm:triple-isotropy-1}, by effecting a quite general
reduction of the problem to the study of certain composed isotropy
factors.

\begin{lem}\label{lem:joint-dist-isotropy-and-something}
Suppose that $T:\G\curvearrowright (X,\mu)$ is a
probability-preserving action of a locally compact secound countable
amenable group $\G$ and that $\pi:(X,\mu,T)\to (Y,\nu,S)$ is a
factor. Then $\zeta_0^T$ and $\pi$ are relatively independent over
$\zeta_0^S\circ \pi = \pi|_{\zeta_0^T}\circ \zeta_0^T$.
\end{lem}

\textbf{Proof}\quad Let $(I_N)_{N\geq 1}$ be a left-F\o lner
sequence in $\G$. If $A \in \pi^{-1}(\S_Y)$ and $B \in \S_X$ is
$T$-invariant, then
\begin{multline*}\mu(A\cap B) = \lim_{N\to
\infty}\int_X\Big(\frac{1}{m_\G(I_N)}\int_{I_N}1_{T^{\g}(A)}\,m_\G(\d
\g)\Big)\cdot 1_B\,\d\mu\\ =\int_X
\sfE_\mu(1_A\,|\,\zeta_0^S\circ\pi)\cdot 1_B\,\d\mu =\int_X
\sfE_\mu(1_A\,|\,\zeta_0^S\circ\pi)\cdot
\sfE_\mu(1_B\,|\,\zeta_0^S\circ\pi)\,\d\mu,
\end{multline*}where the middle equality follows from the mean ergodic theorem. \qed

\textbf{Proof of
Theorem~\ref{thm:triple-isotropy-1}}\quad\textbf{(1)}\quad If $A \in
\S_X$ is $T_1$- and $T_2$-invariant, $B \in \S_X$ is $T_1$- and
$T_3$-invariant and $C \in \S_X$ is $T_2$- and $T_3$-invariant then
averaging first under $T_3$ gives
\[\int_X 1_A\cdot 1_B\cdot 1_C\,\d\mu = \int_X\sfE_\mu(1_A\,|\,\zeta_0^{T_3})\cdot 1_B\cdot 1_C\,\d\mu = \int_X\sfE_\mu(1_A\,|\,\zeta_0^{T_1,T_2,T_3})\cdot 1_B\cdot 1_C\,\d\mu,\]
and now averaging under $T_2$ gives
\begin{eqnarray*}&&\int_X\sfE_\mu(1_A\,|\,\zeta_0^{T_1,T_2,T_3})\cdot 1_B\cdot
1_C\,\d\mu\\
&&= \int_X\sfE_\mu(1_A\,|\,\zeta^{T_1,T_2,T_3})\cdot
\sfE_\mu(1_B\,|\,\zeta_0^{T_2})\cdot 1_C\,\d\mu\\ &&=
\int_X\sfE_\mu(1_A\,|\,\zeta_0^{T_1,T_2,T_3})\cdot
\sfE_\mu(1_B\,|\,\zeta_0^{T_1,T_2,T_3})\cdot
1_C\,\d\mu\\
&&= \int_X\sfE_\mu(1_A\,|\,\zeta_0^{T_1,T_2,T_3})\cdot
\sfE_\mu(1_B\,|\,\zeta_0^{T_1,T_2,T_3})\cdot
\sfE_\mu(1_C\,|\,\zeta_0^{T_1,T_2,T_3})\,\d\mu;\end{eqnarray*}
concatenating these equalities gives the result.

\textbf{(2)}\quad This follows similarly.  For this proof let
$\psi_i := \zeta_0^{T_i}\wedge (\zeta_0^{T_j}\vee \zeta_0^{T_k})$.
If $A\in\S_X$ is $T_1$-invariant, $B\in\S_X$ is $T_2$-invariant and
$C\in\S_X$ is $T_3$-invariant then
Lemma~\ref{lem:joint-dist-isotropy-and-something} applied to the
action $T_i$ and the factors $\zeta_0^{T_i}$ and $\pi:=
\zeta_0^{T_j}\vee \zeta_0^{T_k}$ gives that these are relatively
independent over $\psi_i$, and hence that
\begin{multline*}
\int_X 1_A\cdot 1_B\cdot 1_C\,\d\mu = \int_X
\sfE_\mu(1_A\,|\,\psi_1)\cdot 1_B\cdot 1_C\,\d\mu = \int_X
\sfE_\mu(1_A\,|\,\psi_1)\cdot \sfE_\mu(1_B\,|\,\psi_2)\cdot 1_C\,\d\mu\\
= \int_X \sfE_\mu(1_A\,|\,\psi_1)\cdot \sfE_\mu(1_B\,|\,\psi_3)\cdot
\sfE_\mu(1_C\,|\,\psi_3)\,\d\mu,
\end{multline*}
as required. \qed

\subsection{Some isometric extensions and their associated Mackey data}

To prove Theorem~\ref{thm:triple-isotropy-2} (and so complete the
proof of Theorem~\ref{thm:threefoldisotropy}) we need to understand
the structure of the composite factors $\zeta_0^{T_i}\wedge
(\zeta_0^{T_j}\vee \zeta_0^{T_k})$ as extensions of
$\zeta_0^{T_i,T_j}\vee \zeta_0^{T_i,T_k}$; most of our work will go
into this. As in the statement of the theorem we will treat the case
$(i,j,k) = (3,1,2)$, the others being analogous.  We will first
obtain some isometricity for the extensions
\[(\zeta_0^{T_1,T_3}\vee\zeta_0^{T_2,T_3})\big|_{\zeta_0^{T_3}\wedge (\zeta_0^{T_1}\vee \zeta_0^{T_2})}:\big(\zeta_0^{T_3}\wedge (\zeta_0^{T_1}\vee \zeta_0^{T_2})\big)(\bfX) \to \zeta_0^{T_1,T_3}\vee\zeta_0^{T_2,T_3}(\bfX),\]
and will then see a gradual extraction of finer and finer properties
of these isometric extensions from an analysis of the associated
Mackey data, with an occasional recoordinatization of the extensions
where necessary.

The various isotropy factors stand related as in the following
commutative diagram (where some of the obvious maps have not been
named):
\begin{center}
$\phantom{.}$ \xymatrix{ & \bfX \ar[ddl]_{\zeta_0^{T_1}\vee
\zeta_0^{T_2}}\ar[ddd]\ar[ddr]^{\zeta_0^{T_3}}\ar@/^3pc/[dddrrr]^{\phantom{xx}\zeta_0^{T_1,T_3}\vee \zeta_0^{T_2,T_3}} \\
\\
(\zeta_0^{T_1}\vee
\zeta_0^{T_2})(\bfX)\ar[dr] & & \bfZ_0^{T_3}\ar[dl] \\
& \big(\zeta_0^{T_3}\wedge (\zeta_0^{T_1}\vee
\zeta_0^{T_2})\big)(\bfX)\ar[rrr] & & &
\zeta_0^{T_1,T_3}\vee\zeta_0^{T_2,T_3}(\bfX) }
\end{center}

It will prove helpful to introduce some more notation.  For $i =
1,2$ let $\a_i := \zeta_0^{T_i,T_3}\vee\zeta_0^{T_1,T_2}$ and
$\bfW_i = (W_i,(\a_i)_\#\mu,T|_{\a_i})$ be its target system.

\begin{lem}\label{lem:composite-isotropy-in-isometrics}
We have
\[\zeta_0^{T_3}\wedge (\zeta_0^{T_1}\vee
\zeta_0^{T_2}) \precsim
(\zeta_{1/\a_1}^{T_3}\wedge\zeta_0^{T_1})\vee
(\zeta_{1/\a_2}^{T_3}\wedge\zeta_0^{T_2}).\]
\end{lem}

\textbf{Proof}\quad By
Lemma~\ref{lem:joint-dist-isotropy-and-something} the system
$(\zeta_0^{T_1}\vee\zeta_0^{T_2})(\bfX)$ is a relatively independent
joining of $\bfZ_0^{T_1}$ and $\bfZ_0^{T_2}$ over their further
factors $\zeta_0^{T_1,T_2}|_{\zeta_0^{T_1}}$ and
$\zeta_0^{T_1,T_2}|_{\zeta_0^{T_2}}$.  On the other hand the
$T_3|_{\zeta_0^{T_1}\vee \zeta_0^{T_2}}$-invariant functions on
$(\zeta_0^{T_1}\vee\zeta_0^{T_2})(\bfX)$ are all virtually
measurable with respect to the maximal subextension of
\[\zeta_0^{T_1,T_2}|_{\zeta_0^{T_1}\vee
\zeta_0^{T_2}}:(\zeta_0^{T_1}\vee\zeta_0^{T_2})(\bfX) \to
\bfZ_0^{T_1,T_2}\] that is isometric for the restricted action of
$T_3$, and so Theorem~\ref{thm:rel-ind-joinings} implies that this
in turn is contained in
$(\zeta_{1/\a_1}^{T_3}\wedge\zeta_0^{T_1})\vee
(\zeta_{1/\a_2}^{T_3}\wedge\zeta_0^{T_2})$, as required. \qed

At this point we will introduce some new notation for the basic
systems and factor maps under study.  In addition to lightening the
presentation, this will make our main technical results
simultaneously relevant to this and the next section and so minimize
the duplication of effort.

We have defined $\a_i := \zeta_0^{T_i,T_3}\vee\zeta_0^{T_1,T_2}$
with target $\bfW_i$ above. We define also $\zeta_i :=
\zeta_{1/\a_i}^{T_3}\wedge\zeta_0^{T_i}$ and let $\bfZ_i$ its target
system, and we let $\bfZ$ be the target of $\zeta := \zeta_1\vee
\zeta_2$ (a joining of $\bfZ_1$ and $\bfZ_2$) and $\bfW$ be the
target of $\a := \a_1\vee \a_2$ (a joining of $\bfW_1$ and
$\bfW_2$). As usual the choice of these target systems is arbitrary
up to isomorphism, but in this case it is natural (and notationally
convenient) to pick $\bfW$ to be
\[(W_1\times W_2,(\a_1\vee\a_2)_\#\mu,T|_{\a_1}\times T_{\a_2}),\]
since we will often want to discuss separately the two coordinates
of a point $(w_1,w_2) \in W$. These factors are now arranged as
shown:
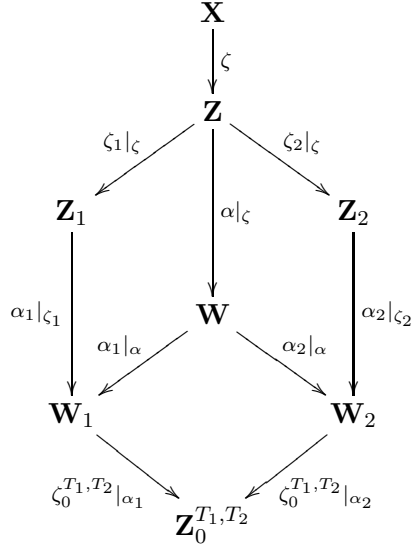
\begin{center}
$\phantom{i}$\xymatrix{ & \bfX\ar[d]^{\zeta}\\
& \bfZ\ar[dl]_{\zeta_1|_{\zeta}}\ar[dd]^{\a|_{\zeta}}\ar[dr]^{\zeta_2|_{\zeta}}\\
\bfZ_1\ar[dd]_{\a_1|_{\zeta_1}} & & \bfZ_2\ar[dd]^{\a_2|_{\zeta_2}}\\
& \bfW\ar[dl]_{\a_1|_{\a}}\ar[dr]^{\a_2|_{\a}}\\
\bfW_1\ar[dr]_{\zeta_0^{T_1,T_2}|_{\a_1}} & & \bfW_2\ar[dl]^{\zeta_0^{T_1,T_2}|_{\a_2}}\\
& \bfZ_0^{T_1,T_2}}
\end{center}
where $\bfW$ and $\bfZ$ are actually generated by all of their
exhibited factors, and the factors on left- and right-hand sides are
relatively independent over their factor maps to $\bfZ_0^{T_1,T_2}$.

In this picture the transformations $T_i = T^{\bf{e}_i}$ have the
following properties:
\begin{itemize}
\item $T_i$ restricts to the identity on $\bfZ_i$ and the factors
beneath it, while acting relatively ergodically on the extension
$\zeta_0^{T_1,T_2}|_{\zeta_{3-i}}:\bfZ_{3-i}\to \bfZ_0^{T_1,T_2}$,
for $i=1,2$;
\item the extensions $\zeta_0^{T_1,T_2}|_{\a_i}:\bfW_i\to \bfZ_0^{T_1,T_2}$ are relatively
invariant for the restriction of $T_3$, and the extensions
$\a_i|_{\zeta_i}:\bfZ_i \to \bfW_i$ are relatively ergodic and
isometric for the restriction of $T_3$.
\end{itemize}
Our goal is to identify the $T_3$-invariant factor of $\bfZ$ (which
we know is also the overall $T_3$-invariant factor by the above
lemma).

In these terms we can now state our main technical result.

\begin{prop}\label{prop:describe-invt-factor}
In the situation described above, there are intermediate factors
\[\bfZ_i \stackrel{\xi_i|_{\zeta_i}}{\longrightarrow} \bfY_i
\stackrel{\a_i|_{\xi_i}}{\longrightarrow} \bfW_i\] factorizing
$\a_i|_{\zeta_i}$ such that there are $T$-invariant compact group
data $G_\bullet$ and cocycle-sections
\begin{itemize}
\item[] $\s:Z_0^{T_1,T_2}\to G_\bullet$ that is
$T_3|_{\zeta_0^{T_1,T_2}}$-relatively ergodic,
\item[] $\tau_1:Z_0^{T_2,T_3}\to G_\bullet$ that is
$T_1|_{\zeta_0^{T_2,T_3}}$-relatively ergodic and
\item[] $\tau_2:Z_0^{T_2,T_3}\to G_\bullet$ that is
$T_2|_{\zeta_0^{T_1,T_3}}$-relatively ergodic
\end{itemize}
so that we can coordinatize
\begin{center}
$\phantom{i}$\xymatrix{
\bfY_1\ar[dr]_{\a_1|_{\xi_1}}\ar@{<->}[rr]^-\cong &
& \bfW_1\ltimes (G_\bullet,m_{G_\bullet},1,(\tau_2\circ\zeta_0^{T_1,T_3}|_{\a_1})^{\rm{op}},\s\circ\zeta_0^{T_1,T_2}|_{\a_1})\ar[dl]^{\rm{canonical}}\\
& \bfW_1, }
\end{center}
and
\begin{center}
$\phantom{i}$\xymatrix{
\bfY_2\ar[dr]_{\a_2|_{\xi_2}}\ar@{<->}[rr]^-\cong &
& \bfW_2\ltimes (G_\bullet,m_{G_\bullet},(\tau_1\circ\zeta_0^{T_2,T_3}|_{\a_2})^{\rm{op}},1,\s\circ\zeta_0^{T_1,T_2}|_{\a_2})\ar[dl]^{\rm{canonical}}\\
& \bfW_2, }
\end{center}
and such that the $T_3$-invariant factor $\zeta_0^{T_3}\wedge
(\zeta_0^{T_1}\vee \zeta_0^{T_2})$ is contained in $\xi_1\vee
\xi_2$.
\end{prop}

\textbf{Proof of Theorem~\ref{thm:triple-isotropy-2} from
Proposition~\ref{prop:describe-invt-factor}}\quad This now follows
simply by unpacking the new notation.  Let $\xi := \xi_1\vee \xi_2$
with target $\bfY$ (a joining of $\bfY_1$ and $\bfY_2$). We know
that $\xi_1$ and $\xi_2$ (as factors of $\zeta_0^{T_1}$ and
$\zeta_0^{T_2}$) are relatively independent over their further
factors $\zeta_0^{T_1,T_2}|_{\xi_i}$, $i=1,2$, and hence certainly
over the intermediate factors $\a_1$ and $\a_2$, and so the
coordinatizations of the extensions $\a_i|_{\xi_i}:\bfY_i\to\bfW_i$
by group data given by Proposition~\ref{prop:describe-invt-factor}
combine to give a coordinatization of $\bfY$ by the group data
$G_\bullet^2$ and the combined cocycles.  We now observe that the
restriction of $T_3$ to this group data extension is described by
the diagonal cocycle-section $(\s,\s)$ corresponding to the
$T_3|_{\zeta_0^{T_1,T_2}}$-ergodic cocycle-section $\s$, and so we
can simply deduce that the Mackey group data can be taken to be the
diagonal subgroup $M_\bullet \cong \{(g,g):\ g\in G_\bullet\}$, and
now the associated Mackey section is trivial by symmetry. This leads
to the coordinatization of the $T_3$-invariant factor
$\zeta_0^{T_3}\wedge(\zeta_0^{T_1}\vee\zeta_0^{T_2}) =
\zeta_0^{T_3}\wedge\zeta = \zeta_0^{T_3}\wedge\xi$ as given by the
location of $\bfZ_0^{T_3|_\xi}$ in the following commutative diagram
\begin{center}
$\phantom{i}$\xymatrix{\bfY\ar[d]^{\zeta_0^{T_3}|_\xi}\ar@{<->}[rr]^-\cong & & \bfW\ltimes (G_\bullet^2,m_{G_\bullet^2},(\tau_1^{\rm{op}},1),(1,\tau_2^{\rm{op}}),(\s,\s))\ar[d]^{\rm{canonical}}\\
\bfZ_0^{T_3|_{\xi}}\ar[dr]_{\a|_{\zeta_0^{T_3|_\xi}}}\ar@{<->}[rr]^-\cong
&
& (\zeta_0^{T_3}\wedge\a)(\bfX)\ltimes (M_\bullet\backslash G^2_\bullet,m_{M_\bullet\backslash G^2_\bullet},(\tau_1^{\rm{op}},1),(1,\tau_2^{\rm{op}}),1)\ar[dl]^{\rm{canonical}}\\
& (\zeta_0^{T_3}\wedge\a)(\bfX) }
\end{center}
(where we have suppressed the need to lift $\tau_i$ through
$\zeta_0^{T_{3-i},T_3}$). Now simply observing that the quotient
$M_\bullet\backslash G_\bullet^2$ is canonically bijective with
$G_\bullet$ under the map $M_\bullet\cdot (g_1,g_2)\leftrightarrow
g_1^{-1}\cdot g_2$ and applying this bijection fibrewise, the
restricted action of $T_3$ on $\zeta_0^{T_3|_\xi}$ is of course
trivial and the restricted actions of $T_1$ and $T_2$ turn into the
respective left- and right-actions by the cocycles $\tau_1$ and
$\tau_2$ asserted in Theorem~\ref{thm:triple-isotropy-2} (where some
additional subscripts `$_3$' from the statement of that theorem have
also been suppressed). \qed

We will prove Proposition~\ref{prop:describe-invt-factor} in several
steps.  First observe that since the extensions
$\a_i|_{\zeta_i}:\bfZ_i\to \bfW_i$ are isometric and relatively
ergodic for $T_3|_{\zeta_i}$, the non-ergodic Furstenberg-Zimmer
Theory of Section~\ref{sec:FZ} enables us to pick coordinatizations
by core-free homogeneous space data for the
$(\bbZ\bf{e}_3)$-subactions
\begin{center}
$\phantom{i}$\xymatrix{ \bfZ_i^{\
\uhr\bf{e}_3}\ar[dr]_{\a_i|_{\zeta_i}}\ar@{<->}[rr]^-\cong & &
\bfW_i^{\ \uhr\bf{e}_3}\ltimes
(G'_{i,\bullet}/K'_{i,\bullet},m_{G'_{i,\bullet}/K'_{i,\bullet}},\s'_i)\ar[dl]^{\rm{canonical}}\\
& \bfW_i^{\ \uhr\bf{e}_3} }
\end{center}
(recall that $\bfZ_i^{\ \uhr\bf{e}_3}$ denotes the subaction system
given by retaining only the action through $T$ of the
one-dimensional subgroup $\bbZ\bf{e}_3 \leq\bbZ^3$), where we may
also choose the cocycle sections $\s_i'$ to be ergodic.

Since $\zeta_1$ and $\zeta_2$ (like $\zeta_0^{T_1}$ and
$\zeta_0^{T_2}$) are relatively independent over $\zeta_0^{T_1,T_2}$
under $\mu$, we can combine the above two coordinatizations to give
\begin{center}
$\phantom{i}$\xymatrix{\bfZ^{\
\uhr\bf{e}_3}\ar[dr]_{\a|_\zeta}\ar@{<->}[rr]^-\cong & & \bfW^{\
\uhr\bf{e}_3}\ltimes
(\vec{G}_\bullet/\vec{K}_\bullet,m_{\vec{G}_\bullet/\vec{K}_\bullet},\vec{\s})\ar[dl]^{\rm{canonical}}\\
& \bfW^{\ \uhr\bf{e}_3} }
\end{center}
where $\vec{G}_\bullet := G'_{1,\pi_1(\bullet)}\times
G'_{2,\pi_2(\bullet)}$, $\vec{K}_\bullet :=
K'_{1,\pi_1(\bullet)}\times K'_{2,\pi_2(\bullet)}$ and $\vec{\s} :=
(\s'_1\circ\pi_1,\s'_2\circ\pi_2)$, and we here write $\pi_i$ for
the obvious factor map $W \to W_i$.

Of course, we do not know that the restrictions of $T_1$ and $T_2$
to the factors $\a_i|_{\zeta_i}:\bfZ_i\to \bfW_i$ are isometric, and
so we have no similar coordinatization of these transformations
using homogeneous space data and cocycles. We will appeal instead
the the Relative Automorphism Structure Theorem~\ref{thm:RAST} to
describe them in terms of cocycles and fibrewise automorphisms.

First, however, an appeal to Corollary~\ref{cor:homo-nonergMackey}
gives a first step towards the more explicit description of the
$T_3|_\xi$-invariant factor in terms of the above coordinatizations:

\begin{prop}\label{prop:intro-the-Mackey-group} There are
$T_3|_\a$-invariant Mackey group data $M'_\bullet \leq
\vec{G}_\bullet$ on $W$ and a measurable section $\vec{b}:W \to
\vec{G}_\bullet$ such that the factor map
\begin{multline*}
Z \to
Z_0^{T_3|_\a}\ltimes(M'_\bullet\backslash\vec{G}_\bullet/\vec{K}_\bullet):\\
(w,g \vec{K}_w) \mapsto
\big((\zeta_0^{T_1,T_3}\vee\zeta_0^{T_2,T_3})(w), M'_w\cdot
\vec{b}(w)\cdot g\cdot \vec{K}_w\big)
\end{multline*}
coordinatizes the $T_3|_{\xi}$-invariant factor of $\bfZ$. \qed
\end{prop}

The remainder of our work will go into analyzing this Mackey group
data $M'_\bullet$ and section $b'$ to deduce properties of the data
$(G'_{i,\bullet}/K'_{i,\bullet},\s'_i)$ that gave rise to them, and
eventually reduce them to the special form promised by
Proposition~\ref{prop:describe-invt-factor}.

We now prove two technical lemmas that will underly our subsequent
analysis, and which it seems easiest to introduce separately.

\begin{lem}\label{lem:rel-inv-Mackey-group}
Suppose that $(X,\mu,T)$ is a $\bbZ$-system, $G_\bullet$ is
$T$-invariant measurable compact group data on $X$ and $\s:X\to
G_\bullet$ a cocycle-section, and that $(Y,\nu)$ is another standard
Borel probability space.  Suppose further that $\l$ is a $(T\times
\id_Y)$-invariant joining of $\mu$ and $\nu$. If $M_\bullet \leq
G_\bullet$ is the Mackey group data of $G_\bullet$ and $\s$ over
$(X,\mu,T)$ and $\pi:X\times Y\to X$ is the coordinate projection,
then the Mackey group data $N_\bullet$ of $G_{\pi(\bullet)}$ and
$\s\circ\pi$ over $(X\times Y,\l,T\times\id_Y)$ is given by
$M_{\pi(\bullet)}$ (up to a $T$-invariant measurable choice of
conjugates) $\l$-almost surely.
\end{lem}

\textbf{Remark}\quad It is easy to see that $N_\bullet \leq
M_{\pi(\bullet)}$; the point to this proposition is that if we
adjoin to $(X,\mu,T)$ a system on a new space $(Y,\nu)$ for which
the action is trivial, then the Mackey group data does not become
any smaller. \fin

\textbf{Proof}\quad We know from Section~\ref{sec:Mackey} that there
is a section $b:X\times Y\to G_{\pi(\bullet)}$ such that
$b(Tx,y)^{-1}\cdot \s(x)\cdot b(x,y) \in N_{(x,y)}$ for $\l$-almost
every $(x,y)$.

Let $A$ be the $\l$-conegligible subset of $X\times Y$ where this
coboundary condition obtains, and let \[B_0 := \{(x,y) \in A:\
\hbox{some conjugate of $N_{(x,y)}$ is properly contained in
$M_x$}\};\] this is easily seen to be Borel and $\l$-almost
$(T\times \id_Y)$-invariant, and so writing $B:=
\bigcap_{n\in\bbZ}T^n(B_0)$ we see that $\l (B) = \l (B_0)$ and that
$B$ is strictly $(T\times\id_Y)$-invariant. It will suffice to show
that $B$ is $\l$-negligible, so suppose otherwise. Then by
Proposition~\ref{prop:invar-meas-select} there are a non-negligible
$T$-invariant subset $C \in \S_X$ and a $T$-invariant measurable
selector $\eta:C\to Y$ such that $(x,\eta(x)) \in B$ almost surely.
We deduce that $b(Tx,\eta(x))^{-1}\cdot \s(x)\cdot b(x,\eta(x)) \in
N_{(x,\eta(x))}$ for every $x \in C$ with $N_{(x,\eta(x))}$ properly
contained in some (clearly measurably-varying) conjugate of the
Mackey group data $M_x$, contradicting the conjugate-minimality of
this latter that was proved part (4) of
Theorem~\ref{thm:nonergMackey}. \qed

\begin{cor}\label{cor:full-1D-proj}
In the notation set up earlier in this section, we have
\[\{g_1:\ \exists g_2 \in G'_{2,\bullet}\ \rm{s.t.}\ (g_1,g_2)\in M'_\bullet\} = G'_{1,\bullet}\]
almost surely, and similarly for the projection of $M'_\bullet$ onto
$G'_{2,\bullet}$.
\end{cor}

\textbf{Proof}\quad We give the argument for $i=1$. Simply observe
that the extension $\bfW^{\ \uhr\bf{e}_3}\to \bfW_1^{\
\uhr\bf{e}_3}$ is relatively invariant, and so the Mackey group data
for our coordinatization of the extension
$\a|_{\zeta_1\vee\a_2}:(\zeta_1\vee\a_2)(\bfZ)^{\ \uhr\bf{e}_3}\to
\bfW^{\ \uhr\bf{e}_3}$ must simply be lifted from the Mackey group
data for $\a_1|_{\zeta_1}:\bfZ_1^{\ \uhr\bf{e}_3}\to \bfW_1^{\
\uhr\bf{e}_3}$ downstairs. Since the former is clearly equal to the
given one-dimensional projection of $M'_\bullet$, and the cocycle
section $\s_1'$ is assumed to be ergodic, this completes the proof.
\qed

The next properties of $M'_\bullet$ that we deduce require a little
more work.  We begin with a useful group-theoretic lemma.

\begin{lem}[Deconstructing a relation between two group
correspondences]\label{lem:slices-respected} Suppose that $G_1$,
$G_2$ are compact groups and that $M_1,M_2\leq G_1\times G_2$ are
two subgroups that both have full one-dimensional projections, and
let their one-dimensional slices be
\[L_{1,1} := \{g \in G_1:\ (g,1_{G_2})\in M_1\},\quad\quad L_{1,2}:= \{g \in G_2:\ (1_{G_1},g)\in M_1\}\]
and similarly $L_{2,1}$, $L_{2,2}$.  Suppose further that
$\Phi_i:G_i\stackrel{\cong}{\longrightarrow} G_i$ and $h_i,k_i \in
G_i$ for $i=1,2$ satisfy
\[(h_1,h_2)\cdot(\Phi_1\times \Phi_2)(M_1)\cdot(k_1,k_2) = M_2.\]
Then $\Phi_i(L_{1,i}) = L_{2,i}$ for $i=1,2$.
\end{lem}

\textbf{Proof}\quad Suppose first that $(g,1_{G_2}) \in L_{1,1}$.
Then the given equation tells us that
\[(h_1\cdot \Phi_1(g)\cdot k_1,h_2\cdot k_2) = (m_1,m_2)\]
for some $m_1,m_2 \in M_2$, and in this case we have that $m_2 =
h_2\cdot k_2$ does not depend on $g$.  Since the above must
certainly hold if $g = 1_{G_1}$, applying it also for any other $g$
and differencing gives
\[(h_1\cdot \Phi_1(g)\cdot k_1)\cdot(h_1\cdot \Phi_1(1_{G_1})\cdot k_1)^{-1} = h_1\cdot\Phi_1(g)\cdot h_1^{-1} \in L_{2,1},\]
so $\Phi_1(L_{1,1}) \subseteq h_1^{-1}\cdot L_{2,1}\cdot h_1$.  An
exactly symmetric argument gives the reverse inclusion, so in fact
$\Phi_1(L_{1,1})$ is a conjugate of $L_{2,1}$.  However, since $M_1$
and $M_2$ have full one-dimensional projections, by
Lemma~\ref{lem:groupcorrespondence} it follows that in fact
$\Phi_1(L_{1,1}) = L_{2,1}$, as required. The case of the other
coordinate is similar. \qed

\begin{lem}\label{lem:1D-slices-restricted-depdce}
If $H_{i,\bullet} \leq G'_{i,\pi_i(\bullet)}$ are the
one-dimensional slices of $M'_\bullet$, then
\begin{enumerate}
\item[(1)] $H_{i,(w_1,w_2)}$ $\a_\#\mu$-almost surely depends only on
$w_i$, so after modifying on a negligible set we may write it as
$H_{i,w_i}$;
\item[(2)] under the above coordinatizations, for $i=1,2$ the map
\begin{multline*}
W_i\ltimes (G'_{i,\bullet}/K'_{i,\bullet})\to W_i\ltimes
(G'_{i,\bullet}/(H_{i,\bullet}K'_{i,\bullet})):\\
(w_i,gK'_{i,w_i})\mapsto (w_i,gH_{i,w_i}K'_{i,w_i})
\end{multline*}
defines a factor for the whole $\bbZ^3$-action $T$ (that is, it is
respected by $T_1$ and $T_2$ as well as $T_3$).
\end{enumerate}
\end{lem}

\textbf{Proof}\quad By symmetry it suffices to treat the case $i=1$
for the first conclusion and $i=2$ for the second (it will turn out
that these come together). First deduce from
Corollary~\ref{cor:full-1D-proj} and
Lemma~\ref{lem:groupcorrespondence} that in fact
$H_{1,(w_1,w_2)}\unlhd G'_{1,w_1}$ for almost every $(w_1,w_2)$.

We will use the presence of the additional transformations of the
factors $\zeta_1$ and $\zeta_2$ given by $T_1$. Of course, $T_1$
just restricts to the identity transformation on $\zeta_1$. On the
other hand, since $T_1|_{\zeta_2}$ commutes with the transformation
$T_3|_{\zeta_2}$ which is relatively ergodic for the extension
$\bfZ_2 \to \bfW_2$, the Relative Automorphism Structure
Theorem~\ref{thm:RAST} allows us to express
\[T_1|_{\zeta_2} \cong T_1|_{\a_2}\ltimes
(L_{\rho'(\bullet)}\circ\Phi'_\bullet)|^{K'_{2,\bullet}}_{K'_{2,T_1|_{\a_2}(\bullet)}}\]
for some $\rho':W_2\to G'_{2,\bullet}$ and $T_3|_{\a_2}$-invariant
section $\Phi'_\bullet:W_2\to
\rm{Isom}(G'_{2,\bullet},G'_{2,T_1|_{\a_2}(\bullet)})$ satisfying
$\Phi'_\bullet(K'_{2,\bullet})= K'_{2,T_1|_{\a_2}(\bullet)}$ almost
surely.

Now observe from Corollary~\ref{cor:homogextiso} that using the
above expression and its partner for $T_2|_{\zeta_1}$ we may extend
all three transformations $T_j|_{\zeta_i}$, $j=1,2,3$, to the
covering group extension $\t{\a_i}:\t{\bfZ_i}\to
\bfZ_i\stackrel{\a_i|_{\zeta_i}}{\longrightarrow} \bfW_i$ arising
from our core-free homogeneous-space-data coordinatization of
$T_3|_{\zeta_i}$, and that these extensions retain commutativity and
all the relative invariance, ergodicity and isometricity properties
listed above. Form the relatively independent joining
\[\t{\bfZ} =
\t{\bfZ_1}\otimes_{\{\zeta_0^{T_1,T_2}|_{\a_1}\circ\t{\a_1} =
\zeta_0^{T_1,T_2}|_{\a_2}\circ\t{\a_2}\}}\t{\bfZ_2}\] with the
coordinate projection factors back onto $\t{\bfZ_1}$ and
$\t{\bfZ_2}$; with the resulting factor map onto $\bfW$ it now
defines a covering group extension of $\a|_{\zeta}:\bfZ\to \bfW$
whose Mackey data are still $M'_\bullet$ and $b'$ (by our initial
construction of these).  Moreover, these new factors $\t{\bfZ}_i$
and $\t{\bfZ}$ are located in a commutative diagram with the factors
$\bfW_i$ and $\bfW$ just as we saw previously for $\bfZ_i$ and
$\bfZ$ (except now not all as factors of the original overall system
$\bfX$, but of some extended overall system).

It follows that for the purpose of proving this proposition, we may
work with these covering group extensions throughout without
disrupting the final conclusions; or, equivalently, that it suffices
to treat the case in which the core-free kernels $K'_\bullet$ are
trivial. Let us therefore make this assumption for the rest of this
proof so as to lighten notation.

Given this assumption, consider the condition that $T_1|_\zeta$
respect $\zeta_0^{T_3}|_\zeta$ in terms of the above expression for
$T_1|_{\zeta_2}$ and the Mackey data.  First, since $M'_\bullet$ has
full one-dimensional projections we may take the Mackey section
$\vec{b}$ of Proposition~\ref{prop:intro-the-Mackey-group} to be of
the form $\vec{b}(w) = (1_{G_{1,w_1}},b'(w))$. Now the above
condition requires, in particular, that
$\zeta_0^{T_3}|_\zeta(T_1|_\zeta(z))$ almost surely depend only on
$\zeta_0^{T_3}|_\zeta(z)$ for $z\in Z$; and on the other hand, in
terms of the above Mackey description, writing points of $Z$ as
$(w,g_1,g_2)$ we know that $\zeta_0^{T_3}|_\zeta(w,g_1,g_2) =
\zeta_0^{T_3}|_\zeta(w,g_1',g_2')$ if and only if
\begin{multline*}
M'_w\cdot (1,b'(w))\cdot (g_1,g_2) = M'_w\cdot (1,b'(w))\cdot
(g_1',g_2')\\
\Leftrightarrow\quad\quad (g_1,g_2) \in (1,b'(w)^{-1})\cdot
M'_w\cdot (1,b'(w))\cdot (g_1',g_2').
\end{multline*}
Therefore the above relation between $T_1|_\zeta$ and
$\zeta_0^{T_3}|_\zeta$ simply asserts that for $\a_\#\mu$-almost
every $(w_1,w_2) \in W$, for Haar-almost every $(g_1',g_2') \in
G'_{1,w_1}\times G'_{2,w_2}$ there is some $(g_1'',g_2'') \in
G'_{1,w_1}\times G'_{2,T_1|_{\a_2}(w_2)}$ such that
\begin{multline*}
(\id_{G'_{1,w_1}}\times
(L_{\rho'(w_2)}\circ\Phi'_{w_2}))\big((1,b'(w_1,w_2)^{-1})\cdot
M'_{(w_1,w_2)}\cdot (1,b'(w_1,w_2))\cdot (g_1',g_2')\big)\\
= (1,b'(w_1,T_1|_{\a_2}(w_2))^{-1})\cdot
M'_{(w_1,T_1|_{\a_2}(w_2))}\cdot (1,b'(w_1,T_1|_{\a_2}(w_2)))\cdot
(g_1'',g_2''),
\end{multline*}
or, re-arranging, that
\begin{eqnarray*}
&&\big(1_{G'_{1,w_1}},b'(w_1,T_1|_{\a_2}(w_2))\rho'(w_2)\Phi'_{w_2}(b'(w_1,w_2)^{-1})\big)\\
&&\quad\quad\quad\quad\cdot
(\id_{G'_{1,w_1}}\times \Phi'_{w_2})(M'_{(w_1,w_2)})\\
&&\quad\quad\quad\quad\quad\quad\quad\quad\cdot \big(g_1'(g_1'')^{-1},\Phi'_{w_2}(b'(w_1,w_2)g_2')(b'(w_1,T_1|_{\a_2}(w_2))g_2'')^{-1}\big)\\
&&\quad\quad\quad\quad\quad\quad\quad\quad\quad\quad\quad\quad\quad\quad\quad\quad\quad\quad\quad\quad\quad\quad\quad\quad=
M'_{(w_1,T_1|_{\a_2}(w_2))}.
\end{eqnarray*}

The two desired conclusions now follow from applying
Lemma~\ref{lem:slices-respected} to this equation for the two
coordinate projections onto $G_{1,w_1}'$ and $G'_{2,w_2}$. Under the
first coordinate projection we obtain
\[H_{1,(w_1,w_2)} = \id_{G'_{1,w_1}}(H_{1,(w_1,w_2)}) = H_{1,(w_1,T_1|_{\a_2}(w_2))},\]
so $H_{1,(w_1,w_2)}$ is a $T_1|_\a$-invariant subgroup of
$G'_{1,w_1}$, and so recalling that $T_1|_\a$ is relatively ergodic
on the extension $\a_1|_\a:W \to W_1$ we deduce that
$H_{1,(w_1,w_2)}$ is virtually a function of $w_1$ alone, as
required for conclusion (1).

For the second coordinate projection we obtain
\[\Phi'_{w_2}(H_{2,(w_1,w_2)}) =
H_{2,(w_1,T_1|_{\a_2}(w_2))}.\] In view of the conclusion (1)
obtained above we can simplify this to
\[\Phi'_{w_2}(H_{2,w_2}) =
H_{2,T_1|_{\a_2}(w_2)},\] and now this is precisely the condition
given by the Relative Automorphism Structure Theorem~\ref{thm:RAST}
for $T_1|_{\zeta_2}$ to respect the given map as a factor map. Since
it is clear that the given map defines a factor for the restrictions
of $T_2$ (since this acts trivially on the whole of $\bfZ_2$) and
$T_3$ (since this acts on this extension by a $G_{2,\bullet}$-valued
cocycle-section, and so our factor map is simply the fibrewise
quotient by the subgroup $H_{i,\bullet}$, recalling our assumption
that $K'_{i,\bullet}$ is almost surely trivial), this completes the
proof. \qed

In view of the above result, we are now able to define our desired
intermediate factors $\bfZ_i
\stackrel{\xi_i|_{\zeta_i}}{\longrightarrow} \bfY_i
\stackrel{\a_i|_{\xi_i}}{\longrightarrow} \bfW_i$ by the commutative
diagrams
\begin{center}
$\phantom{i}$\xymatrix{
\bfZ_1\ar[d]_{\xi_1|_{\zeta_1}}\ar@{<->}[rr]^-\cong & & \bfW_1\ltimes (G'_{1,\bullet}/K'_{1,\bullet},m_{G'_{1,\bullet}/K'_{1,\bullet}},1_{G'_{1,\bullet}},(L_{\rho'_1(\bullet)}\circ\Phi'_{1,\bullet}),\s'_1)\ar[d]^{\rm{canonical}}\\
\bfY_1\ar[dr]_{\a_1|_{\xi_1}}\ar@{<->}[rr]^-\cong & & \bfW_1\ltimes (G_{1,\bullet}/K_{1,\bullet},m_{G_{1,\bullet}/K_{1,\bullet}},1_{G_{1,\bullet}},(L_{\rho_1(\bullet)}\circ\Phi_{1,\bullet}),\s_1)\ar[dl]^{\rm{canonical}}\\
& \bfW_1, }
\end{center}
and similarly for $\bfY_2$, where $G_{i,\bullet} :=
G'_{i,\bullet}/H_{i,\bullet}$, $K_{i,\bullet}:=
(H_{i,\bullet}K'_{i,\bullet})/H_{i,\bullet}$ and $\rho_i$,
$\Phi_{i,\bullet}$ and $\s_i$ are the appropriate quotients or
restrictions of $\rho'_i$, $\Phi'_{i,\bullet}$ and $\s_i'$: part (1)
above gives that $H_{i,\bullet}$ is correctly defined as a function
on $W_i$, and part (2) gives that the above diagram defines a factor
map for our whole $\bbZ^3$-action.

The important feature of these new smaller extensions $\bfY_i \to
\bfW_i$ is that the Mackey group data $M_\bullet$ of their joining
under $\bfX$ takes a particularly simple form: having quotiented out
the one-dimensional slices $H_{i,\bullet}$, it is almost surely the
graph of a continuous isomorphism. Indeed, $M_\bullet$ is clearly
obtained from $M'_\bullet$ simply by quotienting out the normal
subgroup data $H_{1,\pi_1(\bullet)}\times H_{2,\pi_2(\bullet)}$, and
from the definition of $H_{i,\bullet}$ it follows that $M_\bullet$
has full one-dimensional projections and trivial one-dimensional
slices almost everywhere, and so defines almost everywhere the
graphs of some measurably-varying $T_3|_{\a}$-invariant isomorphisms
$\Psi_{(w_1,w_2)}:G_{1,w_1}\stackrel{\cong}{\longrightarrow}
G_{2,w_2}$.  Henceforth we will refer to these as the \textbf{Mackey
isomorphisms}.  In addition we set $b := b'\cdot
H_{2,\pi_2(\bullet)}$, so that $(1_{G_{1,\pi_1(\bullet)}},b)$ is a
Mackey section of the extension $\bfY^{\ \uhr\bf{e}_3} \to \bfW^{\
\uhr\bf{e}_3}$ associated to the choice of Mackey group data
$M_\bullet$.

On the other hand, from the description given in
Proposition~\ref{prop:intro-the-Mackey-group} it follows that the
$T_3$-invariant factor is actually contained in the join of these
smaller isometric extensions $\a_i|_{\xi_i}:\bfY_i \to \bfW_i$, and
so it will suffice to study these new factors. The remaining steps
of this subsection will give a recoordinatization of these new
factors into the form required by
Proposition~\ref{prop:describe-invt-factor}.

\begin{cor} We have $\zeta_0^T|_{\a_1}(w_1)
=\zeta_0^T|_{\a_2}(w_2)$ for $\a_\#\mu$-almost every $(w_1,w_2)$,
and there are compact group data $G_\bullet$ invariant for the whole
action $T$ such that we can recoordinatize the extensions
$\a_i|_{\xi_i}:\bfY_i\to\bfW_i$ so that $G_{1,w_1} =
G_{\zeta_0^T|_{\a_1}(w_1)} = G_{2,w_2}$ for $\a_\#\mu$-almost every
$(w_1,w_2)$.
\end{cor}

\textbf{Proof}\quad The first assertion is clear from the
definitions.

By $T_3|_{\a_i}$-invariance the groups $G_{i,w_i}$ actually depend
only on $z_i:=\zeta_0^{T_3}|_{\a_i}(w_i) \in Z_0^{T_i,T_3}$, and
similarly the isomorphism $\Psi_{(w_1,w_2)}$ depends only on the
image $(z_1,z_2)$ of $(w_1,w_2)$.  In addition, by
Theorem~\ref{thm:triple-isotropy-1} the coordinates $z_1$, $z_2$ of
this image are relatively independent over $\zeta_0^T|_{\a_1}(w_1) =
\zeta_0^T|_{\a_2}(w_w)$ under $\a_\#\mu$.

Now let $P:Z_0^T \stackrel{\rm{p}}{\to} Z_0^{T_2,T_3}$ be a
probability kernel representing the disintegration of
$(\zeta_0^{T_2,T_3})_\#\mu$ over $\zeta_0^T|_{\zeta_0^{T_2,T_3}}$.
For almost every $z_1 \in Z_0^{T_1,T_3}$ we can choose a measurable
family of isomorphisms
$\Theta_{2,z_2}:G_{2,z_2}\stackrel{\cong}{\longrightarrow}
G_{1,z_1}$ defined for
$P(\zeta_0^T|_{\zeta_0^{T_1,T_3}}(z_1),\,\cdot\,)$-almost every $z_2
\in Z_0^{T_2,T_3}$, because the Mackey isomorphisms themselves
witness that these almost surely exist. Making a measurable
selection of such a $z_1$ in each fibre of $Z_0^{T_1,T_3} \to
Z_0^T$, we now take this family as defining a fibrewise isomorphism
recoordinatization of our initial homogeneous-space-data
coordinatization of $\a_2|_{\xi_2}:\bfY_2\to \bfW_2$ obtained above.
This has the effect of adjusting to a coordinatization in which the
covering group of the homogeneous space fibre over $z_2 \in Z_2$ is
$P(\zeta_0^T|_{\zeta_0^{T_1,T_3}}(z_1),\,\cdot\,)$-almost everywhere
equal to $G_{1,z_1}$, and with the kernel of the homogeneous space
fibre given by $\Theta_{2,z_2}(K_{2,z_2})\leq G_{1,z_1}$.

In particular, the covering group data of this new coordinatization
depends only on $\zeta_0^T|_{\zeta_0^{T_1,T_3}}(z_1) =
\zeta_0^T|_{\a_1}(w_1) = \zeta_0^T|_{\a_2}(w_2)$. Exactly similarly
we can now recoordinatize $\a_1|_{\xi_1}:\bfY_1\to\bfW_1$ to have
covering fibre groups also depending only on $\zeta_0^T|_{\a_1}(w_1)
= \zeta_0^T|_{\a_2}(w_2)$.  Since both these recoordinatizations are
by fibrewise isomorphisms that are invariant for the relevant
restrictions of $T_3$, the new coordinatizations of these extensions
that result are still given as cocycle-section extensions for these
restrictions of $T_3$. Finally, in this new coordinatization the
measurable family of Mackey isomorphisms $\Psi_{(w_1,w_2)}$ clearly
shows that after one more fibrewise recoordinatization by a
$T$-invariant isomorphism we are left with the same $T$-invariant
group data $G_\bullet$ everywhere. \qed

Now let us re-apply the Relative Automorphism Structure
Theorem~\ref{thm:RAST} to write
\[T_1|_{\xi_2} = T_1|_{\a_2}\ltimes (L_{\rho_1(\bullet)}\circ\Phi_{1,\bullet})|^{K_{2,\bullet}}_{K_{2,T_1|_{\a_2}(\bullet)}}\]
and
\[T_2|_{\xi_1} = T_2|_{\a_1}\ltimes (L_{\rho_2(\bullet)}\circ\Phi_{2,\bullet})^{H_{1,\bullet}}_{H_{1,T_2|_{\a_1}(\bullet)}}\]
where now $\Phi_{i,w_i}$ is $T_3|_{\a_{3-i}}$-invariant and takes
values in $\Aut\,G_{\zeta_0^T|_{\a_i}(w_i)}$ for $i=1,2$.  In
addition, we recall the notation $\rm{Co}_{\rho(\bullet)}$ for the
fibrewise automorphism of some measurable group data $G_\bullet$
given by fibrewise conjugation by a section $\rho$ of $G_\bullet$.

\begin{prop}\label{prop:final-recoord}
The extensions $\a_i|_{\xi_i}:\bfY_i \to \bfW_i$ can be
recoordinatized by fibrewise affine transformations so that
\begin{enumerate}
\item[(1)] there are cocycles $\tau_i:W_{3-i}\to
G_{\zeta_0^T|_{\a_i}(\bullet)}$ such that $L_{\rho_i(\bullet)}\circ
\Phi_{i,\bullet} = R_{\tau_i(\bullet)}$;
\item[(2)] the Mackey isomorphisms are trivial: $\Psi_\bullet \equiv
\id_{G_{\zeta_0^T|_\a(\bullet)}}$;
\item[(3)] the Mackey section is
trivial: $b \equiv 1_{G_{\zeta_0^T|_\a(\bullet)}}$;
\item[(4)] the cocycle $\tau_i$ is invariant under $T_3|_{\a_i}$ and the cocycle $\s_i$ is invariant under $T_{3-i}|_{\a_i}$.
\end{enumerate}
\end{prop}

\textbf{Proof}\quad This will follow from a careful consideration of
the commutativity conditions relating the expressions for our three
transformations on $\bfY_i$. We make our recoordinatizations in two
steps, the first by fibrewise automorphisms and the second by
fibrewise rotations.  We will construct these so as to guarantee the
asserted properties of the ingredients $\rho$ and $\Phi$, and will
then find that the asserted forms of $\Psi$ and $b$ are an immediate
consequence.

First observe that just as in the proof of
Lemma~\ref{lem:1D-slices-restricted-depdce}, we may lift all of our
commuting transformations $T_j|_{\a_i}$ to the covering group-data
extensions of $\a_i|_{\xi_i}:\bfY_i \to \bfW_i$, and have that
$M_\bullet$ and $b$ will still be Mackey data of their relatively
independent joining over $\bfZ_0^{T_1,T_2}$, and therefore if we
effect our desired fibrewise recoordinatizations on these covering
group-data extensions then simply quotienting will give the desired
recoordinatizations of $\a_i|_{\xi_i}:\bfY_i \to \bfW_i$. As in
Lemma~\ref{lem:1D-slices-restricted-depdce}, this argument reduces
our work to the special case when $K_\bullet \equiv
\{1_{G_\bullet}\}$.

The remainder of our work breaks into five steps.

\quad\textbf{Step 1}\quad We consider the case $i=1$.  First recall
our earlier expression of the fact that $T_1|_{\xi}$ respects
$\zeta_0^{T_3}|_{\xi}$: for $\a_\#\mu$-almost every $(w_1,w_2) \in
W$, setting $s := \zeta_0^T|_{\a_1}(w_1)$, we have that for
Haar-almost any $g'\in G_s$ there is some $g'' \in G_s$ for which
\begin{multline*}
(\id_{G_s}\times
(L_{\rho_1(w_2)}\circ \Phi_{1,w_2}))\big((1,b(w_1,w_2)^{-1})\cdot M_{(w_1,w_2)}\cdot(1,b(w_1,w_2))\cdot (1,g')\big)\\
= (1,b(w_1,T_1|_{\a_2}(w_2))^{-1})\cdot
M_{(w_1,T_1|_{\a_2}(w_2))}\cdot (1,b(w_1,T_1|_{\a_2}(w_2)))\cdot
(1,g'').
\end{multline*}
We can re-write this condition in terms of the Mackey isomorphisms
to give
\begin{multline*}
\rho_1(w_2)\cdot(\Phi_{1,w_2}\circ L_{b(w_1,w_2)^{-1}}\circ\Psi_{(w_1,w_2)})(\bullet)\cdot \Phi_{1,w_2}(b(w_1,w_2)\cdot g')\\
=
b(w_1,T_1|_{\a_2}(w_2))^{-1}\cdot\Psi_{(w_1,T_1|_{\a_2}(w_2))}(\bullet)\cdot
b(w_1,T_1|_{\a_2}(w_2))\cdot g'',
\end{multline*}
and now if we write $\tilde{\Psi}_{\bullet} :=
\rm{Co}_{b(\bullet)^{-1}}\circ\Psi_\bullet$ and
$\tilde{\Phi}_{1,\bullet} :=
\rm{Co}_{\rho_1(\bullet)}\circ\Phi_{1,\bullet}$ this in turn becomes
\[(\tilde{\Phi}_{1,w_2}\circ \tilde{\Psi}_{(w_1,w_2)})(\bullet)\cdot\rho_1(w_2)\cdot \Phi_{1,w_2}(g')
= \tilde{\Psi}_{(w_1,T_1|_{\a_2}(w_2))}(\bullet)\cdot g'',\] (so we
have simply shifted all the `translation' parts of our affine
transformations over to the right).  Finally, this now clearly
requires that
\[\tilde{\Phi}_{1,w_2} =
\tilde{\Psi}_{(w_1,T_1|_{\a_2}(w_2))}\circ\tilde{\Psi}_{(w_1,w_2)}^{-1}\]
almost everywhere: the automorphisms graphed by the Mackey group
data have themselves become a coboundary for the automorphism-valued
cocycle $\tilde{\Phi}_{1,\bullet}$.

We will refer to the above as the `automorphism coboundary equation'
for the remainder of this proof. The condition that this hold for
almost every $(w_1,w_2)$ also gives nontrivial information on the
automorphisms $\tilde{\Psi}_{(w_1,w_2)}$ for different $w_1$, since
$w_1$ is absent from the left-hand side.  Since the extension
$\a_2|_\a:\bfW\to \bfW_2$ is relatively invariant for the
restrictions of $T_1$ and $T_3$, using
Proposition~\ref{prop:invar-meas-select} we can therefore choose a
$T_1$-invariant measurable selector $\eta:W_2 \to W_1$ so that
\[\tilde{\Phi}_{1,w_2} =
\tilde{\Psi}_{(\eta(w_2),T_1|_{\a_2}(w_2))}\circ\tilde{\Psi}_{(\eta(w_2),w_2)}^{-1},\]
holds almost surely and so witnesses that $\tilde{\Phi}_{1,w_2}$ is
a coboundary in $\rm{Aut}(G_s)$ for the transformation
$T_1|_{\a_2}:W_2\to W_2$.

Na\"\i vely we should now like to use the cocycle
$\tilde{\Psi}_{\eta(\bullet),\bullet}$ to make a fibrewise
automorphism recoordinatization of the extension
$\a_2|_{\xi_2}:\bfY_2 \to \bfW_2$ so that the first of our
automorphism-valued coboundary equations above gives a
simplification of $\tilde{\Phi}_{1,w_2}$. However, this idea runs
into difficulties because the new isomorphisms
$\tilde{\Psi}_\bullet$, unlike $\Psi_\bullet$, are not necessarily
$T_3$-invariant, and so applying them fibrewise may disrupt the
coordinatization of $T_3|_{\xi_2}$ as acting by rotations.

\quad\textbf{Step 2}\quad The best we can do at this stage is to
apply fibrewise the automorphisms $\Psi_{(\eta(w_2),w_2)}^{-1}$ to
our coordinatization of $\a_2|_{\xi_2}$.  This gives some
improvement: in the resulting new coordinatization of this
extension, our automorphism coboundary equation above now reads
\[\tilde{\Phi}_{1,w_2} = \rm{Co}_{b(\eta(w_2),T_1|_{\a_2}(w_2))}^{-1}\circ \rm{Co}_{b(\eta(w_2),w_2)}= \rm{Co}_{b(\eta(w_2),T_1|_{\a_2}(w_2))^{-1}\cdot b(\eta(w_2),w_2)}.\]

Recalling that $\tilde{\Phi}_{1,\bullet} =
\rm{Co}_{\rho_1(\bullet)}\circ \Phi_{1,\bullet}$ this unravels to
give
\begin{multline*}
\Phi_{1,w_2} =
\rm{Co}_{\rho_1(w_2)}^{-1}\circ\rm{Co}_{b(\eta(w_2),T_1|_{\a_2}(w_2))^{-1}\cdot b(\eta(w_2),w_2)}\\
= \rm{Co}_{\rho_1(w_2)^{-1}\cdot
b(\eta(w_2),T_1|_{\a_2}(w_2))^{-1}\cdot b(\eta(w_2),w_2)},
\end{multline*}
so we conclude, in particular, that the automorphism-valued cocycle
$\Phi_{1,\bullet}$ takes values in the compact subgroup of inner
automorphisms, and so we may represent it as
$\rm{Co}_{\theta(\bullet)}$ for some $T_3|_{\a_2}$-invariant section
$\theta:W_2 \to G_{\zeta_0^T|_{\a_2}(\bullet)}$. Writing out the
above automorphism coboundary equation in terms of $\theta$ it
becomes
\begin{multline*}
\rm{Co}_{\theta(w_2)} =
\rm{Co}_{\rho_1(w_2)}^{-1}\circ\rm{Co}_{b(\eta(w_2),T_1|_{\a_2}(w_2))^{-1}\cdot b(\eta(w_2),w_2)}\\
= \rm{Co}_{\rho_1(w_2)^{-1}\cdot
b(\eta(w_2),T_1|_{\a_2}(w_2))^{-1}\cdot b(\eta(w_2),w_2)},
\end{multline*}
and so it we now substitute into our original expression for
$T_1|_{\xi_2}$ we obtain
\begin{eqnarray*}
T_1|_{\xi_2} &=& T_1|_{\a_2}\ltimes
(L_{\rho_1(\bullet)}\circ\rm{Co}_{\theta(\bullet)})\\
&=& T_1|_{\a_2}\ltimes
(R_{\rho_1(\bullet)}\circ\rm{Co}_{\rho_1(\bullet)\cdot\theta(\bullet)})\\
&=& T_1|_{\a_2}\ltimes
(R_{\rho_1(\bullet)}\circ\rm{Co}_{b(\eta(\bullet),T_1|_{\a_2}(\bullet))^{-1}\cdot
b(\eta(\bullet),\bullet)}).
\end{eqnarray*}

It follows that if we now make a second fibrewise recoordinatization
of $\a_2|_{\xi_2}:\bfY_2 \to \bfW_2$, this time by rotating each
fibre copy of $G_s$ from the left by $b(\eta(w_2),w_2)$ (which
virtually depends only on $s = \zeta_0^T|_{\a_2}(w_2)$), we are left
with a resulting coordinatization of $T_1|_{\xi_2}$ in the desired
form of an opposite action:
\[T_1|_{\xi_2} = T_1|_{\a_2}\ltimes R_{\tau_1(\bullet)}\]
where
\[\tau_1(\bullet) :=
b(\eta(\bullet),\bullet)^{-1}\cdot
b(\eta(\bullet),T_1|_{\a_2}(\bullet))\cdot \rho_1(\bullet).\]

Of course, both of the above recoordinatizations can be repeated
analogously for the extension $\a_1|_{\xi_1}:\bfY_1\to \bfW_1$ to
put $T_2|_{\xi_1}$ into a similar right-multiplicative form.

It follows that in the coordinatizations of these extensions that we
have now obtained, $T_3|_{\xi_i}$ is still in the form of a
cocycle-section extension $T_3|_{\a_i}\ltimes \s_i$ (with a modified
cocycle-section $\s_i$) and $T_1|_{\xi_2}$ and $T_2|_{\xi_1}$ are in
the desired right-multiplicative form, as for part (1) of the
proposition.

\quad\textbf{Step 3}\quad We now `invert' the above implication to
discover what consequences these improved coordinatizations imply
for the data $\Psi_\bullet$ and $b_\bullet$.

Recall that our first automorphism cocycle equation held for almost
all $(w_1,w_2)$, before we chose the measurable selector $\eta$, and
so in our latest coordinatization this tells us that
\[\rm{Co}_{\tau_1(w_2)} =
\tilde{\Psi}_{(w_1,T_1|_{\a_2\circ\xi_2}(w_2))}\circ\tilde{\Psi}_{(w_1,w_2)}^{-1}.\]
(Note that our first fibrewise recoordinatization above by
automorphisms rendered the cocycle $\Psi_\bullet$ inner at
$(\a_2)_\#\mu$-almost all the points $(\eta(w_2),w_2)$, which depend
on our choice of measurable selector $\eta$, but we have not yet
seen that this cocycle is inner for almost all $(w_1,w_2)$ as a
result of this recoordinatization, hence our need to go back to the
above form of this equation for this stage of the argument.)

It follows that the class $\tilde{\Psi}_{(w_1,w_2)}\circ
\rm{Inn}(G_s) \in \rm{Out}(G_s)$ is invariant under the action of
$\id_{W_1}\times T_1|_{\a_2} = T_1|_\a$, and it follows similarly
that it is invariant under $T_2|_\a$. On the other hand, we have
$\tilde{\Psi}_{(w_1,w_2)}\circ \rm{Inn}(G_s) = \Psi_{(w_1,w_2)}\circ
\rm{Inn}(G_s)$ and this latter is clearly invariant under $T_3|_\a$,
since it arises from the Mackey group data. Therefore it is actually
$T|_\a$-invariant, and so since $\rm{Inn}(G_s) \unlhd \rm{Aut}(G_s)$
is compact, and so the resulting space of equivalence classes
$\rm{Out}(G_s)$ is smooth, it is almost surely equal to
$\Psi_s\circ\rm{Inn}(G)$ for some Borel map $\Psi_\bullet:Z_0^T\to
\rm{Aut}(G_\bullet)$. Therefore one last fibrewise automorphism
recoordinatization of $\a_2|_{\xi_2}:\bfY_2 \to \bfW_2$ by
$\Psi_{\zeta_0^T|_{\a_2}(\bullet)}$ (which still does not disrupt
any of the properties guaranteed previously, provided we replace
$\rho_1$ and $\tau_1$ with
$\Psi_{\zeta_0^T|_{\a_2}(\bullet)}(\rho_1)$ and
$\Psi_{\zeta_0^T|_{\a_2}(\bullet)}(\tau_1)$) now gives Mackey group
data of the form
\[M_\bullet \equiv \{(g,g_0(\bullet)gg_0(\bullet)^{-1}):\ g\in
G_{\zeta_0^T|_\a(\bullet)}\}\] for some $T_3|_{\a}$-invariant
section $g_0:W\to G_{\zeta_0^T|_\a(\bullet)}$, and now we can simply
adjust the Mackey section $b$ so that $g_0 \equiv 1$, and so the
Mackey group data can be taken to be the diagonal subgroup almost
everywhere.

\quad\textbf{Step 4}\quad Having removed all the nontrivial outer
automorphisms and adjusted the joining Mackey group data, our
automorphism coboundary equation has now simplified down to
\[\t{\Phi}_{1,w_2} = \rm{Co}_{\tau_1(\bullet)}\circ \rm{Co}_{\tau_1(\bullet)}^{-1} = \rm{id} = \rm{Co}_{b(w_1,T_1|_{\a_2}(w_2))^{-1}\cdot b(w_1,w_2)},\]
and hence we deduce that
$b(\bullet)\cdot\rm{C}(G_{\zeta_0^T|_\a(\bullet)})$, where
$\rm{C}(G_{\zeta_0^T|_\a(\bullet)})$ is the centre of
$G_{\zeta_0^T|_\a(\bullet)}$, is $T_1|_\a$-invariant, and similarly
that it is $T_2|_\a$-invariant. Making another measurable selection
and recoordinatizing each fibre of $\a_2|_{\xi_2}:\bfY_2 \to \bfW_2$
by a left-rotation by $b(\eta(\bullet),\bullet)$ therefore preserves
the structure of $T_1|_{\a}$ as an \emph{opposite} rotation (since
the resulting additional cocycle $b(w_1,T_1|_{\a_2}(w_2))^{-1}\cdot
b(w_1,w_2)$ acting on the left takes values in $\rm{C}(G_s)$, and so
may in fact be taken to act on either side); and after making this
recoordinatization we find that the Mackey section has also
trivialized.

\quad\textbf{Step 5}\quad Finally, let us look back at the relation
between $\s_1$ and $\s_2$ that is implied by the cocycle equation
satisfied by the Mackey data given by part (3) of
Theorem~\ref{thm:nonergMackey} in light of this newly-simplified
Mackey group and section: this now becomes simply that
\[\s_2(w_2) = \s_1(w_1)\]
$\a_\#\mu$-almost surely, and hence in this coordinatization it
follows that each $\s_i$ virtually depends only on
$\zeta_0^{T_1,T_2}|_{\a_i}$, or, equivalently, is
$T_{3-i}|_{\a_i}$-invariant.  Given this, the condition that
$T_1|_{\xi_2}$ and $T_2|_{\xi_2}$ commute simply reads that for
almost every $w_2 \in W_2$ we have
\[\s_2(w_2)\cdot g\cdot \tau_1(T_3|_{\a_2}(w_2)) = \s_2(T_1|_{\a_2}(w_2))\cdot g\cdot \tau_1(w_2) = \s_2(w_2)\cdot g\cdot \tau_1(w_2)\quad\forall g\in G_s,\]
and so we must also have that $\tau_1$ is $T_3|_{\a_2}$-invariant,
and similarly that $\tau_2$ is $T_3|_{\a_1}$-invariant. This
completes the proof. \qed

The recoordinatization of the preceding proposition leaves only one
detail remaining for the proof of
Proposition~\ref{prop:describe-invt-factor}.

\begin{cor}
In our homogeneous-space data coordinatizations of
$\a_i|_{\xi_i}:\bfY_i\to \bfW_i$ the core-free kernels
$K_{i,\bullet}$ are almost surely trivial.
\end{cor}

\textbf{Proof}\quad As remarked at the beginning of the preceding
proof, we can lift to the covering group extensions and make the
adjustments of Proposition~\ref{prop:final-recoord} there, and they
will then quotient back down to well-defined recoordinatizations of
the original extensions, because at each stage we have only applied
either fibrewise automorphism or fibrewise left-rotations.  From
these we have obtained expressions
\[T_1|_{\xi_2} = T_1|_{\a_2}\ltimes R_{\tau_1}\]
and similarly for $T_2|_{\xi_1}$ at the level of the covering group
extensions, and so for these to have well-defined quotient it is
necessary that for almost every $w_1$ all two-sided cosets $g\cdot
K_{2,w_2}\cdot \tau_1(w_2)$ for $g \in G_{\zeta_0^T|_{\a_2}(w_2)}$
actually be left-cosets of $K_{2,w_2}$. This, in turn, requires that
$\tau_1$ almost surely take values in the normalizer
$\rm{N}_{G_{\zeta_0^T|_{\a_2}(w_2)}}(K_{2,w_2})$, which is a closed
measurably-varying subgroup of $G_{\zeta_0^T|_{\a_2}(w_2)}$.

Now we recall that $T_1$ restricts to a \emph{relatively ergodic}
action on the extension $\a_2|_{\xi_2}:\bfY_2\to \bfW_2$ --- a
condition we have not exploited so far --- and so we must have
\[\rm{N}_{G_{\zeta_0^T|_{\a_2}(w_2)}}(K_{2,w_2})\cdot K_{2,w_2} =
G_{\zeta_0^T|_{\a_2}(w_2)}\] almost surely, for otherwise the
homogeneous space fibres of the extension $\a_2|_{\xi_2}:\bfY_2\to
\bfW_2$ would decompose into cosets of the closed subgroups
$\rm{N}_{G_{\zeta_0^T|_{\a_2}(w_2)}}(K_{2,w_2})\cdot K_{2,w_2}$ to
give additional nontrivial invariant sets under the restriction
$T_1|_{\a_2}$. However, since
\[\rm{N}_{G_{\zeta_0^T|_{\a_2}(w_2)}}(K_{2,w_2})\supseteq
K_{2,w_2},\] this requires in fact that
\[\rm{N}_{G_{\zeta_0^T|_{\a_2}(w_2)}}(K_{2,w_2}) =
G_{\zeta_0^T|_{\a_2}(w_2)}\] almost surely, and since $K_{2,w_2}$ is
core-free this is possible only if $K_{2,w_2} = \{1\}$ almost
surely. An exactly similar argument treats $K_{2,\bullet}$. \qed

As remarked previously, this completes the proof of
Proposition~\ref{prop:describe-invt-factor}. \qed

\subsection{Application to characteristic factors}\label{subs:charfactors}

We will finish this section by offering a second application of our
machinery (although in truth it is largely a corollary of the
above).

Since Furstenberg's ergodic theoretic proof of Szemer\'edi's Theorem
in~\cite{Fur77} and his extension with Katznelson of this result to
the multi-dimensional setting in~\cite{FurKat78}, considerable
interest has been attracted by the `non-conventional' ergodic
averages
\[\frac{1}{N}\sum_{n=1}^N\prod_{i=1}^df_i\circ T_i^n\]
associated to a commuting $d$-tuple of probability-preserving
transformations $T_1$, $T_2$, \ldots, $T_d:\bbZ\curvearrowright
(X,\mu)$, that emerge naturally in the course of those proofs. That
these averages converge in $L^2(\mu)$ as $N\to \infty$ in the case
$d=2$ was first shown by Conze and Lesigne in~\cite{ConLes84}, and
their result has since been extended in many
directions~\cite{ConLes88.1,ConLes88.2,Zha96,HosKra01,HosKra05,Zie05},
culminating in the first proof of the fully general case by Tao
in~\cite{Tao08(nonconv)}.  We refer the reader to those papers and
to~\cite{Aus--nonconv} for more thorough motivation and historical
discussion of this problem.

Conze and Lesigne's proof of convergence is comparatively soft,
using only quite weak structural information about the above
averages to show that they converge (in particular, using only the
structure of certain finite-rank modules rather than their concrete
coordinatizations). More recently, other convergence results for
nonconventional ergodic averages have been based on a similar but
more detailed analysis, resting on the notion of a `characteristic
tuple of factors'. A tuple of factors $\xi_i:\bfX\to\bfY_i$ is
\textbf{characteristic} if
\[\frac{1}{N}\sum_{n=1}^N\prod_{i=1}^df_i\circ T_i^n - \frac{1}{N}\sum_{n=1}^N\prod_{i=1}^d\sfE_\mu(f_i\,|\,\xi_i)\circ T_i^n \to 0\]
in $L^2(\mu)$ as $N\to\infty$ for any $f_1$, $f_2$, \ldots, $f_d \in
L^\infty(\mu)$. Starting with the Conze-Lesigne proof, most
convergence proofs in this area require at some stage the
identification of a characteristic tuple of factors (or a suitable
finitary analog of them in the case of Tao's proof) on which the
restricted actions of each $T_i$ take simplified forms, so that the
right-hand averages above can be analyzed to prove convergence more
easily.

A precise description of these characteristic factors in the special
case when $T_i = T^i$ for some fixed ergodic transformation $T$ has
now been given in terms of the special class of `pronilsystems' in
work of Host and Kra~\cite{HosKra05} (see also the subsequent
approach of Ziegler~\cite{Zie07}).  Frantzikinakis and Kra have
extended this description to more general commuting tuples subject
to some additional ergodicity assumption in~\cite{FraKra05}, but a
description for arbitrary tuples of commuting transformations,
without those ergodicity assumptions, seems to be more difficult.
Indeed, it may be that no comparably clean and useful description is
available in the general case. However, at least when $d=2$ a
reasonably simple coordinatization of a characteristic pair of
factors seems to have been folklore knowledge in ergodic theory for
some time, and in this subsection we will show how our theory
enables a careful proof of it.

\begin{thm}[Characteristic factors for double nonconventional
averages]\label{thm:charfactors2} Given a $\bbZ^2$-system $\bfX =
(X,\mu,T_1,T_2)$, let $\bfW_i$ be the target system of the joined
factor $\zeta_0^{T_i}\vee\zeta_0^{T_1 = T_2}$.  Then $\bfX$ admits a
characteristic pair of factors $\xi_i:\bfX\to\bfY_i$ that extend the
factors $\zeta_0^{T_i}\vee\zeta_0^{T_1 = T_2}$ and can be described
as follows: there are $T$-invariant compact group data $G_\bullet$,
a $T_1|_{\zeta_0^{T_1 = T_2}}$-ergodic cocycle $\s:\bfZ_0^{T_1 =
T_2} \to G$, and a pair of $T_i|_{\zeta_0^{T_{3-i}}}$-ergodic
cocycles $\tau_i:\bfZ_0^{T_{3-i}}\to G$ such that we can
coordinatize
\begin{center}
$\phantom{i}$\xymatrix{
(Y_1,(\xi_1)_\#\mu)\ar[dr]_{\a_1|_{\xi_1}}\ar@{<->}[rr]^-\cong &&
(W_1,(\a_1)_\#\mu)\ltimes (G_\bullet,m_{G_\bullet})\ar[dl]^{\rm{canonical}}\\
& (W_1,(\a_1)_\#\mu) }
\end{center}
with
\[T_1|_{\xi_1} \cong T_1|_{\a_1}\ltimes \s\circ\zeta_0^{T_1 = T_2}|_{\a_1}\quad\hbox{and}\quad T_2|_{\xi_1} \cong T_2|_{\a_1}\ltimes (L_{\s\circ\zeta_0^{T_1 = T_2}|_{\a_1}}\circ R_{\tau_2\circ\zeta_0^{T_1}|_{\a_1}}),\]
and similarly
\begin{center}
$\phantom{i}$\xymatrix{
(Y_2,(\xi_2)_\#\mu)\ar[dr]_{\a_2|_{\xi_2}}\ar@{<->}[rr]^-\cong &&
(W_2,(\a_2)_\#\mu)\ltimes
(G_\bullet,m_{G_\bullet})\ar[dl]^{\rm{canonical}}\\
& (W_2,(\a_2)_\#\mu) }
\end{center}
with
\[T_1|_{\xi_2} \cong T_1|_{\a_2}\ltimes (L_{\s\circ\zeta_0^{T_1 = T_2}|_{\a_2}}\circ R_{\tau_1\circ\zeta_0^{T_2}|_{\a_2}})\quad\hbox{and}\quad T_2|_{\xi_2} \cong T_2|_{\a_2}\ltimes \s\circ\zeta_0^{T_1 = T_2}|_{\a_2}.\]
\end{thm}

\textbf{Remarks}\quad\textbf{1.}\quad The form of the
coordinatizations given above with one action extended by a cocycle
and the other by an opposite cocycle, similarly to
Theorem~\ref{thm:triple-isotropy-2}, is a special feature of the
case of two commuting transformations. It would be possible to
replace it with a coordinatization by homogeneous space space data
(but not group data) in which both extensions are by cocycles acting
on fibres on the left, for example by enlarging $G_\bullet$ to
$G_\bullet\times G_\bullet$, quotienting by the diagonal subgroup
$\{(g,g):\ g\in G_\bullet\}$ and having $\rho$ rotate
$G_\bullet\times G_\bullet$ only in the first coordinate and $\s_1$,
$\s_2$ only in the second. It is presumably this more canonical but
more fiddly representation, if any, that would admit generalization
to larger numbers of commuting transformations.

\quad\textbf{2.}\quad The above result describes the possible
structures of the two characteristic factors individually, but some
opacity remains as to how they can be joined inside $\bfX$. While we
suspect that the methods of the present section can be brought to
bear on this question also, we will not pursue this analysis in
detail here. \fin

The first steps of our analysis, which are essentially contained in
Conze and Lesigne~\cite{ConLes84} (as well as many subsequent
papers; see, for example, Furstenberg and Weiss~\cite{FurWei96} for
a nice treatment of this stage of the proof), give control over the
asymptotic behaviour of our averages in terms of a certain two-fold
self-joining of $\bfX$. We will then complete the proof essentially
by re-applying Proposition~\ref{prop:describe-invt-factor} to
certain factors of that self-joining.

We observe from the mean ergodic theorem that
\begin{multline*}
\int_X\frac{1}{N}\sum_{n=1}^N(f_1\circ T_1^n)(f_2\circ T_2^n)\,\d\mu
= \int_Xf_1\cdot \frac{1}{N}\sum_{n=1}^N(f_2\circ
(T_2T_1^{-1})^n)\,\d\mu\\ \to \int_X
f_1\cdot\sfE_\mu(f_1\,|\,\zeta_0^{T_1 = T_2})\,\d\mu = \int_{X^2}
f_1\otimes f_2\,\d\mu^{\rm{F}}
\end{multline*}
where $\mu^{\rm{F}}$ is the \textbf{Furstenberg self-joining}, which
in this case equals $\mu\otimes_{\zeta_0^{T_1 = T_2}}\mu$ (it has a
much more complicated structure for larger numbers of commuting
transformations which is not yet well understood;
see~\cite{Aus--nonconv,Aus--newmultiSzem} for further discussion of
this matter).  It is easy to check that $\mu^{\rm{F}}$ is invariant
under the lifted transformations $T_1^{\times 2}$ and $T_2^{\times
2}$, and also under the \textbf{diagonal transformation} $\vec{T} :=
T_1\times T_2$.

This self-joining now helps control our averages through the
following consequence of the van der Corput estimate (for which see,
for example, Bergelson~\cite{Ber87}).

\begin{lem}
If the pair of factors $\xi_i:\bfX\to\bfY_i$ is such that
$\xi_i\succsim \zeta_0^{T_1 = T_2}$ and $\zeta_0^{\vec{T}} \precsim
\xi_1\times \xi_2$ then this pair is characteristic.
\end{lem}

\textbf{Proof}\quad This follows from a routine application of the
van der Corput estimate. Suppose that $f_1,f_2 \in L^\infty(\mu)$;
clearly by symmetry and iterating our argument, it suffices to prove
that
\[\frac{1}{N}\sum_{n=1}^N(f_1\circ T_1^n)(f_2\circ T_2^n) - \frac{1}{N}\sum_{n=1}^N(\sfE_\mu(f_1\,|\,\xi_1)\circ T_1^n)(f_2\circ T_2^n)\to 0\]
in $L^2(\mu)$ as $N\to\infty$, and hence (taking the difference of
the two sides above) that
\[\frac{1}{N}\sum_{n=1}^N(f_1\circ T_1^n)(f_2\circ T_2^n) \to 0\]
in $L^2(\mu)$ if $\sfE_\mu(f_1\,|\,\xi_1) = 0$.

Letting $F_n := (f_1\circ T_1^n)(f_2\circ T_2^n)$, by the van der
Corput estimate this will follow if we show that
\[\frac{1}{M}\sum_{h=1}^M\frac{1}{N}\sum_{n=1}^N\langle F_n,F_{n+h}\rangle \to 0\]
as $N \to \infty$ and then $M\to\infty$. Now we simply compute
\begin{eqnarray*}
\frac{1}{N}\sum_{n=1}^N\langle F_n,F_{n+h}\rangle &=&
\frac{1}{N}\sum_{n=1}^N\int_X(f_1\circ T_1^n)(f_2\circ
T_2^n)(f_1\circ T_1^{n+h})(f_2\circ T_2^{n+h})\,\d\mu\\ &=& \int_X
\frac{1}{N}\sum_{n=1}^N((f_1\cdot (f_1\circ T_1^h))\circ T_1^n)\cdot
((f_2\cdot(f_2\circ T_2^h))\circ
T_2^n)\,\d\mu\\
&\to& \int_{X^2} (f_1\otimes f_2)((f_1\circ T_1^h)\otimes (f_2\circ
T_2^h))\,\d\mu^{\rm{F}},
\end{eqnarray*}
so that if we now average also in $h$ this converges by the mean
ergodic theorem to
\[\int_{X^2}(f_1\circ\pi_1)\cdot(f_2\circ\pi_2)\cdot g\,\d\mu^{\rm{F}}\]
for some $\vec{T}$-invariant function $g$. Finally, this last
integral is zero if
$\sfE_{\mu^{\rm{F}}}(f_1\circ\pi_1\,|\,\pi_2\vee\zeta_0^{\vec{T}}) =
0$, and this follows from our assumptions and the relative
independence of $\pi_1$ and $\pi_2$ over $\zeta_0^{T_1 =
T_2}\circ\pi_1$ under $\mu^{\rm{F}}$. \qed

Since on the other hand Theorem~\ref{thm:rel-ind-joinings} tells us
that $\zeta_0^{\vec{T}} \precsim
(\zeta_{1/\zeta_0^{T_1=T_2}}^{T_1}\circ\pi_1)\vee(\zeta_{1/\zeta_0^{T_1=T_2}}^{T_2}\circ\pi_2)$,
we can deduce the following at once.

\begin{cor}[Reduction to isometric extensions of isotropy factors]\label{prop:CL-reduction-1}
There is a characteristic pair of factors satisfying $\xi_i \precsim
\zeta_{1/\zeta_0^{T_1 = T_2}}^{T_i}$ for $i=1,2$. \qed
\end{cor}

We can now present the factors introduced above as another instance
of the situation described before
Proposition~\ref{prop:describe-invt-factor}, and have that
proposition do the heavy lifting we need again here.

To see this, we define a system of \emph{three} commuting
transformations on the Furstenberg self-joining. Let $\bfX^{\rm{F}}$
be the $\bbZ^3$-system $(X^2,\mu^{\rm{F}},S_1,S_2,S_3)$ obtained by
setting $S_1 := T_1^{\times 2}\vec{T}^{-1} = \id_X\times
(T_1T_2^{-1})$, $S_2 := T_2^{\times 2}\vec{T}^{-1} =
(T_2T_1^{-1})\times \id_X$ and $S_3 := \vec{T}$. We observe directly
from the definition of $\mu^{\rm{F}}$ that for $i=1,2$ the
coordinate projection $\pi_i:X^2\to X$ is equivalent to
$\zeta_0^{S_i}$.

Now let $\bfW_i$ be $(\zeta_0^{T_1 = T_2}\vee \zeta_0^{T_i})(\bfX)$
for $i=1,2$ and let $\bfZ_i$ be the target of the maximal
subextension of $\zeta_0^{T_1 = T_2}\vee
\zeta_0^{T_i}:\bfX\to\bfW_i$ that is isometric for the restriction
of $T_i$. Let
\[\a_i:= (\zeta_0^{T_1 = T_2}\vee
\zeta_0^{T_i})\circ\pi_i =
\zeta_0^{S_1,S_2}\vee\zeta_0^{S_i,S_3}:\bfX^\rm{F} \to \bfW_i\] and
\[\zeta_i:= \zeta_{1/\a_i}^{T_i}\circ\pi_i:\bfX^{\rm{F}} \to \bfZ_i,\]
and let $\zeta:\bfX^{\rm{F}}\to \bfZ$ and $\a:\bfX^{\rm{F}}\to \bfW$
be the joinings $\zeta_1 \vee \zeta_2$ and $\a_1\vee \a_2$
respectively.

It is now routine to check from the basic results above that these
data satisfy the same conditions as were needed for
Proposition~\ref{prop:describe-invt-factor} with $\bfX^{\rm{F}}$ in
place of $\bfX$ and $S_i$ in place of $T_i$: these factors are once
again arranged as in the commutative diagram
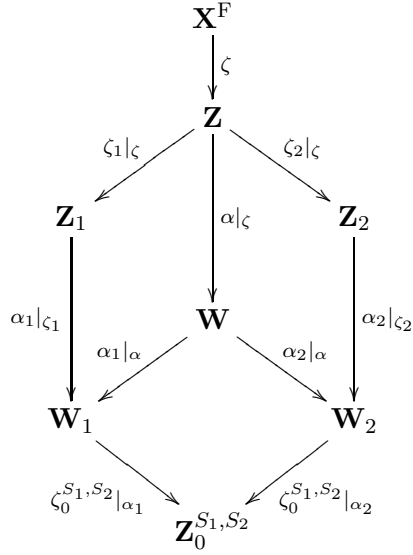
\begin{center}
$\phantom{i}$\xymatrix{ & \bfX^{\rm{F}}\ar[d]^{\zeta}\\
& \bfZ\ar[dl]_{\zeta_1|_{\zeta}}\ar[dd]^{\a|_{\zeta}}\ar[dr]^{\zeta_2|_{\zeta}}\\
\bfZ_1\ar[dd]_{\a_1|_{\zeta_1}} & & \bfZ_2\ar[dd]^{\a_2|_{\zeta_2}}\\
& \bfW\ar[dl]_{\a_1|_{\a}}\ar[dr]^{\a_2|_{\a}}\\
\bfW_1\ar[dr]_{\zeta_0^{S_1,S_2}|_{\a_1}} & & \bfW_2\ar[dl]^{\zeta_0^{S_1,S_2}|_{\a_2}}\\
& \bfZ_0^{S_1,S_2}}
\end{center}
where we observe easily from the structure of $\mu^{\rm{F}} =
\mu\otimes_{\zeta_0^{T_1 = T_2}}\mu$ that $\zeta_0^{S_1 = S_2}
\simeq \zeta_0^{T_1 = T_2}\circ\pi_1 \simeq \zeta_0^{T_1 =
T_2}\circ\pi_2$.  In addition, the transformations $S_i$ enjoy the
following properties:
\begin{itemize}
\item $S_i$ restricts to the identity on $\bfZ_i$ and the factors
beneath it, while acting relatively ergodically on the extension
$\zeta_0^{S_1,S_2}|_{\zeta_{3-i}}:\bfZ_{3-i}\to \bfZ_0^{S_1,S_2}$,
for $i=1,2$;
\item the extensions $\zeta_0^{S_1,S_2}|_{\a_i}:\bfW_i\to \bfZ_0^{S_1,S_2}$ are relatively
invariant for the restriction of $S_3$, and the extensions
$\a_i|_{\zeta_i}:\bfZ_i \to \bfW_i$ are relatively ergodic and
isometric for the restriction of $S_3$.
\end{itemize}

We can therefore apply Proposition~\ref{prop:describe-invt-factor}
to these systems and maps to deduce the following.

\begin{prop}
There are intermediate factors $\bfZ_i
\stackrel{\xi_i|_{\zeta_i}}{\longrightarrow} \bfY_i
\stackrel{\a_i|_{\xi_i}}{\longrightarrow} \bfW_i$ factorizing
$\a_i|_{\zeta_i}$ such that there are $S$-invariant compact group
data $G_\bullet$ and cocycle-sections
\begin{itemize}
\item[] $\s:Z_0^{S_1,S_2}\to G_\bullet$ that is
$S_3|_{\zeta_0^{S_1,S_2}}$-relatively ergodic,
\item[] $\tau_1:Z_0^{S_2,S_3}\to G_\bullet$ that is
$S_1|_{\zeta_0^{S_2,S_3}}$-relatively ergodic and
\item[] $\tau_2:Z_0^{S_2,S_3}\to G_\bullet$ that is
$S_2|_{\zeta_0^{S_1,S_3}}$-relatively ergodic
\end{itemize}
so that we can coordinatize the actions of the transformations $S$
as
\begin{center}
$\phantom{i}$\xymatrix{
\bfY_1\ar[dr]_{\a_1|_{\xi_1}}\ar@{<->}[rr]^-\cong &
& \bfW_1\ltimes (G_\bullet,m_{G_\bullet},1,(\tau_2\circ\zeta_0^{S_1,S_3}|_{\a_1})^{\rm{op}},\s\circ\zeta_0^{S_1,S_2}|_{\a_1})\ar[dl]^{\rm{canonical}}\\
& \bfW_1, }
\end{center}
and
\begin{center}
$\phantom{i}$\xymatrix{
\bfY_2\ar[dr]_{\a_2|_{\xi_2}}\ar@{<->}[rr]^-\cong &
& \bfW_2\ltimes (G_\bullet,m_{G_\bullet},(\tau_1\circ\zeta_0^{S_2,S_3}|_{\a_2})^{\rm{op}},1,\s\circ\zeta_0^{S_1,S_2}|_{\a_2})\ar[dl]^{\rm{canonical}}\\
& \bfW_2, }
\end{center}
and such that the $(S_3 = \vec{T})$-invariant factor of
$\bfX^{\rm{F}}$ is contained in $\xi_1\vee \xi_2$. \qed
\end{prop}

\textbf{Proof of Theorem~\ref{thm:charfactors2}}\quad Again this
follows simply by unpacking the notation of the above result: the
tower of factors $\a_i|_{\xi_i}:\bfY_i\to \bfW_i$ of $\bfX^{\rm{F}}$
are all actually contained within the $i^{\rm{th}}$ coordinate
projection $\pi_i:\bfX^{\rm{F}}\to \bfX_i$, and by definition we
have $T_i = (S_iS_3)|_{\pi_i}$ for $i=1,2$, and therefore the above
coordinatization of the action of $S$ restricted to the tower of
factors $\a_1|_{\xi_1}:\bfY_1 \to\bfW_1$ converts into a
coordinatization description of $T_1|_{\xi_1}$ as
\[(w,g)\mapsto (T_1|_{\a_1}(w),\s(\zeta_0^{T_1 = T_2}(w))\cdot g)\]
and of $T_2|_{\xi_1}$ as
\[(w,g)\mapsto (T_1|_{\a_1}(w),\s(\zeta_0^{T_1 = T_2}(w))\cdot g\cdot \tau_2(\zeta_0^{T_1}(w))),\]
(since $\zeta_0^{S_1,S_3}\simeq \zeta_0^{T_1}\circ\pi_1$) and
similarly for $T_1|_{\xi_2}$ and $T_2|_{\xi_2}$. This completes the
proof. \qed

\textbf{Remarks}\quad\textbf{1.}\quad In fact, it is relatively easy
to see by checking functions $f_1$, $f_2$ that are constructed from
measurable selections of representative functions on the compact
fibre groups $G_\bullet$ that the characteristic pair of factors
$\xi_i$ that we have now isolated is minimal, in that any other
characteristic pair $\xi'_1$, $\xi_2'$ satisfies $\xi_i \precsim
\xi_i'$.

\quad\textbf{2.}\quad The results of the preceding subsection also
give a precise picture of the $\vec{T}$-invariant factor of
$\bfX^\rm{F}$ in terms of a diagonal Mackey group and trivial Mackey
section for the joining of the above coordinatizations of $\bfY_1$
and $\bfY_2$ inside $\bfX^\rm{F}$; we omit these details here. \fin

\section{Further questions}\label{sec:furtherques}

This paper leaves open the obvious question of how to generalize the
analysis of Section~\ref{sec:app} to describe in similar detail
\begin{itemize}
\item the possible joint distributions among a larger collection of
isotropy factors $\zeta_0^{T^{\ \uhr\G_i}}:\bfX\to\bfZ_0^{T^{\
\uhr\G_i}}$ for $\G_1$, $\G_2$, \ldots, $\G_d \leq \bbZ^d$;
\item the possible structures of characteristic factors (and,
relatedly, Furstenberg self-joinings) for larger commuting tuples of
transformations (or commuting actions of some other fixed group).
\end{itemize}

On the one hand, it seems likely that the machinery of extensions by
measurably-varying compact homogeneous spaces will be quite
essential to any further developments in this area.  On the other, I
suspect that even the next cases to consider in the natural
hierarchy (joint distributions of four isotropy factors, or
characteristic factors for three commuting transformations) become
much more complicated, and it may be in general too much to ask for
the kind of precision that we obtained in
Theorems~\ref{thm:threefoldisotropy} and~\ref{thm:charfactors}.

There is an alternative viewpoint on questions such as these that
may be more tractable.  Instead of asking about exact joint
distributions or characteristic factors for an initially-given
system, if we allow ourselves the freedom to pass to any extension
of that system matters sometimes improve considerably.  Indeed, the
new proofs of convergence for linear nonconventional averages
in~\cite{Aus--nonconv} and of Furstenberg and Katznelson's
associated multidimensional multiple recurrence theorem
in~\cite{Aus--newmultiSzem} both relied on procedures for passing
from an initially-given system to some extension in which the
relevant characteristic factors and their joint distributions could
be described much more simply.

The constructions of those papers were abstract enough to work
without any of the machinery of homogeneous-space-data extensions.
However, for further applications of this idea, in particular to the
problem of convergence of related `polynomial nonconventional
averages' such as
\[\frac{1}{N}\sum_{n=1}^N(f_1\circ T_1^{n^2})(f_1\circ T_1^{n^2}T_2^n)\]
(discussed, for example, by Bergelson and Leibman
in~\cite{BerLei99}), it seems likely that these more delicate tools
will be necessary. In the forthcoming
works~\cite{Aus--lindeppleasant1,Aus--lindeppleasant2,Aus--lindeppleasant3}
we will make such an analysis allowing ourselves to pass to
extensions, focusing on what improved characteristic factors can be
found while retaining some given algebraic relations among the
transformations involved, and will then use this to prove
convergence of the above polynomial averages in $L^2(\mu)$ as
$N\to\infty$.

In addition to these quite specialized applications, let us also
mention that there seem to be further issues on the general
behaviour of extensions by homogeneous space data to be explored.
For example, in~\cite{AusLem--Rokhlincocycs} it is shown that some
of the machinery of Furstenberg and Zimmer concerning finite-rank
modules of an extension can be extended to the setting in which the
extension is relatively finite measure-preserving, but the base,
while ergodic, is only assumed to be non-singular (that is, its
measure is only quasi-invariant). In that paper this machinery is
needed for the proof of a result about the lifting of the
`multiplier property' through certain kinds of extension, and this
will not require that these general results on finite-rank modules
be pushed very far. However, it might be interesting to examine
whether that development can be easily recovered without the
assumption of ergodicity of the base, using a version of the
formalism of the present paper.

\parskip 0pt

\bibliographystyle{abbrv}
\bibliography{bibfile}

\begin{thebibliography}{10}

\bibitem{Arv76}
W.~Arveson.
\newblock {\em An Invitation to C$^\ast$-Algebras}.
\newblock Springer, 1976.

\bibitem{AusGreHah63}
L.~Auslander, L.~Green, and F.~Hahn.
\newblock {\em Flows on homogeneous spaces}.
\newblock With the assistance of L. Markus and W. Massey, and an appendix by L.
  Greenberg. Annals of Mathematics Studies, No. 53. Princeton University Press,
  Princeton, N.J., 1963.

\bibitem{Aus--newmultiSzem}
T.~Austin.
\newblock Deducing the multidimensional {S}zemer\'edi {T}heorem from an
  infinitary removal lemma.
\newblock To appear, \emph{{J}. d'{A}nalyse {M}ath.}, 2008.

\bibitem{Aus--nonconv}
T.~Austin.
\newblock On the norm convergence of nonconventional ergodic averages.
\newblock To appear, \emph{{E}rgodic {T}heory {D}ynam. {S}ystems}, 2008.

\bibitem{Aus--lindeppleasant1}
T.~Austin.
\newblock Pleasant extensions retaining algebraic structure, {I}.
\newblock Preprint, available online at \verb|arXiv.org|: 0905.0518, 2009.

\bibitem{Aus--lindeppleasant2}
T.~Austin.
\newblock Pleasant extensions retaining algebraic structure, {II}.
\newblock Preprint, available online at \verb|arXiv.org|: 0910.0907, 2009.

\bibitem{Aus--lindeppleasant3}
T.~Austin.
\newblock Pleasant extensions retaining algebraic structure, {III}.
\newblock Preprint, available online at \verb|arXiv.org|: 0910.0909, 2009.

\bibitem{AusLem--Rokhlincocycs}
T.~Austin and M.~Lema\'nczyk.
\newblock Relatively finite measure-preserving extensions and lifting
  multipliers by {R}okhlin cocycles.
\newblock To appear, \emph{{J}. {F}ixed {P}oint {T}heory {A}ppl.}

\bibitem{Ber87}
V.~Bergelson.
\newblock Weakly mixing {PET}.
\newblock {\em Ergodic Theory Dynam. Systems}, 7(3):337--349, 1987.

\bibitem{BerLei99}
V.~Bergelson and A.~Leibman.
\newblock Set-polynomials and polynomial extension of the {H}ales-{J}ewett
  theorem.
\newblock {\em Ann. of Math. (2)}, 150(1):33--75, 1999.

\bibitem{BerTaoZie09}
V.~Bergelson, T.~Tao, and T.~Ziegler.
\newblock An inverse theorem for the uniformity seminorms associated with the
  action of $\mathbb{F}_p^\infty$.
\newblock Preprint, available online at \verb|arXiv.org|: 0901.2602, 2009.

\bibitem{Brotom85}
T.~Br\"ocker and T.~tom Dieck.
\newblock {\em Representations of Compact Lie Groups}.
\newblock Springer, 1985.

\bibitem{ConLes84}
J.-P. Conze and E.~Lesigne.
\newblock Th\'eor\`emes ergodiques pour des mesures diagonales.
\newblock {\em Bull. Soc. Math. France}, 112(2):143--175, 1984.

\bibitem{ConLes88.1}
J.-P. Conze and E.~Lesigne.
\newblock Sur un th\'eor\`eme ergodique pour des mesures diagonales.
\newblock In {\em Probabilit\'es}, volume 1987 of {\em Publ. Inst. Rech. Math.
  Rennes}, pages 1--31. Univ. Rennes I, Rennes, 1988.

\bibitem{ConLes88.2}
J.-P. Conze and E.~Lesigne.
\newblock Sur un th\'eor\`eme ergodique pour des mesures diagonales.
\newblock {\em C. R. Acad. Sci. Paris S\'er. I Math.}, 306(12):491--493, 1988.

\bibitem{ConRau07}
J.-P. Conze and A.~Raugy.
\newblock On the ergodic decomposition for a cocycle.
\newblock To appear, \emph{Colloq. Math.}

\bibitem{Dow06}
T.~Downarowicz.
\newblock Minimal models for noninvertible and not uniquely ergodic systems.
\newblock {\em Israel J. Math.}, 156:93--110, 2006.

\bibitem{Eng89}
R.~Engelking.
\newblock {\em General topology}, volume~6 of {\em Sigma Series in Pure
  Mathematics}.
\newblock Heldermann Verlag, Berlin, second edition, 1989.
\newblock Translated from the Polish by the author.

\bibitem{FisMorWhy04}
D.~Fisher, D.~W. Morris, and K.~Whyte.
\newblock Nonergodic actions, cocycles and superrigidity.
\newblock {\em New York J. Math.}, 10:249--269 (electronic), 2004.

\bibitem{FraKra05}
N.~Frantzikinakis and B.~Kra.
\newblock Convergence of multiple ergodic averages for some commuting
  transformations.
\newblock {\em Ergodic Theory Dynam. Systems}, 25(3):799--809, 2005.

\bibitem{FreVol3}
D.~H. Fremlin.
\newblock {\em Measure Theory, Volume 3: Measure Algebras}.
\newblock Torres Fremlin, Colchester, 2004.

\bibitem{FreVol4}
D.~H. Fremlin.
\newblock {\em Measure Theory, Volume 4: Topological Measure Theory}.
\newblock Torres Fremlin, Colchester, 2005.

\bibitem{Fur77}
H.~Furstenberg.
\newblock Ergodic behaviour of diagonal measures and a theorem of {S}zemer\'edi
  on arithmetic progressions.
\newblock {\em J. d'Analyse Math.}, 31:204--256, 1977.

\bibitem{Fur81}
H.~Furstenberg.
\newblock {\em Recurrence in Ergodic Theory and Combinatorial Number Theory}.
\newblock Princeton University Press, Princeton, 1981.

\bibitem{FurKat78}
H.~Furstenberg and Y.~Katznelson.
\newblock An ergodic {S}zemer\'edi {T}heorem for commuting transformations.
\newblock {\em J. d'Analyse Math.}, 34:275--291, 1978.

\bibitem{FurKat85}
H.~Furstenberg and Y.~Katznelson.
\newblock An ergodic {S}zemer\'edi theorem for {I}{P}-systems and combinatorial
  theory.
\newblock {\em J. d'Analyse Math.}, 45:117--168, 1985.

\bibitem{FurKat91}
H.~Furstenberg and Y.~Katznelson.
\newblock A {D}ensity {V}ersion of the {H}ales-{J}ewett {T}heorem.
\newblock {\em J. d'Analyse Math.}, 57:64--119, 1991.

\bibitem{FurWei96}
H.~Furstenberg and B.~Weiss.
\newblock A mean ergodic theorem for
  $\frac{1}{N}\sum_{n=1}^{N}f({T}^nx)g({T}^{n^2}x)$.
\newblock In V.~Bergleson, A.~March, and J.~Rosenblatt, editors, {\em
  Convergence in Ergodic Theory and Probability}, pages 193--227. De Gruyter,
  Berlin, 1996.

\bibitem{Gla03}
E.~Glasner.
\newblock {\em Ergodic {T}heory via {J}oinings}.
\newblock American {M}athematical {S}ociety, {P}rovidence, 2003.

\bibitem{Gui72}
A.~Guichardet.
\newblock {\em Symmetric Hilbert Spaces and Related Topics}.
\newblock Springer, 1972.

\bibitem{HosKra01}
B.~Host and B.~Kra.
\newblock Convergence of {C}onze-{L}esigne averages.
\newblock {\em Ergodic Theory Dynam. Systems}, 21(2):493--509, 2001.

\bibitem{HosKra05}
B.~Host and B.~Kra.
\newblock Nonconventional ergodic averages and nilmanifolds.
\newblock {\em Ann. Math.}, 161(1):397--488, 2005.

\bibitem{Kal02}
O.~Kallenberg.
\newblock {\em Foundations of modern probability}.
\newblock Probability and its Applications (New York). Springer-Verlag, New
  York, second edition, 2002.

\bibitem{Mac66}
G.~W. Mackey.
\newblock Ergodic theory and virtual groups.
\newblock {\em Math. Ann.}, 166:187--207, 1966.

\bibitem{Mac76}
G.~W. Mackey.
\newblock {\em The theory of unitary group representations}.
\newblock University of Chicago Press, Chicago, Ill., 1976.
\newblock Based on notes by James M. G. Fell and David B. Lowdenslager of
  lectures given at the University of Chicago, Chicago, Ill., 1955, Chicago
  Lectures in Mathematics.

\bibitem{Mei90}
D.~Meiri.
\newblock Generalized correlation sequences.
\newblock Master's thesis, Tel Aviv University; available online at\newline
  \verb|http://taalul.com/David/Math/ma.pdf|, 1990.

\bibitem{Men91}
M.~K. Mentzen.
\newblock Ergodic properties of group extensions of dynamical systems with
  discrete spectra.
\newblock {\em Studia Math.}, 101(1):19--31, 1991.

\bibitem{New79}
D.~Newton.
\newblock On canonical factors of ergodic dynamical systems.
\newblock {\em J. London Math. Soc. (2)}, 19(1):129--136, 1979.

\bibitem{Rud93}
D.~J. Rudolph.
\newblock Eigenfunctions of {$T\times S$} and the {C}onze-{L}esigne algebra.
\newblock In {\em Ergodic theory and its connections with harmonic analysis
  ({A}lexandria, 1993)}, volume 205 of {\em London Math. Soc. Lecture Note
  Ser.}, pages 369--432. Cambridge Univ. Press, Cambridge, 1995.

\bibitem{Sch95}
K.~Schmidt.
\newblock {\em Dynamical systems of algebraic origin}, volume 128 of {\em
  Progress in Mathematics}.
\newblock Birkh\"auser Verlag, Basel, 1995.

\bibitem{Sch94}
R.~Schmidt.
\newblock {\em Subgroup lattices of groups}, volume~14 of {\em de Gruyter
  Expositions in Mathematics}.
\newblock Walter de Gruyter \& Co., Berlin, 1994.

\bibitem{Sta00}
A.~N. Starkov.
\newblock {\em Dynamical systems on homogeneous spaces}, volume 190 of {\em
  Translations of Mathematical Monographs}.
\newblock American Mathematical Society, Providence, RI, 2000.
\newblock Translated from the 1999 Russian original by the author.

\bibitem{Tao08(nonconv)}
T.~Tao.
\newblock Norm convergence of multiple ergodic averages for commuting
  transformations.
\newblock {\em Ergodic Theory and Dynamical Systems}, 28:657--688, 2008.

\bibitem{Tow07}
H.~P. Towsner.
\newblock Convergence of {D}iagonal {E}rgodic {A}verages.
\newblock Preprint, available online at \verb|arXiv.org|: 0711.1180, 2007.

\bibitem{Vee82}
W.~A. Veech.
\newblock A criterion for a process to be prime.
\newblock {\em Monatsh. Math.}, 94(4):335--341, 1982.

\bibitem{Zha96}
Q.~Zhang.
\newblock On convergence of the averages {$(1/N)\sum\sp N\sb {n=1}f\sb 1(R\sp
  nx)f\sb 2(S\sp nx)f\sb 3(T\sp nx)$}.
\newblock {\em Monatsh. Math.}, 122(3):275--300, 1996.

\bibitem{Zie05}
T.~Ziegler.
\newblock A non-conventional ergodic theorem for a nilsystem.
\newblock {\em Ergodic Theory Dynam. Systems}, 25(4):1357--1370, 2005.

\bibitem{Zie07}
T.~Ziegler.
\newblock Universal characteristic factors and {F}urstenberg averages.
\newblock {\em J. Amer. Math. Soc.}, 20(1):53--97 (electronic), 2007.

\bibitem{Zim76.2}
R.~J. Zimmer.
\newblock Ergodic actions with generalized discrete spectrum.
\newblock {\em Illinois J. Math.}, 20(4):555--588, 1976.

\bibitem{Zim76.1}
R.~J. Zimmer.
\newblock Extensions of ergodic group actions.
\newblock {\em Illinois J. Math.}, 20(3):373--409, 1976.

\bibitem{Zim84}
R.~J. Zimmer.
\newblock {\em Ergodic theory and semisimple groups}, volume~81 of {\em
  Monographs in Mathematics}.
\newblock Birkh\"auser Verlag, Basel, 1984.

\end{thebibliography}

\vspace{10pt}

\small{\textsc{Department of Mathematics, University of California,
Los Angeles CA 90095-1555, USA}}

\vspace{7pt}

\small{Email: \verb|timaustin@math.ucla.edu|}

\vspace{7pt}

\small{URL: \verb|http://www.math.ucla.edu/~timaustin|}

\end{document}